\documentstyle[amssymb,righttag]{amsart}

\advance\textheight by \topskip

\newtheorem{teo}{\sc Theorem}[section]

\newtheorem{cor}{\sc Corollary}[section]
\newtheorem{lemma}{\sc Lemma}[section]
\newtheorem{pro}{\sc Proposition}[section]

\theoremstyle{definition}
\newtheorem{defi}{\sc Definition}[section]

\theoremstyle{remark}
\newtheorem{rem}[teo]{\sc Remark}

\newtheorem{notation}{\sc Notation}

\newcommand{\bte}{\begin{teo}}
\newcommand{\ete}{\end{teo}}
\newcommand{\bc}{\begin{cor}}
\newcommand{\ec}{\end{cor}}
\newcommand{\bp}{\begin{pro}}
\newcommand{\ep}{\end{pro}}
\newcommand{\bl}{\begin{lemma}}
\newcommand{\el}{\end{lemma}}
\newcommand{\bd}{\begin{defi}}
\newcommand{\ed}{\end{defi}}
\newcommand{\bno}{\begin{notation}}
\newcommand{\eno}{\end{notation}}
\newcommand{\fp}{\hfill $\Box$}

\newcommand{\bca}{\begin{cases}}
\newcommand{\eca}{\end{cases}}
\newcommand{\la}{\langle}
\newcommand{\ra}{\rangle}

\newcommand{\bq}{\begin{equation}}
\newcommand{\eq}{\end{equation}}
\newcommand{\btabu}{\begin{table}}
\newcommand{\etabu}{\end{table}}
\newcommand{\bt}{\begin{tabular}}
\newcommand{\et}{\end{tabular}}
\newcommand{\ba}{\begin{array}}
\newcommand{\ea}{\end{array}}
\newcommand{\br}{\begin{eqnarray}}
\newcommand{\er}{\end{eqnarray}}
\newcommand{\brn}{\begin{eqnarray*}}
\newcommand{\ern}{\end{eqnarray*}}
\newcommand{\benu}{\begin{enumerate}}
\newcommand{\eenu}{\end{enumerate}}
\newcommand{\bite}{\begin{itemize}}
\newcommand{\eite}{\end{itemize}}

\newcommand{\supp}{\operatorname{supp }}

\title[Nikishin systems are perfect] {Nikishin systems are perfect}

\thanks{Dedicated to the memory of the outstanding russian mathematician E.M. Nikishin who died on December 17, 1987, at the early age of 42. See \cite{obi} for a brief account of his results and list of publications.}

\thanks{The work of both authors was
supported by Ministerio de Ciencia y Tecnolog\'{\i}a under grants
MTM2006-13000-C03-02  and MTM2009-12740-C03-01.}

\author[Fidalgo]{U. Fidalgo Prieto}

\address[Fidalgo]{Departamento de Matem\'aticas,
Universidad Carlos III de Madrid, c/ Universidad 30, 28911 Legan\'es, Spain.}

\email[Fidalgo]{ufidalgo@@math.uc3m.es}

\author[L\'{o}pez]{G. L\'opez Lagomasino}
\address[L\'{o}pez]{Departamento de Matem\'aticas,
Universidad Carlos III de Madrid, c/ Universidad 30, 28911 Legan\'es, Spain.}

\email[L\'{o}pez]{lago@@math.uc3m.es}

\begin{document}

\maketitle

\begin{abstract}
K. Mahler introduced the concept of perfect systems in the general theory he developed for the simultaneous Hermite-Pad\'{e} approximation of analytic functions.  We prove that Nikishin systems are perfect providing, by far, the largest class of systems of functions for which this important property holds. As consequences, in the context of Nikishin systems, we obtain: an extension of Markov's theorem to simultaneous Hermite-Pad\'{e} approximation, a general result on the convergence of simultaneous quadrature rules of Gauss-Jacobi type, the logarithmic asymptotics of general sequences of multiple orthogonal polynomials, and an extension of the Denisov-Rakhmanov theorem for the ratio asymptotics of mixed type multiple orthogonal polynomials.
\end{abstract}

\vspace{1cm}


{\it A.M.S. Subject Classification.} Primary: 30E10, 42C05;
Secondary: 41A20.

\section{Introduction} \label{sec:HP}

\subsection{Some historical remarks.}
In 1873, Charles Hermite publishes in \cite{Her} his proof of the transcendence of $e$ making use of simultaneous rational approximation of systems of exponentials. That paper marked the beginning of the modern analytic theory of numbers.

The formal theory of simultaneous rational approximation for general systems of analytic functions was initiated by K. Mahler in lectures delivered at the University of Groningen in 1934-35. These lectures were published years later in \cite{Mah}. Important contributions in this respect are also due to his students J. Coates and H. Jager, see \cite{Coa} and \cite{Jag}. K. Mahler's approach to the simultaneous approximation of finite systems of analytic functions may be reformulated in the following terms.

Let ${\bf f} = (f_0,\ldots,f_m)$ be a  family of analytic functions in some domain $D$ of the extended complex plane containing $\infty$. Fix a non-zero multi-index ${\bf n} = (n_0,\ldots,n_m) \in {\mathbb{Z}}_+^{m+1}, |{\bf n}| = n_0+\ldots,n_m$. There exist polynomials $a_{{\bf n},0}, \ldots, a_{{\bf n},m}$, not all identically equal to zero, such that
\begin{itemize}
\item[i)] $\deg a_{{\bf n},j} \leq n_j -1, j=0,\ldots,m$ $(\deg a_{{\bf n},j} \leq -1$ means that $a_{{\bf n},j} \equiv 0$),
\item[ii)] $\sum_{j=0}^m a_{{\bf n},j}(z) f_j(z) - d_{\bf n}(z) = {\mathcal{O}}(1/z^{|{\bf n}|}), z \to \infty,$
\end{itemize}
for some polynomial $d_{\bf n}$. Analogously, there exists $Q_{\bf n}$, not identically equal to zero, such that
\begin{itemize}
\item[i)] $\deg Q_{{\bf n}} \leq |{\bf n}|,$
\item[ii)] $Q_{\bf n}(z)f_j(z) - P_{{\bf n},j}(z)= {\mathcal{O}}(1/z^{ n_j+1}), z \to \infty, j=0,\ldots,m,$
\end{itemize}
for some polynomials $P_{{\bf n},j}, j=0,\ldots,m.$

Initially, the polynomials $a_{{\bf n},0},\ldots,a_{{\bf n},m}$ were called latin and $Q_{\bf n}$ german polynomials, due to the letters employed in denoting them (see the papers of Mahler, Coates and Jager cited above). The polynomials $d_{\bf n}$ and $P_{{\bf n},j},j=0,\ldots,m,$ are uniquely determined from ii) once their partners are found. Later, the two constructions have been called type I and type II polynomials (approximants) of the system $(f_0,\ldots,f_m)$. Algebraically, they are closely related. This is clearly exposed in $\cite{Coa}, \cite{Jag},$ and $\cite{Mah}$. When $m=0$ both definitions coincide with that of the well-known Pad\'{e} approximation in its linear presentation.

Apart from Hermite's result, type I, type II, and a combination of the two (called mixed type), have been employed in the proof of the irrationality of other numbers. For example, in \cite{Beu} F. Beukers shows that Apery's  proof (see \cite{Ape}) of the irrationality of $\zeta(3)$  can be placed in the context of mixed type Hermite-Pad\'{e} approximation. See \cite{Ass} for a brief introduction and survey on the subject. More recently, mixed type approximation has appeared in random matrix and non-intersecting brownian motions theories (see, for example, \cite{BleKui} and \cite{DK}).

In applications in the areas of number theory, convergence of simultaneous rational approximation, and asymptotic properties of type I and type II polynomials, a central question is that of uniqueness, to within a constant factor.

\begin{defi}\label{def1}
A multi-index ${\bf n}$ is said to be {\bf normal} for the system ${\bf f}$ for type I approximation (respectively, for type II,) if $\deg a_{{\bf n},j} = n_j -1, j=0,\ldots,m$ (respectively, $\deg Q_{\bf n} = |{\bf n}|$). A system of functions ${\bf f}$ is said to be {\bf perfect} if all multi-indices are normal.
\end{defi}

It is easy to see that under normality $(a_{{\bf n},0},\ldots,a_{{\bf n},m})$ and $Q_{\bf n}$ are uniquely determined to within a constant factor.

Considering the construction at the origin (instead of $z=\infty$ which we chose for convenience), the system of exponentials considered by Hermite, $(e^{w_0z},\ldots,e^{w_mz}), w_i \neq w_j, i \neq j, i,j = 0,\ldots,m, $ is known to be perfect for type I and type II. A second example of a perfect system for both types is that given by the binomial functions $(1-z)^{w_0},\ldots,(1-z)^{w_m}, w_i -w_j \not\in {\mathbb{Z}}$. All multi-indices $\bf n$ such that $n_0 \geq \cdots \geq n_m$ are known to be type I and type II normal for $(\log^{m}(1-x),\ldots,\log(1-x),1)$. When normality occurs for multi-indices with decreasing components the system is said to be {\bf weakly perfect}. Except for systems formed by Cauchy transforms of measures, basically these are the only examples known of perfect or weakly perfect systems.

\subsection{Markov systems and orthogonality.} Let $s$ be a finite Borel measure with constant sign whose compact support consists of infinitely many points and is contained in the real line.  In the sequel, we only consider such measures. By $\Delta$ we denote the smallest interval which contains the support, $\supp{s},$ of $s$. We denote this class of measures by ${\mathcal{M}}(\Delta)$. Let
\[ \widehat{s}(z) = \int \frac{ds(x)}{z-x}
\]
denote the Cauchy transform of $s$. Obviously, $\widehat{s} \in {\mathcal{H}}(\overline{\mathbb{C}}\setminus \Delta);$  that is, it is analytic in $\overline{\mathbb{C}}\setminus \Delta.$ If we apply the construction above to the system formed by $\widehat{s}\,\, (m=0)$, it is easy to verify that $Q_{\bf n}$ turns out to be orthogonal to all polynomials of degree less than $n \in {\mathbb{Z}}_+$. Consequently, $\deg Q_{\bf n} = n,$ all its zeros are simple and lie in the open convex hull $\mbox{Co}(\supp s)$ of $\supp s$. Therefore, such systems are perfect. Let us see two other examples much more illustrating.

\subsubsection{Angelesco systems.} In \cite{Ang}, A. Angelesco considered the following systems of functions. Let $\Delta_j, j=0,\ldots,m,$ be pairwise disjoint bounded intervals contained in the real line and $s_j, j=0,\ldots,m,$ a system of measures such that $\mbox{Co}(\supp s_j) = \Delta_j.$

Fix ${\bf n} \in {\mathbb{Z}}_+^{m+1}$ and consider the type II approximant of the so called Angelesco system of functions $(\widehat{s}_0,\ldots,\widehat{s}_m)$ relative to $\bf n$. It turns out that
\[ \int x^{\nu} Q_{\bf n}(x) ds_j(x) =0, \quad \nu = 0,\ldots,n_j -1,\quad j=0,\ldots,m.
\]
Therefore, $Q_{\bf n}$ has $n_j$ simple zeros in the interior (with respect to the euclidean topology of ${\mathbb{R}}$) of $\Delta_j$. In consequence, since the intervals $\Delta_j$ are pairwise disjoint, $\deg Q_{\bf n} = |{\bf n}|$ and Angelesco systems are type II perfect. Type I perfectness for Angelesco systems has not been studied.

Unfortunately, Angelesco's paper received little attention and such systems  reappear many years later in \cite{Nik1} where E.M. Nikishin  deduces some of their formal properties.

Though type II normality for Angelesco systems is so easy to deduce, the multiple orthogonal polynomials and the rational approximations associated with them do not have good asymptotic behavior. In \cite{GR} and \cite{Apt1}, their logarithmic and strong asymptotic formulas, respectively, are given. In this respect, a different  system of Markov functions turns out to be much more interesting and foundational from the geometric and analytic points of view.

\subsubsection{Nikishin systems.}  In an attempt to construct general classes of functions for which normality takes place, in \cite{Nik} E.M. Nikishin introduced the concept of MT-system. Let $\Delta_{\alpha}, \Delta_{\beta}$ be two non intersecting bounded intervals contained in the real line and $\sigma_{\alpha} \in {\mathcal{M}}(\Delta_{\alpha}), \sigma_{\beta} \in {\mathcal{M}}(\Delta_{\beta})$. With these two measures we define a third one as follows (using the differential notation)
\[ d \la \sigma_{\alpha},\sigma_{\beta} \ra (x) = \widehat{\sigma}_{\beta}(x) d\sigma_{\alpha}(x);
\]
that is, one multiplies the first measure by a weight formed by the Cauchy transform of the second measure. Certainly, this product of measures is non commutative.

Above, $\widehat{\sigma}_{\beta}$ denotes the Cauchy transform of the measure $\sigma_{\beta}$. The reader may argue, and we agree, that the appropriate notation is $\widehat{\sigma_{\beta}}$. However, throughout the paper, we will need Cauchy transforms of measures with several sub-indices and supra-indices; for example $\widehat{s_{1,j}^2}$ (and much more extended). The correct notation causes space consumption and aesthetic inconveniences. So, be cautious, $\widehat{s}_{1,j}^2$ is not  the Cauchy transform of $s$ sub-indexed with $1,j$ and then squared, but precisely the Cauchy transform of a measure denoted $s_{1,j}^2$. The good news is that powers rarely appear in the paper and they are clear from the context (for example, $(-1)^j$ or $z^2$).

\begin{defi} \label{Nikishin} Take a collection  $\Delta_j, j=0,\ldots,m,$ of intervals such that
\[ \Delta_j \cap \Delta_{j+1} = \emptyset, \qquad j=0,\ldots,m-1.
\]
Let $(\sigma_0,\ldots,\sigma_m)$ be a system of measures such that $\mbox{Co}(\supp \sigma_j) = \Delta_j, \sigma_j \in {\mathcal{M}}(\Delta_j), j=0,\ldots,m.$
We say that $(s_0,\ldots,s_m) = {\mathcal{N}}(\sigma_0,\ldots,\sigma_m)$, where
\[ s_0 = \sigma_0, \quad s_1 = \la \sigma_0,\sigma_1 \ra, \ldots \quad , s_m = \la \sigma_0,\la \sigma_1,\ldots,\sigma_m \ra \ra
\]
is the Nikishin system of measures generated by $(\sigma_0,\ldots,\sigma_m)$.
\end{defi}

Fix ${\bf n} \in {\mathbb{Z}}_+^{m+1}$ and consider the type II approximant of the Nikishin system of functions $(\widehat{s}_0,\ldots,\widehat{s}_m)$ relative to $\bf n$. It is easy to prove that
\[ \int x^{\nu} Q_{\bf n}(x) ds_j(x) =0, \quad \nu = 0,\ldots,n_j -1,\quad j=0,\ldots,m.
\]
All the measures $s_j$ have the same support; therefore, it is not immediate to conclude that $\deg Q_{\bf n} = |{\bf n}|$. Nevertheless, if we denote
\[s_{j,k} = \la \sigma_j,\sigma_{j+1},\ldots,\sigma_k \ra,\qquad j < k, \quad s_{j,j} = \la \sigma_j \ra = \sigma_j,
\]
the previous orthogonality relations may be rewritten as follows
\begin{equation}
\label{eq:a}
\int (p_0 (x) + \sum_{k=1}^m p_j(x)\widehat{s}_{1,k}(x))Q_{\bf n}(x) d \sigma_0(x) =0,
\end{equation}
where $p_0,\ldots,p_m$ are arbitrary polynomials such that $\deg p_k \leq n_k -1, k=0,\ldots,m.$

\begin{defi} \label{def:AT} A system of real continuous functions $u_0,\ldots,u_m$ defined on an interval $\Delta$ is called an AT-system on $\Delta$ for the multi-index ${\bf n} \in {\mathbb{Z}}_+^{m+1}$ if for any choice of real polynomials (that is, with real coefficients) $p_0,\ldots,p_m, \deg p_k \leq n_k -1,$ the function
\[ \sum_{k=0}^m p_k(x) u_k(x)
\]
has at most $|{\bf n}| -1$ zeros on $\Delta$. If this is true for all ${\bf n} \in {\mathbb{Z}}_+^{m+1}$ we have an AT system on $\Delta$.
\end{defi}

In other words, $u_0,\ldots,u_m$ forms an AT-system for ${\bf n}$ on $\Delta$ when the system of functions
\[ (u_0,\ldots, x^{n_0 -1}u_0,u_1,\ldots,x^{n_m-1}u_m)
\]
is a Tchebyshev system on $\Delta$ of order $|{\bf n}| -1$. From the properties of Tchebyshev systems (see \cite[Theorem 1.1]{KN}), it follows that given
$x_1,\ldots,x_N, N < |{\bf n}|,$ points in the interior of $\Delta$ one can find polynomials $h_0,\ldots,h_m,$ conveniently, with $\deg h_k \leq n_k -1,$ such that $\sum_{k=0}^m h_k(x) u_k(x)$
changes  sign at $x_1,\ldots,x_N,$ and has no other points where it changes sign on $\Delta.$

In \cite{Nik}, Nikishin stated without proof that the system of functions $(1,\widehat{s}_{1,1},\ldots,\widehat{s}_{1,m})$ forms an AT-system for all multi-indices ${\bf n}$ such that $n_0 \geq \cdots \geq n_m$ (he proved it when additionally $n_0 -n_m \leq 1$). Due to (\ref{eq:a}) this implies that Nikishin systems are type II weakly perfect.

The proof of Nikishin's assertion is a consequence of \cite[Theorem 4.1]{DrSt0}. Actually, Driver and Stahl proved normality for the wider class of multi-indices
\[ {\mathbb{Z}}_+^{m+1}(\circledast) = \{{\bf n} \in {\mathbb{Z}}_+^{m+1}: 0 \leq j < k \leq m \Rightarrow n_k \leq n_j +1\}
\]
(see also \cite{GRS}). In \cite[Theorem 4.2]{DrSt0},  for the same class of multi-indices, the authors also proved type I weak perfectness for Nikishin systems. In the same paper (see remark on page 171), it is shown that when $m=1$ Nikishin systems are type II perfect. Improvements are also contained in \cite[Theorem 1]{BBFG} (see also \cite[Theorem 2]{LIF}) where normality is proved for all multi-indices in
\[ {\mathbb{Z}}_+^{m+1}(*) = \{{\bf n} \in {\mathbb{Z}}_+^{m+1}: \not\exists\,\, 0 \leq i < j < k \leq m\,\, \mbox{such that}\,\, n_i <  n_j < n_k \}
\]
and \cite[Theorem 1]{LF1} containing the proof that for $m=2$ Nikishin systems are type II perfect.

Ever since the appearance of \cite{Nik}, a subject of major interest for those involved in simultaneous approximation was to determine whether or not Nikishin systems are perfect. The main result of this paper gives a positive answer to this question. Moreover, we will prove perfectness for mixed type Nikishin systems, containing type I and type II as particular cases. The proof is based on the reduction of the problem to the case of multi-indices with decreasing components (that is, to weak perfectness). In the sequel,
\[  {\mathbb{Z}}_+^{m+1}(\bullet) = \{{\bf n} \in {\mathbb{Z}}_+^{m+1}:  n_0 \geq \cdots \geq n_m\}.
\]
Notice that
\[{\mathbb{Z}}_+^{m+1}(\bullet) \subset {\mathbb{Z}}_+^{m+1}(\circledast) \subset {\mathbb{Z}}_+^{m+1}(*) \subset {\mathbb{Z}}_+^{m+1}.
\]
When $m=0$ these sets are equal. For $m=1$ the last two coincide. If $m\geq 2$ they are all distinct.

The proof of the main result relies on interesting reduction formulas concerning products and ratios of Cauchy transforms of measures. We will see numerous consequences of the perfectness of Nikishin systems in: convergence of simultaneous Pad\'e approximation, convergence of simultaneous quadrature rules, and asymptotic properties of multiple orthogonal polynomials.

\subsubsection{Mixed type Nikishin systems.} In \cite{Sor0}, Sorokin introduced the following construction. Let
\[ {\mathbb{F}} = (f_{j.k}),
\]
be an $(m_2+1)\times (m_1+1)$ dimensional matrix  of analytic functions in some domain $D$ of the extended complex plane containing $\infty$. Fix a multi-index ${\bf n} = ({\bf n}_1;{\bf n}_2) \in {\mathbb{Z}}_+^{m_1 +1}\times {\mathbb{Z}}_+^{m_2 +1}$, such that $|{\bf n}_1| = |{\bf n}_2| +1$. We denote ${\bf n}_i =(n_{i,0},\ldots,n_{i,m_i}), i=1,2.$ There exists a  vector polynomial ${\mathbb{A}}_{\bf n}=(a_{{\bf n},0}, \ldots, a_{{\bf n},m_1})$,  such that
\begin{itemize}
\item[a)] ${\mathbb{A}}_{\bf n} \not \equiv {\bf 0}, \deg a_{{\bf n},k} \leq n_{1,k} -1, k=0,\ldots,m_1, $
\item[b)] $({\mathbb{F}}{\mathbb{A}}_{\bf n}^t - {\mathbb{D}}_{\bf n}^t)(z) = ({\mathcal{O}}(1/z^{n_{2,0} +1}),\ldots,{\mathcal{O}}(1/z^{n_{2,m_2} +1} )^t =: {\mathcal{O}}(1/z^{{\bf n}_2 +1}), z \to \infty.$
\end{itemize}
for some $m_2 +1$ dimensional vector polynomial ${\mathbb{D}}_{\bf n}$ (the super-index $t$ means taking transpose and {\bf 0} denotes the zero vector). Finding ${\mathbb{A}}_{\bf n}$ reduces to solving a linear homogeneous system of $|{\bf n}_2|$ equations determined by the conditions b) on $|{\bf n}_1|$ unknowns (the total number of coefficients of the polynomials $a_{{\bf n},k}, k=0,\ldots,m_1$). Since $|{\bf n}_2| +1 = |{\bf n}_1|$ a non trivial solution exists.

\begin{defi}
\label{defiSor}
 A non zero vector ${\mathbb{A}}_{\bf n}$ satisfying a)-b) is called mixed type vector polynomial relative to ${\mathbb{F}}$ and ${\bf n} \in {\mathbb{Z}}_+^{m_1 +1}\times {\mathbb{Z}}_+^{m_2 +1}, |{\bf n}_1| = |{\bf n}_2| +1$. If $\deg a_{{\bf n},k} = n_{1,k} -1, k=0,\ldots,m_1$, the multi-index ${\bf n}$ is called mixed type normal. ${\mathbb{F}}$ is mixed type perfect when all multi-indices in ${\mathbb{Z}}_+^{m_1 +1}\times {\mathbb{Z}}_+^{m_2 +1}$ such that $|{\bf n}_1| = |{\bf n}_2| +1$ are normal.
\end{defi}

This construction has as particular cases type I $(m_2 = 0)$ and type II $(m_1 =0)$ polynomials.
\vspace{0,2cm}

Let $S^{1} = (s^1_{0,0},\ldots,s^1_{0,m_1}) = {\mathcal{N}}(\sigma_0^1,\ldots,\sigma_{m_1}^1), S^{2} = (s^2_{0,0},\ldots,s^2_{0,m_2})= {\mathcal{N}}(\sigma_0^2,\ldots,\sigma_{m_2}^2), \sigma_0^1 = \sigma_0^2,$ be two given Nikishin systems generated by $m_1+1$ and $m_2+1$ measures, respectively. We underline the fact that both Nikishin systems stem from the same basis measure $\sigma_0^1 = \sigma_0^2$, but there is no other restriction on them.
Let us  introduce the row vectors
\[ {\mathbb{U}} =  (1,\widehat{s}^{2}_{1,1},\ldots,\widehat{s}^{2}_{1,m_2}),
\qquad {\mathbb{V }} =
(1,\widehat{s}^{1}_{1,1},\ldots,\widehat{s}^{1}_{1,m_1})
\]
and the $(m_2+1) \times (m_1 +1)$ dimensional  matrix function
\[ {\mathbb{W}} = {\mathbb{U}}^t {\mathbb{V}}.
\]
Define the matrix Markov type function
\[ \widehat{\mathbb{S}}(z) = \int \frac{{\mathbb{W}}(x)d\sigma_0^2 (x)}{z-x}
\]
understanding that integration is carried out entry by entry on
the matrix $\mathbb{W}$. We say that $\widehat{\mathbb{S}}$ is a mixed type Nikishin system of functions.

In the rest of the paper, we will study mixed type Nikishin systems and their mixed type polynomials. Occasionally, we reduce the study to type II $(m_1 =0)$. In such cases, for simplicity, we reduce the notation. Namely, ${\bf n} = (n_0,\ldots,n_m)$, the vector function will be ${\bf f} = (\widehat{s}_0, \ldots, \widehat{s}_m)$, where $m = m_2$, and $(s_0,\ldots,s_m) = {\mathcal{N}}(\sigma_0,\ldots,\sigma_m)$. The mixed type polynomials ${\mathbb{A}}_{\bf n}$ will then be denoted $Q_{\bf n}$.

\subsection{Statement of the main results.} Mixed type Nikishin systems and their associated mixed type polynomials satisfy many interesting properties. Let us begin with

\begin{teo}
\label{teo:1}
Let $(s_{1,1},\ldots,s_{1,m}) = {\mathcal{N}}(\sigma_{1},\ldots,\sigma_m)$ be given. Then, the system  $(1,\widehat{s}_{1,1},\ldots,\widehat{s}_{1,m})$ forms an AT-system on any interval $\Delta$ disjoint from $\Delta_1 = \mbox{\rm Co}(\supp \sigma_1)$. Moreover, for each ${\bf n} \in {\mathbb{Z}}_+^{m+1},$ and arbitrary polynomials with real coefficients $p_k, \deg p_k \leq n_k-1, k=0,\ldots,m,$ the linear form $p_0 + \sum_{k=1}^m p_k \widehat{s}_{1,k}, $ has at most $|{\bf n}|-1$ zeros in $\overline{\mathbb{C}} \setminus \Delta_1$.
\end{teo}

From here, we can prove

\begin{teo}\label{teo:2} The matrix $\widehat{\mathbb{S}}$ is mixed type perfect. For each ${\bf n} \in {\mathbb{Z}}_+^{m_1 +1}\times {\mathbb{Z}}_+^{m_2 +1}, |{\bf n}_1| = |{\bf n}_2| +1,$ the vector polynomial ${\mathbb{A}}_{\bf n}$ is uniquely determined up to a constant factor.
\end{teo}

In particular, this means that Nikishin systems are type I and type II perfect.

An easy consequence of Theorem \ref{teo:2} (more precisely, of Lemma \ref{lem:4} below) is that the linear form
\[ {\mathcal{A}}_{{\bf n}}:= a_{{\bf n},0} + \sum_{k=1}^{m_1} a_{{\bf n},k} \widehat{s}_{1,k}^1,
\]
has at least $|{\bf n}_2|$ sign changes in the interior (with respect to the euclidean topology of ${\mathbb{R}}$) of the interval $\Delta_0$.  In particular, for any ${\bf n} \in {\mathbb{Z}}_+^{m+1}$, the type II polynomial $Q_{\bf n}(= a_{{\bf n},0} )$ has all its zeros located inside  $\mbox{Co}(\supp \sigma_0)$. This has a striking consequence in terms of the convergence of type II rational approximation of Nikishin systems.

In \cite{Nik}, it was proved that for $(s_0,s_1)= {\mathcal{N}}(\sigma_0,\sigma_1)$
\[ \lim_{n\to \infty} \frac{P_{{\bf n},k}}{Q_{\bf n}} = \widehat{s}_k, \qquad k=0,1,
\]
uniformly on compact subsets of $\overline{\mathbb{C}} \setminus \mbox{\rm Co}(\supp \sigma_0),$ where the limit is taken along the sequence ${\bf n} = (n,n), n \in {\mathbb{Z}}_+$. When $m_2=0$ the corresponding result is Markov's classical theorem on the convergence of Pad\'{e} approximants, see \cite{Mar}. The extension to diagonal sequences for arbitrary $m$  and ${\bf n} \in \Lambda \subset {\mathbb{Z}}_{+}^{m+1}(\bullet)$ such that $|{\bf n}|\to \infty$ and $\max\{n_{0}-n_{m}: {\bf n} \in \Lambda\} < \infty$ is contained in \cite[Corollary 1]{Bus} (which also includes the case when the measures have unbounded support and a Carleman type condition is satisfied). Combining Theorem \ref{teo:2} and \cite[Theorem 1]{LF3} we obtain

\begin{cor}
\label{cor:1} Let $(s_{0,0},\ldots,s_{0,m}) = {\mathcal{N}}(\sigma_0,\ldots,\sigma_m)$ and $\Lambda \subset {\mathbb{Z}}_+^{m+1}$ be given. Assume that there exist constants $c > 0, \kappa <1,$ such that
\[ n_j \geq \frac{|{\bf n}|}{m+1}  - c |{\bf n}|^{\kappa},  \qquad j=0,\ldots,m.
\]
Then,
\[ \lim_{{\bf n} \in \Lambda} \frac{P_{{\bf n},k}}{Q_{\bf n}} = \widehat{s}_{0,k} \qquad {\mathcal{K}} \subset \overline{\mathbb{C}} \setminus \mbox{\rm Co}(\supp \sigma_0), \qquad k=0,\ldots,m.
\]
Moreover,
\[ \limsup_{{\bf n}\in \Lambda} \|\widehat{s}_{0,k} - \frac{P_{{\bf n},k}}{Q_{\bf n}}\|_{\mathcal{K}}^{1/2|{\bf n}|} \leq
\delta_{\mathcal{K}} <1, \qquad k=0,\ldots,m,
\]
where $\|\cdot\|_{\mathcal{K}}$ denotes the uniform norm on $\mathcal{K}$,
\[ \delta_{\mathcal{K}} = \max \{|\varphi_t(z)|: z \in {\mathcal{K}}, t \in \mbox{\rm Co}(\supp \sigma_1) \cup\{\infty\}\},
\]
and $\varphi_t$ represents conformally $\overline{\mathbb{C}} \setminus \mbox{\rm Co}(\supp \sigma_0)$
onto the unit circle with $\varphi_t(t) = 0, \varphi_t'(t) > 0$.
\end{cor}

Throughout the paper, the notation
\[ \lim_{n \in \Lambda} g_n(z) = g(z), \qquad
{\mathcal{K}} \subset \Omega, \] stands for uniform convergence of
the sequence of functions $\{g_n\}, n\in \Lambda,$ to the function $g$ on each compact
subset $\mathcal{K}$ contained in the indicated region (in this
case $\Omega$).

This result gives a very general extension of Markov's theorem. Notice that the sequences are not required to be close to diagonal (equal components in the multi-indices).

Since the zeros $x_{{\bf n},j}, j=1,\ldots,|{\bf n}|,$ of $Q_{\bf n}$ are simple, we can decompose $P_{{\bf n},k}/Q_{\bf n}$ as follows
\[\frac{P_{{\bf n},k}(z)}{Q_{\bf n}(z)} = \sum_{i=1}^{|{\bf n}|} \frac{\lambda_{{\bf n},k,i}}{z - x_{{\bf n},i}}, \qquad \lambda_{{\bf n},k,i} = \lim_{z \to x_{{\bf n},i}} (z - x_{{\bf n},i})\frac{P_{{\bf n},k}(z)}{Q_{\bf n}(z)} = \frac{P_{{\bf n},k}(x_{{\bf n},i})}{Q_{\bf n}'(x_{{\bf n},i})}.
\]

The following analogue of the Gauss-Jacobi quadrature formula takes place.

\begin{cor} \label{cor:2} Let $(s_{0,0},\ldots,s_{0,m}) = {\mathcal{N}}(\sigma_0,\ldots,\sigma_m)$ and ${\bf n} = {\mathbb{Z}}_+^{m+1}$ be given. Then, for each $k=0,\ldots,m$ and every polynomial $p, \deg p \leq |{\bf n}| + n_k -1$
\[ \int p(x) ds_{0,k}(x) = \sum_{i=1}^{|{\bf n}|} \lambda_{{\bf n},k,i} p(x_{{\bf n},i}).
\]
If ${\bf n} = (n,n+1,\ldots,n+1)$, then
\[  \mbox{\rm sign} (\lambda_{{\bf n},k,i}) =  \mbox{\rm sign}(s_{0,k}), \qquad i= 1,\ldots,|{\bf n}|.
\]
Consequently, for the sequence of multi-indices $\{(n,n+1,\ldots,n+1)\}_{n \in {\mathbb{Z}}_+} \subset {\mathbb{Z}}_+^{m+1}$, for any bounded Riemann-Stieltjes integrable function $f$ on $\mbox{\rm Co}(\supp \sigma_0)$ and each $k=0,\ldots,m$
\[ \int f(x) ds_{0,k}(x) = \lim_{n \to \infty}\sum_{i=1}^{|{\bf n}|} \lambda_{{\bf n},k,i} f(x_{{\bf n},i}).
\]
\end{cor}

This result provides convergence of the quadrature formulas simultaneously for all the measures in the Nikishin system taking the same nodes in all the quadrature formulas. Simultaneous quadrature formulas were studied in \cite{Bor} in connection with certain application to computer graphics illuminating bodies. Whenever feasible, simultaneous quadrature formulas are more efficient, from the computational point of view, compared with the use of Gauss-Jacobi quadrature independently on each measure. In \cite{LIF}, a more detailed study of simultaneous quadrature formulas for Nikishin systems of measures may be found. We wish to point out that for the class of multi-indices considered in Corollary \ref{cor:2} all the statements in \cite[Corollary 2]{LIF} hold true for all $k=0,\ldots,m.$

Theorem \ref{teo:1} follows easily from

\begin{teo}
\label{teo:3}
Let $ (s_{1,1},\ldots,s_{1,m}) = {\mathcal{N}}(\sigma_{1},\ldots,\sigma_m)$ and ${\bf n} = (n_0,\ldots,n_m) \in {\mathbb{Z}}_+^{m+1}$ be given. Then, there exists a permutation $\lambda$ of $(0,\ldots,m)$ which reorders the components of ${\bf n}$ decreasingly, $n_{\lambda(0)} \geq \cdots \geq n_{\lambda(m)},$ and an associated Nikishin system $S(\lambda) = (r_{1,1},\ldots,r_{1,m}) = {\mathcal{N}}(\rho_{1},\ldots,\rho_m)$ such that for any real polynomials
$p_k, \deg p_k \leq n_k -1$, there exist real polynomials $q_k$ such that
\[ p_0 + \sum_{k=1}^m p_k \widehat{s}_{1,k} = (q_0 + \sum_{k=1}^m q_k \widehat{r}_{1,k})\widehat{s}_{1,\lambda(0)}, \qquad \deg q_k \leq n_{\lambda(k)} -1, \qquad k=0,\ldots,m.
\]
\end{teo}

We wish to point out that here $\widehat{s}_{1,0}$ denotes the function identically equal to $1$; this is relevant when $\lambda(0) =0$.  There may be several permutations $\lambda$ for which the statement holds, each one with an associated $S(\lambda)$. We do not know if there is an $S(\lambda)$ for each $\lambda$ which reorders the components of ${\bf n}$ decreasingly. As reference, we can say that there exists $S(\lambda)$ (but not exclusively) for that $\lambda$ which  additionally satisfies that for all $0 \leq j <k \leq n$ with $n_j = n_k$ then also $\lambda(j) < \lambda(k)$.

\begin{teo}\label{teo:4} Suppose that  $S^1 = (s^1_{0,0},\ldots,s^1_{0,m_1}) = {\mathcal{N}}(\sigma_0^1,\ldots,\sigma_{m_1}^1),  S^2=(s^2_{0,0},\ldots,s^2_{0,m_2})= {\mathcal{N}}(\sigma_0^2,\ldots,\sigma_{m_2}^2), \sigma_0^1 = \sigma_0^2,$ and ${\bf n} = ({\bf n}_1;{\bf n}_2)\in {\mathbb{Z}}_+^{m_1 +1}\times {\mathbb{Z}}_+^{m_2 +1}, |{\bf n}_1| = |{\bf n}_2| +1,$ be given.  Let $\lambda_2$ and $S(\lambda_2)= {\mathcal{N}}(\rho_1^2,\ldots,\rho_{m_2}^2)$ be a permutation and a Nikishin system associated with ${\mathcal{N}}(\sigma_1^2,\ldots,\sigma_{m_2}^2)$ and ${\bf n}_2$ by Theorem $\ref{teo:3}$.  Construct $(r^2_{0,0},\ldots,r^2_{0,m_2})= {\mathcal{N}}(\rho_0^2,\ldots,\rho_{m_2}^2)$, where  $\rho^2_0 = \widehat{s}^2_{1,\lambda_2 (0)}\sigma_0^2$.  Then
\[ \int x^{\nu} {\mathcal{A}}_{{\bf n}}(x) d r_{0,k}^2(x) =0,\qquad \nu =0,\ldots,n_{2,\lambda_2 (k)} -1, \qquad k=0,\ldots,m_2.
\]
\end{teo}

Using Theorem \ref{teo:4} and results from \cite{FLLS}, we can obtain logarithmic and ratio asymptotics for sequences $\{{\mathcal{A}}_{{\bf n}}\}_{{\bf n} \in \Lambda}, \Lambda \subset {\mathbb{Z}}_+^{m_1+1} \times {\mathbb{Z}}_+^{m_2+1}, |{\bf n}_1| = |{\bf n}_2| +1,$ under appropriate assumptions on the measures generating $S^1,S^2,$ and $\Lambda$.

A positive measure $\sigma$ is
said to be regular if
\[ \lim_{n \to \infty} \kappa_n^{1/n} = 1/\mbox{cap}(\supp \sigma),
\]
where $\mbox{cap}(\cdot)$ denotes the logarithmic capacity of the
Borel set $(\cdot)$ and $\kappa_n$ denotes the leading coefficient
of the $n$-th orthonormal polynomial with respect to $\sigma$. A negative measure $\sigma$ is regular if $-\sigma$ is regular. In either cases we write $\sigma \in \mbox{\bf Reg}$. For
equivalent forms of defining regular measures, see
sections 3.1 to 3.3 in \cite{stto} (in particular Theorem 3.1.1). For short, we write $ (S^1,S^2) \in \mbox{\bf Reg}$ to mean that
all the measures which generate both Nikishin systems $(S^1,S^2)$ are regular.

A region $\Omega$ of the extended complex plane which has a compact complement $E$ is said to be regular if the Dirichlet problem has a solution on $\Omega$ for any continuous function defined on $\partial U = \partial E$. This is equivalent to proving that the Green's function on $\Omega$ with singularity at $\infty$ can be extended continuously to all $\mathbb{C}$ (for details on Green's function and regular domains see \cite[Theorem 10.12]{He}). In this case it is usual to say also that $E$ is a regular compact set.

\begin{defi}
\label{quasi} We say that a compact set $E$ is quasi-regular when $E
= \widetilde{E} \cup e$, where $\widetilde{E}$ is a regular compact set
and $e$ is at most a denumerable set whose accumulation points lie
in $\widetilde{E}.$
\end{defi}

Let $\Lambda =
\Lambda(p_{1,0},\ldots,p_{1,m_1};p_{2,0},\ldots,p_{2,m_2}) \subset
{\mathbb{Z}}_+^{m_1+1} \times {\mathbb{Z}}_+^{m_2
+1} $ be an infinite sequence of distinct multi-indices
such that for each ${\bf n} = ({\bf n}_1,{\bf n}_2) \in \Lambda, |{\bf n}_1| = |{\bf n}_2| +1$, and
\[
\lim_{ {\bf n} \in \Lambda} \frac{n_{1,j}}{|{\bf n}_1|} = p_{1,j}
\in (0,1), \quad j=0,\ldots,m_1, \quad \lim_{{\bf n} \in \Lambda}
\frac{n_{2,j}}{|{\bf n}_2|} = p_{2,j} \in (0,1), \quad
j=0,\ldots,m_2.
\]

The following two results require some normalization on the sequence of linear forms under consideration. By Theorem \ref{teo:2}, for each ${\bf n} \in {\mathbb{Z}}_+^{m_1+1} \times
{\mathbb{Z}}_+^{m_2+1},$ ${\mathcal{A}}_{\bf n}$ is uniquely determined except for a constant factor.

\begin{defi} \label{monico} Let $\min\{n_{1,0},\ldots,n_{1,m_1}\} \geq 1$. We say that ${\mathcal{A}}_{\bf n}$ is monic if the leading coefficient of $a_{{\bf n},j}$ is $1$, where $j$ is either $m_1$ when $n_{1,m_1} = \min\{n_{1,0},\ldots,n_{1,m_1}\}$, or if $n_{1,m_1} > \min\{n_{1,0},\ldots,n_{1,m_1}\}$ it is such that $n_{1,j} = \min\{n_{1,0},\ldots,n_{1,m_1}\}$ and $n_{1,j} < n_{1,k}, j < k \leq m_1$.
\end{defi}

We do not need to normalize ${\mathcal{A}}_{\bf n}$ when ${\bf n}_1$ has components equal to zero. We have

\begin{teo}\label{teo:5}
Let $\Lambda =
\Lambda(p_{1,0},\ldots,p_{1,m_1};p_{2,0},\ldots,p_{2,m_2}) \subset
{\mathbb{Z}}_+^{m_1+1} \times
{\mathbb{Z}}_+^{m_2+1}, (S^1,S^2)\in\mbox{\bf {Reg}},$ $
S^1= {\mathcal{N}}(\sigma^1_0,\ldots,\sigma^1_{m_1}),$ and $S^2=
{\mathcal{N}}(\sigma^2_0,\ldots,\sigma^2_{m_2})$ be given. Assume that the supports of the measures which generate $S^1, S^2$ are quasi-regular. Then, the associated sequence of monic mixed type multiple orthogonal linear forms $\{{\mathcal{A}}_{\bf n}\}, {\bf n} \in \Lambda$, satisfies
\[
\lim_{{\bf n} \in \Lambda}\left|{\mathcal{A}}_{{\bf n}}(z)\right|^{1/|{\bf n}_1|}=G(z), \qquad {\mathcal{K}}
\subset {{\mathbb{C}}}\setminus (\Delta^1_0\cup\Delta^1_1),
\]
where $\Delta^1_i = \mbox{\rm Co}(\supp \sigma_i^1), i=0,1.$
\end{teo}

A formula for $G$ is given in (\ref{eq:lim5}) during the proof of Theorem \ref{teo:5}. It is expressed in terms of the solution of a vector equilibrium problem for the logarithmic potential. The  matrix governing the interaction between the different potentials in the system depends on $(p_{1,0},\ldots,p_{1,m_1};p_{2,0},\ldots,p_{2,m_2})$.

By allowing quasi-regularity of the supports, this theorem is already novel for standard orthogonal polynomials $(m_1=m_2=0)$ (see Lemma \ref{lem:asint} below). Theorem $\ref{teo:5}$ unifies the study of the logarithmic asymptotics of type I and type II multiple orthogonal polynomials under the most general conditions on the measures, their supports, and the behavior of the sequence of multi-indices on which the limit is taken. It is motivated by the results of \cite{busta}, \cite{GRS} and \cite{nik2} where type I and type II were considered separately, and the generating measures are  supported on intervals on which their Radon Nikodym derivative is positive almost everywhere. In \cite[Theorem 1.3]{FLLS} an analogue was obtained assuming that the components of ${\bf n}_1, {\bf n}_2$ are decreasing and the supports of the generating measures are regular sets. The first study of the logarithmic asymptotics of mixed type multiple orthogonal polynomials of Nikishin systems was carried out in  \cite{Sor}.

For the next result, we assume that $\supp (\sigma^{i}_j) =
\widetilde{\Delta}^{i}_j \cup e^{i}_j, j=0,\ldots,m_i, i=1,2$, where
$\widetilde{\Delta}^{i}_j$ is a bounded interval of the real line,
$|(\sigma^{i}_j)^{\prime}| > 0$ a.e. on
$\widetilde{\Delta}_j^{i}$, and $e^{i}_j$ is at most a denumerable
set without accumulation points in ${\mathbb{R}} \setminus
\widetilde{\Delta}_j^{i}$. We denote this writing
\[S^1=
{\mathcal{N}}'(\sigma^1_0,\ldots,\sigma^1_{m_1}), \qquad S^2=
{\mathcal{N}}'(\sigma^2_0,\ldots,\sigma^2_{m_2}).\]
Notice the ``prime'' on ${\mathcal{N}}$.
In the context of this paper, this condition is the analogue of the one imposed by S. A. Denisov (see \cite{Denisov}) in his extension of E. A. Rakhmanov's celebrated theorem on ratio asymptotics of orthogonal polynomials. The original proof of Rakhmanov's theorem is in \cite{kn:Rak1}-\cite{kn:Rak2}. An improved and reduced version of the proof by the author may be found in \cite{kn:Rak3}.

Fix a vector
$l:=(l_{1};l_{2})$ where $0\leq l_{1}\leq m_{1}$ and $0\leq
l_{2}\leq m_{2}$. We define the multi-index ${\bf n}^{l}:=({\bf
n}_{1}+{\bf e}^{l_{1}};{\bf n}_{2}+{\bf e}^{l_{2}})=({\bf
n}_{1}^{l_{1}};{\bf n}_{2}^{l_{2}})$, where ${\bf e}^{l_{i}}$
denotes the unit vector of length $m_{i}+1$ with all components
equal to zero except the component $(l_{i}+1)$ which equals $1$.

Fix two permutations $\lambda_1,\lambda_2,$ of $(0,\ldots,m_1)$ and $(0,\ldots,m_2)$, respectively, and a positive number $C$. By $\Lambda(\lambda_1,\lambda_2,C)$ we denote the set of all  multi-indices ${\bf n} = ({\bf n}_1;{\bf n}_2) \in {\mathbb{Z}}_+^{m_1 +1}\times {\mathbb{Z}}_+^{m_2 +1}$ such that
\begin{itemize}
\item[a)] $|{\bf n}_1| = |{\bf n}_2| +1,$
\item[b)] $\lambda_1, S(\lambda_1),$ and $\lambda_2, S(\lambda_2),$ are solutions given by Theorem \ref{teo:3} to ${\bf n}_1, {\mathcal{N}}(\sigma_1^1,\ldots,\sigma_{m_1}^1),$ and ${\bf n}_2, {\mathcal{N}}(\sigma_1^2,\ldots,\sigma_{m_2}^2),$ respectively.
\item[c)] $n_{i,\lambda_i(0)} - n_{i,\lambda_i(m_i)}  \leq C, i=1,2,$
\end{itemize}

Any sequence $\Lambda \subset {\mathbb{Z}}_+^{m_1 +1}\times {\mathbb{Z}}_+^{m_2 +1}$ of distinct multi-indices satisfying a) and
\[ \sup\{\max\{n_{i,0},\ldots,n_{i,m_i}\}-\min\{n_{i,0},\ldots,n_{i,m_i}\}: {\bf n} \in \Lambda, i=1,2\} < \infty
\]
is contained in $\cup_{\lambda_1,\lambda_2} \Lambda(\lambda_1,\lambda_2,C)$ for some sufficiently large $C$, where the union is taken over all possible pairs of permutations. Thus, any such $\Lambda$ can be partitioned in a finite number of sequences of indices satisfying a)-c) for the same pair $\lambda_1,\lambda_2$ of permutations, plus a set containing a finite number of multi-indices.

\begin{teo} \label{teo:6}
Let $S^1= {\mathcal{N}}'(\sigma^1_0,\ldots,\sigma^1_{m_1}), S^2=
{\mathcal{N}}'(\sigma^2_0,\ldots,\sigma^2_{m_2})$ and $\lambda_1, \lambda_2,$ be given. Fix  $l =(l_{1};l_{2}), 0\leq l_{1}\leq m_{1},
0 \leq l_{2}\leq m_{2}.$ Let $\Lambda \subset {\mathbb{Z}}_+^{m_1 +1}\times {\mathbb{Z}}_+^{m_2 +1}$ be an infinite sequence of distinct multi-indices such that for all ${\bf n} \in \Lambda$, ${\bf n}, {\bf n}^l \in \Lambda(\lambda_1,\lambda_2,C)$ for some sufficiently large $C$.    Then,
the associated sequence of monic mixed type multiple orthogonal linear forms $\{{\mathcal{A}}_{\bf n}\}, {\bf n} \in \Lambda$, verifies
\[
\lim_{{\bf n}\in\Lambda}\,\frac{{\mathcal{A}}_{{\bf
n}^l}(z)}{{\mathcal{A}}_{{\bf n}}(z)}={\mathcal{A}}^{(l)}(z),\qquad  {\mathcal{K}}
\subset{\mathbb{C}}\setminus(\supp{\sigma_0^1} \cup \mbox{\rm Co}(\supp{\sigma_1^1}))\,,
\]
where ${\mathcal{A}}^{(l)} $ is a one to one analytic
function in ${\mathbb{C}} \setminus
(\widetilde{\Delta}_{0}^{1}\cup\widetilde{\Delta}_{1}^{1})$.
\end{teo}

An expression for ${\mathcal{A}}^{(l)} $ will be given in (\ref{eq:lim6}) at the end of the proof of this result. The answer is given in terms of a conformal representation of an associated Riemann surface with $m_1+m_2+2$
sheets and genus zero onto the extended complex plane. It also depends on $l,\lambda_1,$ and $\lambda_2$. We wish to point out that if $\Lambda \subset \Lambda(\lambda_1,\lambda_2,C)$, taking $l_1 = \lambda_1(0)$ and $l_2 = \lambda_2(0)$ then ${\bf n}^{l} \in \Lambda(\lambda_1,\lambda_2,C+1)$, for all ${\bf n} \in \Lambda$. Therefore, for any sequence $\Lambda \subset \Lambda(\lambda_1,\lambda_2,C)$ we can always find at least one pair $(l_1;l_2),$ for which ratio asymptotics can be proved. In general, $(l_1;l_2)$ is admissible if $\lambda_1$ and $\lambda_2$ applied to ${\bf n}_1^{l_1}$ and ${\bf n}_2^{l_2}$, respectively, are also decreasing for all ${\bf n} \in \Lambda$ except at most a finite number of multi-indices.

For type II multiple orthogonal polynomials, with generating measures supported on intervals, and multi-indices in ${\mathbb{Z}}_+^{m+1}(\bullet)$, ratio asymptotics was proved in $\cite{AptLopRoc}$. Different extensions followed in \cite{LL} and \cite{FLLS}. Theorem \ref{teo:6} is an immediate consequence of Theorem \ref{teo:4} and \cite[Theorem 6.4]{FLLS}.

Let $M$ be the least common multiple of $m_1+1$ and $m_2+1$. Denote $\widetilde{\bf n} = (\widetilde{\bf n}_1;\widetilde{\bf n}_2)$ which is obtained adding $M/(m_1+1)$ to each component of ${\bf n}_1$ and $M/(m_2+1)$ to each component of ${\bf n}_2$. We have

\begin{cor}
\label{cor:3}
Let $S^1= {\mathcal{N}}'(\sigma^1_0,\ldots,\sigma^1_{m_1}), S^2=
{\mathcal{N}}'(\sigma^2_0,\ldots,\sigma^2_{m_2})$ be given. Let $\Lambda \subset {\mathbb{Z}}_+^{m_1 +1}\times {\mathbb{Z}}_+^{m_2 +1}$ be an infinite sequence of distinct multi-indices such that $|{\bf n}_1| = |{\bf n}_2| +1$ for all ${\bf n} \in \Lambda$ and
\[ \sup \{\max\{n_{i,0},\ldots,n_{i,m_i}\} - \min\{n_{i,0},\ldots,n_{i,m_i}\}: {\bf n} \in \Lambda \} < \infty, \qquad i=1,2.
\]
Then
\[
\lim_{{\bf n}\in\Lambda}\,\frac{{\mathcal{A}}_{\widetilde{\bf
n}}(z)}{{\mathcal{A}}_{{\bf n}}(z)}={\mathcal{A}}(z),\qquad  {\mathcal{K}}
\subset{\mathbb{C}}\setminus(\supp \sigma_0^1 \cup \mbox{\rm Co}(\supp \sigma_1^1)).
\]
\end{cor}

An expression for ${\mathcal{A}}$ appears in (\ref{eq:51}).

The strong asymptotics of type II multiple orthogonal polynomials for Nikishin systems was given by A. I. Aptekarev in \cite{Apt2} for diagonal sequences of multi-indices and systems of generating measures formed by weights satisfying Szeg\H{o}'s condition. It remains the best result in this respect. To conclude the Introduction, we call the reader's attention to the excellent survey by J. Nuttall \cite{Nut} on Hermite--Pad\'{e} polynomials. Here, in the form of a general conjecture, the author draws the general picture of the asymptotic behavior of Hermite--Pad\'{e} polynomials in terms of functions which are solutions of boundary value problems on associated Riemann surfaces. Our asymptotic results once more confirm his (at that time somewhat bold) predictions.

\section{Proof of Theorems \ref{teo:1}, \ref{teo:2} and Corollary \ref{cor:1}}

We begin with some auxiliary Lemmas.

\bl \label{reduc} Let $s_k,k=1,\ldots,m,$ be finite signed Borel
measures with compact support such that $ \mbox{\rm Co}(\supp s_k) = \Delta \subset
{\mathbb{R}}$. Let $F(z)= f_0(z) + \sum_{k=1}^m f_k(z)\widehat{
s}_k(z)  \in {\mathcal{H}}(\overline{\mathbb{C}} \setminus \Delta),$
where $f_k \in {\mathcal{H}}(V), k=0,\ldots,m,$ and $V$ is a
neighborhood of $\Delta$. If $F(z) = {\mathcal{O}}(1/z^2), z \to
\infty,$ then
\begin{equation} \label{eq:f}
\sum_{k=1}^m \int  f_k(x)d { s_k}(x) =0,
\end{equation}
whereas $F(z) = {\mathcal{O}}(1/z), z \to \infty,$ implies that
\begin{equation} \label{eq:g}
F(z) = \sum_{k=1}^m \int \frac{f_k(x)d { s_k}(x)}{z-x}.
\end{equation}
\el

{\bf Proof.} Let $\Gamma \subset V$ be a positively oriented closed
smooth Jordan curve that surrounds $\Delta$. If $F(z) =
{\mathcal{O}}(1/z^2), z \to \infty,$ from Cauchy's theorem, Fubini's
theorem and Cauchy's integral formula, it follows that
\[ 0 = \int_{\Gamma} F(z) dz = \sum_{k=1}^m \int_{\Gamma} f_k(z)\widehat{ s}_k(z) dz = \sum_{k=1}^m \int
\int_{\Gamma} \frac{f_k(z)dz}{z-x}ds_k(x) = 2\pi i \sum_{k=1}^m \int
f_k(x) ds_k(x),
\]
and we obtain (\ref{eq:f}). On the other hand, if $F(z) =
{\mathcal{O}}(1/z), z \to \infty,$  and we assume that $z$ is in the
unbounded connected component of the complement of $\Gamma$, Cauchy's
integral formula and Fubini's theorem render
\[ F(z) = \frac{1}{2\pi i} \int_{\Gamma} \frac{F(\zeta)d\zeta}{z-\zeta} = \frac{1}{2\pi i} \sum_{k=1}^m \int_{\Gamma}
\frac{f_k(\zeta)\widehat{s}_k(\zeta)  d\zeta}{z-\zeta} = \] \[
\sum_{k=1}^m \int \frac{1}{2\pi i} \int_{\Gamma} \frac{f_k(\zeta)
d\zeta}{(z-\zeta)(\zeta -x)} ds_k(x) = \sum_{k=1}^m \int
\frac{f_k(x)ds_k(x)}{z-x}
\]
which is (\ref{eq:g}). \fp

\bl
\label{lem:3}
Let $(s_{1,1},\ldots,s_{1,m}) =
{\mathcal{N}}(\sigma_1,\ldots,\sigma_m)$ and ${\bf n} \in
{\mathbb{Z}}_+^{m+1}$ be given. Consider the linear form
\[ {\mathcal{L}}_{\bf n}  = p_0 + \sum_{k=1}^m p_k \widehat{s}_{1,k}, \quad \deg p_k \leq n_k-1, \quad k=0,\ldots,m,
\]
where the polynomials $p_k$ have real coefficients. Assume that $n_0 = \max\{n_0,n_1-1,\ldots,n_m-1\}$. If
${\mathcal{L}}_{\bf n}$ had at least $|{\bf n}|$ zeros in $ {\mathbb{C}} \setminus \Delta_1$ the reduced form
$p_1 + \sum_{k=2}^m p_k \widehat{s}_{2,k}$ would have at least $|{\bf n}| - n_0$ zeros in $ {\mathbb{C}} \setminus \Delta_2$.
\el

{\bf Proof.} The function ${\mathcal{L}}_{\bf n}$ is symmetric with respect to the real line ${\mathcal{L}}_{\bf n}(\overline{z}) = \overline{{\mathcal{L}}_{\bf n}(z)}$; therefore, its zeros come in conjugate pairs. Thus, if ${\mathcal{L}}_{\bf n}$ has at least $|{\bf n}|$ zeros in $ {\mathbb{C}} \setminus \Delta_1$, there exists a polynomial $w_{\bf n}, \deg w_{\bf n}
\geq |{\bf n}|,$  with real coefficients and zeros contained in
${\mathbb{C}}\setminus \Delta_1$ such that ${\mathcal{L}}_{\bf
n}/w_{\bf n} \in {\mathcal{H}}({\mathbb{C}}\setminus \Delta_1)$. This function has a zero of order $ \geq |{\bf n}| - n_0 +1$
at $\infty$. Consequently, for all $ \nu
=0,\ldots,|{\bf n}| - n_0 -1,$
\[ \frac{z^{\nu} {\mathcal{L}}_{\bf n}}{w_{\bf n}} =
{\mathcal{O}}(1/z^2) \in {\mathcal{H}}(\overline{\mathbb{C}} \setminus
\Delta_1), \qquad z \to \infty, \qquad
\]
and
\[ \frac{z^{\nu} {\mathcal{L}}_{\bf n}}{w_{\bf n}} =  \frac{z^{\nu}p_0}{w_{\bf
n}} +  \sum_{k=1}^m \frac{z^{\nu}  p_k }{w_{\bf n}}
\widehat{s}_{1,k}\,.
\]
From (\ref{eq:f}), it follows that
\[ 0 = \int x^{\nu}  ( p_1 + \sum_{k=2}^m p_k
\widehat{s}_{2,k})(x) \frac{d \sigma_1(x)}{w_{\bf n}(x)}, \qquad \nu
=0,\ldots,|{\bf n}| - n_0 -1,
\]
taking into consideration that $s_{1,1} = \sigma_1$ and $ds_{1,k}(x) =
\widehat{s}_{2,k}(x) d\sigma_1(x), k=2,\ldots,m$.

These orthogonality relations imply that $p_1 + \sum_{k=2}^m p_k
\widehat{s}_{2,k}$ has at least $|{\bf n}| - n_0$ sign changes in the interior of $\Delta_1$.  In fact, if there were at most $|{\bf n}| - n_0 -1$ sign changes one can easily construct a polynomial $p$ of degree $\leq |{\bf n}| - n_0 -1$ such that $p(p_1 + \sum_{k=2}^m p_k
\widehat{s}_{2,k})$ does not change sign on $\Delta_1$ which contradicts the orthogonality relations. Therefore, already in the interior of $\Delta_1 \subset {\mathbb{C}} \setminus \Delta_2$, the reduced form would have the number of zeros claimed. \hfill $\Box$ \vspace{0,2cm}

Using induction, this lemma already allows to prove the AT property for multi-indices in ${\mathbb{Z}}_+^{m+1}(\circledast)$. That result is due to Driver and Stahl (see \cite[Theorem
2.4.1]{DrSt1}).

We reduce the general case to the one with $n_0 = \max\{n_0,n_1 -1,\ldots,n_m -1\}$ with

\bl \label{lem:4}
Let $(s_{1,1},\ldots,s_{1,m}) =
{\mathcal{N}}(\sigma_1,\ldots,\sigma_m), m \geq 1,$ and ${\bf n} \in
{\mathbb{Z}}_+^{m+1}$ be given. Consider the linear form ${\mathcal{L}}_{\bf n}$ defined in Lemma $\ref{lem:3}$. Assume that
$n_j = \max\{n_0+1,n_1,\ldots,n_m\}$. Then, there exist a Nikishin system $(s_{1,1}^*,\ldots,s_{1,m}^*) =
{\mathcal{N}}(\sigma_1^*,\ldots,\sigma_m^*)$, a multi-index ${\bf n}^* = (n_0^*,\ldots,n_m^*) \in {\mathbb{Z}}_+^{m+1}$
which is a permutation of ${\bf n}$ with $n_0^* = n_j$, and polynomials with real coefficients $p_k^*, \deg p_k^* \leq n_k^* -1, k=0,\ldots,m$, such that
\[ {\mathcal{L}}_{\bf n} = p_0 + \sum_{k=1}^m p_k \widehat{s}_{1,k} = (p_0^* + \sum_{k=1}^m p_k^* \widehat{s}_{1,k}^*)\widehat{s}_{1,j} = {\mathcal{L}}_{\bf n}^*\widehat{s}_{1,j}.
\]
\el

The proof is quite intricate and we leave it to the next section. Instead let us prove Theorem \ref{teo:1} assuming that the Lemma \ref{lem:4} is true. \vspace{0,2cm}

{\bf Proof of Theorem 1.} Obviously, the first statement of the theorem follows from the second. We prove the second one using induction on $m$. For $m=0$ the linear form reduces to a polynomial of degree $\leq n_0-1$ and thus has at most $n_0 -1$ zeros in the complex plane as claimed.

Assume that the result is true for any Nikishin system with $m-1 (\geq 0) $ measures and let us show that it is also valid for Nikishin systems with $m$ measures. To the contrary, let us suppose that ${\mathcal{L}}_{\bf n}$ has at least $|{\bf n}|$ zeros on ${\mathbb{C}}\setminus \Delta_1$.

Should $n_0 = \max\{n_0,n_1,\ldots,n_m\}$, by Lemma \ref{lem:3} the linear form $p_1 + \sum_{k=2}^m p_k \widehat{s}_{2,k}$ would have at least $|{\bf n}| - n_0$ zeros in ${\mathbb{C}} \setminus \Delta_2$. Now, $|{\bf n}| - n_0$ is the norm of the multi-index $(n_1,\ldots,n_m)$ which together with the Nikishin system ${\mathcal{N}}(\sigma_2,\ldots,\sigma_m)$ define the reduced form. This contradicts the induction hypothesis.

Suppose that $n_j = \max\{n_0+1,n_1,\ldots,n_m\}$. According to Lemma \ref{lem:4}, the linear form ${\mathcal{L}}_{\bf n}^* $ has the same zeros as ${\mathcal{L}}_{\bf n}$ in ${\mathbb{C}}\setminus \Delta_1$, since $s_{1,j}$ is never zero on that region. The multi-index ${\bf n}^*$ which determines ${\mathcal{L}}_{\bf n}^*$ has the same norm as ${\bf n}$  and its first component satisfies the assumptions of Lemma $\ref{lem:3}$. Following the same arguments as before we arrive to a contradiction. The proof is complete. \hfill $\Box$ \vspace{0,2cm}

For the proof of Theorem \ref{teo:2} and Corollary \ref{cor:1}, we also use

\bl \label{lem:orto} Let  $\widehat{\mathbb{S}}$ and ${\bf n} \in {\mathbb{Z}}_+^{m_1 +1}\times {\mathbb{Z}}_+^{m_2 +1}, |{\bf n}_1| = |{\bf n}_2| +1,$ be given. Then  ${\mathcal{A}}_{{\bf n}}$ satisfies
\begin{equation} \label{orto}
  \int {\mathcal{L}}_{{\bf n}_2}(x) {\mathcal{A}}_{{\bf n}}(x) d\sigma_0^2(x) = 0,
\end{equation}
for any linear form
\[ {\mathcal{L}}_{{\bf n}_2}(x) = p_0(x) + \sum_{j=1}^{m_2} p_j(x) \widehat{s}^2_{1,j}(x),
\]
where the $p_j, j=0,\ldots,m_2,$ denote arbitrary polynomials such that $\deg p_j \leq n_{2,j} -1$.    ${\mathcal{A}}_{{\bf n}}$ has exactly  $|{\bf n}_2|$ zeros in  ${\mathbb{C}} \setminus \mbox{\rm Co}(\supp \sigma_1^1)$, they are simple, and lie in the interior of $\mbox{\rm Co}(\supp \sigma_0^1)$.
\el

{\bf Proof}. In fact, from the condition b) of Definition \ref{defiSor} it follows that there exists a polynomial $d_{{\bf n},j}$ such that for any polynomial $p_j, \deg p_j \leq n_{2,j} -1, j\in \{0,\ldots,m_2\},$
\[ p_j(z) \left(\sum_{k=0}^{m_1} a_{{\bf n},k}(z) \int \frac{\widehat{s}^2_{1,j}(x)\widehat{s}^1_{1,k}(x)d\sigma^2_0(x)}{z-x}  - d_{{\bf
n},j}(z)\right) = {\mathcal{O}}\left(1/z^2\right), \qquad z \to
\infty,
\]
(here $\widehat{s}^2_{1,0}  \equiv 1)$ and the function on the left
hand side is holomorphic in $\overline{\mathbb{C}} \setminus
\mbox{Co}(\supp \sigma^2_{0})$. Using Lemma \ref{reduc}, it follows that
\[
 \int p_j(x) \widehat{s}^2_{1,j}(x)\sum_{k=0}^{m_1} a_{{\bf n},k}(x)\widehat{s}^1_{1,k}(x) d \sigma^2_{0}(x)  =0.
\]
Adding these relations for $j=0,\ldots,m_2$, we obtain (\ref{orto}).

From Theorem \ref{teo:1}, we know that ${\mathcal{A}}_{\bf n}$ has at most $|{\bf n}_1| -1 = |{\bf n}_2|$ zeros on
${\mathbb{C}} \setminus \mbox{\rm Co}(\supp \sigma_1^1)$. From (\ref{orto}) it follows that this form has at least $|{\bf n}_2|$ sign changes in the interior of $\mbox{\rm Co}(\supp \sigma_0^2) = \mbox{\rm Co}(\supp \sigma_0^1)$. Therefore, the last statement is obtained.
\fp \vspace{0,2cm}

Let us prove Theorem \ref{teo:2} assuming that Lemma \ref{lem:4} (Theorem \ref{teo:1}) is true.
\vspace{0,2cm}

{\bf Proof of Theorem \ref{teo:2}.} Suppose that for some ${\bf n}, {\mathbb{A}}_{\bf n}$ is not normal. That is,  some component $a_{{\bf n},k}$ of ${\mathbb{A}}_{\bf n}$ has $\deg a_{{\bf n},k} \leq n_{1,k} -2$. According to Theorem \ref{teo:1}, ${\mathcal{A}}_{\bf n}$ can have on the interval $\mbox{Co}(\supp \sigma^2_0)$ at most $|{\bf n}_1| -2 =  |{\bf n}_2| -1$ zeros. Consequently, this function can have in the interior of $\mbox{Co}(\supp \sigma^2_0)$ at most $N \leq |{\bf n}_2| -1$ sign changes. Suppose this is the case and let $x_1,\ldots,x_N$ be the points where it changes sign. According to Theorem \ref{teo:1}, $(1,\widehat{s}^2_{1,1},\ldots,\widehat{s}^2_{1,m_2})$ is also an AT system. Using the properties of Tchebyshev systems, we can find polynomials $p_0,\ldots,p_{m_2},$  with $\deg p_j \leq n_{2,j} -1,$ such that ${\mathcal{L}}_{\bf n}(x) = p_0(x) + \sum_{j=1}^{m_2} p_{j}(x) \widehat{s}_{1,j}^2(x)$
changes  sign at $x_1,\ldots,x_N,$ and has no other points where it changes sign in the interior of $\mbox{Co}(\supp \sigma^2_0)$. Therefore, the function ${\mathcal{L}}_{\bf n}(x) {\mathcal{A}}_{{\bf n}}(x)$ has constant sign on $\mbox{Co}(\supp \sigma^2_0)$ but this contradicts (\ref{orto}) since $\sigma^2_0$ is a measure with constant sign whose support contains infinitely many points.   Thus, $\deg a_{{\bf n},k} = n_{1,k} -1, k=0,\ldots,m_1$, and perfectness has been established.

Let us assume that there are two non collinear solutions ${\mathbb{A}}_{\bf n}, {\mathbb{A}}_{\bf n}^*,$ to a)-b).   Then, there exists a real constant $C \neq 0$ such that ${\mathbb{A}}_{\bf n} - C {\mathbb{A}}_{\bf n}^* \not \equiv 0$ and at least one of the components  of ${\mathbb{A}}_{\bf n} - C {\mathbb{A}}_{\bf n}^* $ satisfies $\deg (a_{{\bf n},k } - C a_{{\bf n},k }^*) \leq n_{1,k} -2.$ This is not possible since ${\mathbb{A}}_{\bf n} - C{\mathbb{A}}_{\bf n}^*$ also solves a)-b) and according to what was proved above all its components must have maximum possible degree. \fp

\begin{defi}
\label{def:hausdroff} Let $E$ be a subset of the complex plane and $\mathcal{U}$ the class of all coverings of $E$ by disks $U_n$. The radius of $U_n$ is denoted $|U_n|$. The (one dimensional) Hausdorff content of $E$ is
\[ h(E) = \inf \{\sum|U_n|: \{U_n\} \in {\mathcal{U}}\}.
\]
\end{defi}

Let $\{f_n\}_{n \in \Lambda}$ be a sequence of functions defined on a region $D \subset {\mathbb{C}}$. We say that $\{f_n\}_{n \in \Lambda}$ converges to $f$ in Hausdorff content on $D$ if for every compact set ${\mathcal{K}}  \subset D$ and any $\varepsilon > 0$
\[ \lim_{n \in \Lambda} h(\{z \in {\mathcal{K}}: |f_n(z) - f(z)| > \varepsilon\}) =0.
\]
We denote this by
\[ {\mathcal{H}}-\lim_{n \to \infty} f_n = f, \qquad {\mathcal{K}} \subset D.
\]

In \cite[Lemma 1]{Gon}, A.A. Gonchar proved that if the functions $f_n$ are holomorphic in $D$ and they converge in Hausdoff content to $f$ in $D$, then $f$ is in fact holomorphic in $D$ (more precisely, differs from a holomorphic function on a set of zero Hausdorff content) and the convergence (to the equivalent holomorphic function) is  uniform on each compact subset of $D$.
\vspace{0.2cm}

{\bf Proof of Corollary \ref{cor:1}.} In \cite[Theorem 1]{LF3} it was proved that under the assumptions of the corollary, for each $k=0,\ldots,m,$
\[ {\mathcal{H}}-\lim_{{\bf n} \in \Lambda} R_{{\bf n},k} = \widehat{s}_k, \qquad {\mathcal{K}} \subset \overline{\mathbb{C}} \setminus \mbox{Co}(\supp \sigma_0).
\]
Due to Gonchar's lemma and the last assertion of Lemma \ref{lem:orto}, it follows that convergence is uniform on each compact subset of $\overline{\mathbb{C}} \setminus \mbox{Co} (\supp \sigma_0)$. Regarding the proof of the rate of convergence, we refer to the last sentence on page 104 of \cite{LF3} (see also \cite[Corollary 1]{LF3}). \hfill $\Box$

\section{Proof of Lemma \ref{lem:4} and Corollary \ref{cor:2}.}

It is well known (see appendix in \cite{KN} and \cite[Theorem 6.3.5]{stto}) that for each $s
\in {\mathcal{M}}(\Delta),$ there exists a measure $\tau \in
{\mathcal{M}}(\Delta)$ and ${\ell}(z)=a z+b, a = 1/|s|, b \in {\mathbb{R}},$ such that
\begin{equation} \label{s22}
{1}/{\widehat{s}(z)}={\ell}(z)+ \widehat{\tau}(z),
\end{equation}
where $|s|$ is the total variation of the measure $s.$  For
convenience, we call $\tau$  the inverse measure of $s.$ Such
measures will appear frequently in our reasonings, so we will fix a
notation to distinguish them. They will always refer to inverses of
measures denoted with $s$ and will carry over to them the
corresponding sub-indices. The same goes for the  polynomials
$\ell$. For instance, if $s_{\alpha,\beta} = \langle
\sigma_{\alpha},\sigma_{\beta} \rangle $
\[
{1}/{\widehat{s}_{\alpha,\beta}(z)}  ={\ell}_{\alpha,\beta}(z)+
\widehat{\tau}_{\alpha,\beta}(z).
\]
For convenience, sometimes we write $\langle
\sigma_{\alpha},\sigma_{\beta} \widehat{\rangle}$ in place of $\widehat{s}_{\alpha,\beta}$. This is specially useful later on where we need the Cauchy transforms of complicated expressions of products of measures for which we do not have a short hand notation. Since $s_{\alpha,\alpha} = \sigma_{\alpha}$, we also write
\[
{1}/{\widehat{\sigma}_{\alpha}(z)} ={\ell}_{\alpha,\alpha}(z)+
\widehat{\tau}_{\alpha,\alpha}(z).
\]

\bl \label{alfabeta} Let $\sigma_{\alpha} \in
{\mathcal{M}}(\Delta_{\alpha}), \sigma_{\beta} \in
{\mathcal{M}}(\Delta_{\beta}),$ and $\Delta_{\alpha}
\cap \Delta_{\beta} =\emptyset.$ Then:

\begin{equation} \label{2.1}\widehat{\sigma}_{\alpha} (z)\widehat{\sigma}_{\beta}(z)
=\langle \sigma_{\alpha},\sigma_{\beta} \widehat{\rangle}(z)+\langle
\sigma_{\beta},{\sigma}_{\alpha}\widehat{\ra}(z), \quad z \in
{\mathbb{C}} \setminus \left(\Delta_{\alpha} \cup
\Delta_{\beta}\right),
\end{equation}

\begin{equation} \label{2.2}
\frac{\widehat{\sigma}_{\alpha}(z)}{\langle
\sigma_{\alpha},\sigma_{\beta}\widehat{\rangle}(z)}=
\frac{|\sigma_{\alpha}|}{|\langle\sigma_{\alpha}\sigma_{\beta}\rangle|}+
\int \frac{\langle
\sigma_{\beta},{\sigma}_{\alpha}\widehat{\ra}(x_{\alpha})}{
\widehat{\sigma}_{\beta} (x_{\alpha})}\frac{d \tau_{\alpha,\beta}
(x_{\alpha})}{z-x_{\alpha}} =
\frac{|\sigma_{\alpha}|}{|\langle\sigma_{\alpha},\sigma_{\beta}\rangle|}+\langle
 \frac{{\tau}_{\alpha,\beta}}{\widehat{\sigma}_{\beta}},\sigma_{\beta}, {\sigma}_{\alpha}\widehat{\ra}(z),
\end{equation}
\begin{equation} \label{2.3}
\frac{\langle
\sigma_{\alpha},\sigma_{\beta}\widehat{\rangle}(z)}{\widehat
{\sigma}_{\alpha} (z)}=\frac{|\langle
\sigma_{\alpha},\sigma_{\beta}\rangle|}{|\sigma_{\alpha}|}- \int
\frac{{\langle \sigma_{\beta},
{\sigma}_{\alpha}\widehat{\ra}(x_{\alpha})}d \tau_{\alpha,\alpha}
(x_{\alpha})}{z-x_{\alpha}}=\frac{|\langle
\sigma_{\alpha},\sigma_{\beta}\rangle|}{|\sigma_{\alpha}|} -\langle
 {\tau}_{\alpha,\alpha},\sigma_{\beta}, {\sigma}_{\alpha}\widehat{\ra}(z).
\end{equation}
\el

{\bf Proof.} In fact, (\ref{2.1}) follows from the  chain of equalities
\[
\widehat{\sigma}_{\alpha} (z)\widehat{\sigma}_{\beta} (z)=\int \int
\frac{d \sigma_{\alpha}(x_{\alpha}) d
\sigma_{\beta}(x_{\beta})}{(z-x_{\alpha})(z-x_{\beta})}=\int \int
\left(\frac{1}{z-x_{\alpha}}-\frac{1}{z-x_{\beta}}\right)\frac{d
\sigma_{\alpha}(x_{\alpha}) d
\sigma_{\beta}(x_{\beta})}{x_{\alpha}-x_{\beta}}
\]
\[
\qquad =\int \widehat{\sigma}_{\alpha} (x_{\beta}) \frac{d
\sigma_{\beta}(x_{\beta})}{z-x_{\beta}} + \int
\widehat{\sigma}_{\beta} (x_{\alpha}) \frac{d
\sigma_{\alpha}(x_{\alpha})}{z -x_{\alpha}}=\langle
\sigma_{\alpha},\sigma_{\beta} \widehat{\rangle}(z)+\langle
\sigma_{\beta}, {\sigma}_{\alpha}\widehat{\ra}(z).
\]

Notice that
\[
\frac{\widehat{\sigma}_{\alpha}(z)}{\langle
\sigma_{\alpha},\sigma_{\beta}\widehat{\rangle}(z)}-
\frac{|\sigma_{\alpha}|}{|\langle\sigma_{\alpha}\sigma_{\beta}\rangle|}={\mathcal{O}}\left(\frac{1}{z}\right)
\in {\mathcal{H}} \left(\overline{\mathbb{C}} \setminus
\Delta_{\alpha}\right), \qquad z \to \infty.
\]
From (\ref{s22}) and (\ref{2.1}), it follows that
\[
\frac{\widehat{\sigma}_{\alpha}(z)}{\langle
\sigma_{\alpha},\sigma_{\beta}\widehat{\rangle}(z)} = \frac{\widehat
{\sigma}_{\beta} (z) \widehat{\sigma}_{\alpha} (z)}{\widehat{
\sigma}_{\beta}(z)\langle
\sigma_{\alpha},\sigma_{\beta}\widehat{\rangle}(z)}  =\frac{\langle
\sigma_{\alpha},\sigma_{\beta}\widehat{\rangle}(z)+\langle
\sigma_{\beta}, {\sigma}_{\alpha}\widehat{\ra}(z)}{\widehat
{\sigma}_{\beta} (z)\langle
\sigma_{\alpha},\sigma_{\beta}\widehat{\rangle}(z)} =
\]
\[
\frac{1}{\widehat {\sigma}_{\beta} (z)} +  \frac{\langle
\sigma_{\beta}, {\sigma}_{\alpha}\widehat{\ra}(z)}{\widehat
{\sigma}_{\beta} (z)}\ell_{\alpha,\beta} + \frac{\langle
\sigma_{\beta}, {\sigma}_{\alpha}\widehat{\ra}(z)}{\widehat
{\sigma}_{\beta} (z)} \widehat{\tau}_{\alpha,\beta} (z).
\]
Since $-
\frac{|\sigma_{\alpha}|}{|\langle\sigma_{\alpha}\sigma_{\beta}\rangle|}
+ \frac{1}{\widehat {\sigma}_{\beta}  } +  \frac{\langle
\sigma_{\beta}, {\sigma}_{\alpha}\widehat{\ra} }{\widehat
{\sigma}_{\beta} }\ell_{\alpha,\beta} $ and $ \frac{\langle
\sigma_{\beta}, {\sigma}_{\alpha}\widehat{\ra} }{\widehat
{\sigma}_{\beta} }$ are analytic on a neighborhood of the interval
$\Delta_{\alpha}$, which contains the support of
$\tau_{\alpha,\beta}$, relation (\ref{eq:g}) implies (\ref{2.2}).

The proof of (\ref{2.3}) is similar but somewhat more direct. Again,
we have that
\[
\frac{\langle\sigma_{\alpha},\sigma_{\beta}\widehat{\rangle}(z)}{
\widehat{\sigma}_\alpha(z)}-\frac{|\langle
\sigma_{\alpha},\sigma_{\beta}\rangle|}{|\sigma_{\alpha}|}={\mathcal{O}}\left(\frac{1}{z}\right)
\in {\mathcal{H }}\left(\overline{\mathbb{C}} \setminus
\Delta_{\alpha} \right), \qquad z \to \infty.
\]
From (\ref{s22}) and (\ref{2.1}), we get
\[
\frac{\langle
\sigma_{\alpha},\sigma_{\beta}\widehat{\rangle}(z)}{\widehat{
\sigma}_{\alpha}(z)} =\frac{\widehat{\sigma}_{\alpha}
(z)\widehat{\sigma}_{\beta} (z)-\langle \sigma_{\beta},
{\sigma}_{\alpha}\widehat{\ra}(z)}{\widehat{ \sigma}_{\alpha}(z)} =
\widehat{\sigma}_{\beta} (z) - \langle \sigma_{\beta},
{\sigma}_{\alpha}\widehat{\ra}(z)\ell_{\alpha,\alpha}(z) - \langle
\sigma_{\beta},
{\sigma}_{\alpha}\widehat{\ra}(z)\widehat{\tau}_{\alpha,\alpha}(z).
\]
But $-\frac{|\langle
\sigma_{\alpha},\sigma_{\beta}\rangle|}{|\sigma_{\alpha}|} +
\widehat{\sigma}_{\beta}   - \langle \sigma_{\beta},
{\sigma}_{\alpha}\widehat{\ra} \ell_{\alpha,\alpha} $ and $\langle
\sigma_{\beta}, {\sigma}_{\alpha}\widehat{\ra}$ are analytic in a
neighborhood of $\Delta_{\alpha}$; therefore, (\ref{eq:g}) implies
(\ref{2.3}).
 \hfill $\Box$ \vspace{0,2cm}

Formulas (\ref{s22})-(\ref{2.1}) are the building blocks for (\ref{2.2})-(\ref{2.3}) and many more interesting relations.
Let us further extend Lemma \ref{alfabeta}. The new formulas
may be grouped in two since (\ref{4.4}) may be regarded a special
case of (\ref{4.5}) and (\ref{4.6})-(\ref{4.7}) as special cases of
(\ref{4.2}). Putting each group in one formula causes some
notational incongruence which we prefer to avoid for the benefit of
the reader.

\bl \label{cocientes} Let $(s_{1,1},\ldots,s_{1,m}) =
{\mathcal{N}}(\sigma_1,\ldots,\sigma_m)$ be given. Then:
\begin{equation} \label{4.4}
\frac{\widehat{s}_{1,k}}{\widehat{s}_{1,1}} =
\frac{|s_{1,k}|}{|s_{1,1}|} - \la \tau_{1,1},\la s_{2,k},\sigma_1
\ra \widehat{\ra}  , \qquad  1=j < k \leq m,
\end{equation}
\begin{equation} \label{4.5} \frac{\widehat{s}_{1,k}}{\widehat{s}_{1,j}} =
\frac{|s_{1,k}|}{|s_{1,j}|} + (-1)^j \la
\tau_{1,j},\la\tau_{2,j},s_{1,j} \ra,\ldots,
\la\tau_{j,j},s_{j-1,j}\ra,\la s_{j+1,k},\sigma_j\ra \widehat{\ra},
\quad 2\leq j < k \leq m,
\end{equation}
\begin{equation} \label{4.6}
\frac{\widehat{s}_{1,1}}{\widehat{s}_{1,j}} =
\frac{|s_{1,1}|}{|s_{1,j}|} + \la
\frac{\tau_{1,j}}{\widehat{s}_{2,j}}, \la s_{2,j},\sigma_1\ra
\widehat{\ra}=
\end{equation}
\[ \frac{|s_{1,1}|}{|s_{1,j}|} + \frac{|\la
s_{2,j},\sigma_1\ra|}{|s_{2,j}|} \widehat{\tau}_{1,j} - \la
\tau_{1,j},\la \tau_{2,j},s_{1,j}\ra \widehat{\ra}, \qquad 1 =k <
j\leq m,
\]
\begin{equation} \label{4.7}
\frac{\widehat{s}_{1,2}}{\widehat{s}_{1,j}} =
\frac{|s_{1,2}|}{|s_{1,j}|} - \la \tau_{1,j}, \frac{\la
\tau_{2,j},s_{1,j}\ra}{\widehat{s}_{3,j}}, \la s_{3,j},\sigma_2\ra
\widehat{\ra}=
\end{equation}
\[ \frac{|s_{1,2}|}{|s_{1,j}|} -
\frac{|\la s_{3,j},\sigma_2\ra|}{|s_{3,j}|}  \la {\tau}_{1,j},\la
\tau_{2,j},s_{1,j}\ra \widehat{\ra} + \la
\tau_{1,j},\la \tau_{2,j},s_{1,j}\ra , \la \tau_{3,j},s_{2,j}\ra
\widehat{\ra},\qquad 2 =k < j \leq m,
\]
\begin{equation} \label{4.2}
\frac{\widehat{s}_{1,k}}{\widehat{s}_{1,j}} =
\frac{|s_{1,k}|}{|s_{1,j}|} + (-1)^{k-1} \la
\tau_{1,j},\la\tau_{2,j},s_{1,j} \ra,\ldots,
\la\tau_{k-1,j},s_{k-2,j}\ra, \frac{\la \tau_{k,j},
{s}_{k-1,j}\ra}{\widehat{s}_{k+1,j}},\la s_{k+1,j},\sigma_k\ra
\widehat{\ra} =
\end{equation}
\[ \frac{|s_{1,k}|}{|s_{1,j}|} + (-1)^{k-1} \frac{|\la s_{k+1,j},\sigma_k\ra|}{|s_{k+1,j}|}\la
\tau_{1,j},\la\tau_{2,j},s_{1,j} \ra,\ldots,
\la\tau_{k-1,j},s_{k-2,j}\ra, \la\tau_{k,j},s_{k-1,j}\ra
\widehat{\ra} \,+
\]
\[ (-1)^{k}  \la \tau_{1,j},\la\tau_{2,j},s_{1,j} \ra,\ldots,
 \la\tau_{k,j},s_{k-1,j}\ra, \la\tau_{k+1,j},s_{k,j}\ra
\widehat{\ra}\,. \qquad 3=k < j \leq m.
\]
\el

{\bf Proof.} Cauchy transforms equal zero at infinity; therefore,
the constants appearing on the right hand sides in each of the first
equalities of (\ref{4.4})-(\ref{4.2}) must be as indicated, if in
fact the other term is a Cauchy transform. Consequently, we will not
pay attention to the constants coming out of the consecutive
transformations we make in our deduction and simply denote them with
consecutive constants $C_j$.

Obviously,  (\ref{4.4}) is deduced from (\ref{2.3}) taking
$\sigma_{\alpha} = \sigma_1 = s_{1,1} $ and $\sigma_{\beta} = \la
\sigma_2,\cdots,\sigma_k \ra = s_{2,k}$. Formula (\ref{4.5}) is
obtained applying (\ref{2.3}) inside out several times  as we will
indicate.

Let $2 \leq j < k \leq m$. Using (\ref{2.3}) on $
{\widehat{s}_{j,k}}/{\widehat{s}_{j,j}}$, we have that
\begin{equation} \label{4.3}
\la \sigma_{j-1},\sigma_j,\ldots,\sigma_k\widehat{\ra} = \la
\frac{s_{j-1,j}}{\widehat{s}_{j,j}}, s_{j,k}\widehat{\ra} =   \la
\frac{\widehat{s}_{j,k}}{\widehat{s}_{j,j}}s_{j-1,j}\widehat{\ra} =
C_1 \widehat{s}_{j-1,j} - \la s_{j-1,j},\tau_{j,j}, \la
s_{j+1,k},\sigma_j\ra \widehat{\ra}.
\end{equation}
In particular, if $j=2$ we get
\[ \la \sigma_{1},\sigma_2,\ldots,\sigma_k\widehat{\ra} = \la \frac{s_{1,2}}{\widehat{s}_{2,2}},
s_{2,k}\widehat{\ra} =   \la
\frac{\widehat{s}_{2,k}}{\widehat{s}_{2,2}}s_{1,2}\widehat{\ra} =
C_1 \widehat{s}_{1,2} - \la s_{1,2},\tau_{2,2}, \la
s_{3,k},\sigma_2\ra \widehat{\ra},
\]
and applying (\ref{2.3}) on ${\la s_{1,2},\tau_{2,2}, \la
s_{3,k},\sigma_2\ra \widehat{\ra}}/{\widehat{s}_{1,2}}$, it follows
that
\[ \frac{\widehat{s}_{1,k}}{\widehat{s}_{1,2}} = C_1 - \frac{1}{\widehat{s}_{1,2}} \la s_{1,2},\tau_{2,2}, \la
s_{3,k},\sigma_2\ra \widehat{\ra} = \frac{|s_{1,k}|}{|s_{1,2}|} +
\la \tau_{1,2},\la \tau_{2,2}, s_{1,2} \ra, \la s_{3,k},\sigma_2\ra
\widehat{\ra}
\]
which is (\ref{4.5}) for $j=2$.

Assume that $j \geq 3$. We write then
\[ \widehat{s}_{1,k} = \la \sigma_1,\ldots,\sigma_{j-2},\frac{\widehat{s}_{j,k}}{\widehat{s}_{j,j}}s_{j-1,j}\widehat{\ra}
\]
and on account of (\ref{4.3}), we obtain
\[ \frac{\widehat{s}_{1,k}}{\widehat{s}_{1,j}} = C_1 - \frac{1}{\widehat{s}_{1,j}} \la
\sigma_1,\ldots,\sigma_{j-2}, s_{j-1,j},\tau_{j,j}, \la
s_{j+1,k},\sigma_j \ra \widehat{\ra}.
\]
This means that
\[ \frac{\widehat{s}_{1,k}}{\widehat{s}_{1,3}} = C_1 - \frac{1}{\widehat{s}_{1,3}} \la
 \frac{s_{1,3}}{\widehat{s}_{2,3}},s_{2,3},\tau_{3,3}, \la
s_{4,k},\sigma_3 \ra \widehat{\ra}, \qquad j=3,
\]
or
\[ \frac{\widehat{s}_{1,k}}{\widehat{s}_{1,j}} = C_1 - \frac{1}{\widehat{s}_{1,j}} \la
\sigma_1,\ldots,\sigma_{j-3},\frac{s_{j-2,j}}{\widehat{s}_{j-1,j}},s_{j-1,j},\tau_{j,j},
\la s_{j+1,k},\sigma_j \ra \widehat{\ra}, \qquad j \geq 4.
\]
Using (\ref{2.3}) again, it follows that
\[
 \frac{\la s_{j-1,j},\tau_{j,j}, s_{j+1,k},\sigma_j \widehat{\ra}}
{\widehat{s}_{j-1,j}}   = C_2  - \la \tau_{j-1,j}, \la
\tau_{j,j},s_{j-1,j}\ra, \la s_{j+1,k},\sigma_j\ra \widehat{\ra}.
\]
Substituting above, we have
\[ \frac{\widehat{s}_{1,k}}{\widehat{s}_{1,3}} = C_3 + (-1)^2 \frac{1}{\widehat{s}_{1,3}} \la
s_{1,3}, \tau_{2,3}, \la \tau_{3,3},s_{2,3}\ra, \la s_{4,k},\sigma_3
\ra \widehat{\ra}, \qquad j=3,
\]
or
\[ \frac{\widehat{s}_{1,k}}{\widehat{s}_{1,j}} = C_3 + (-1)^2 \frac{1}{\widehat{s}_{1,j}} \la
\sigma_1,\ldots,\sigma_{j-3}, s_{j-2,j}, \tau_{j-1,j}, \la
\tau_{j,j},s_{j-1,j}\ra, \la s_{j+1,k},\sigma_j \ra \widehat{\ra},
\qquad j \geq 4.
\]
If $j=3$ one more use of (\ref{2.3}) brings us to (\ref{4.5}). If $j
\geq 4$ we keep on applying (\ref{2.3}) inside out until we arrive
at
\[ \frac{\widehat{s}_{1,k}}{\widehat{s}_{1,j}} = C_4 + (-1)^{j-1} \frac{1}{\widehat{s}_{1,j}}
\la s_{1,j},\tau_{2,j},\la \tau_{3,j},s_{2,j}\ra, \ldots,
\la\tau_{j,j},s_{j-1,j}\ra,\la s_{j+1,k},\sigma_j\ra \widehat{\ra},
\]
which is just one step away  from  (\ref{4.5}) through  (\ref{2.3})
taking
\[\sigma_{\alpha} = s_{1,j}, \qquad \sigma_{\beta} = \la
\tau_{2,j},\la \tau_{3,j},s_{2,j}\ra, \ldots,
\la\tau_{j,j},s_{j-1,j}\ra,\la s_{j+1,k},\sigma_j\ra \ra  .
\]

Now, let us prove formulas (\ref{4.6})-(\ref{4.2}). The second
equality in each one of these relations is an immediate consequence
of (\ref{2.3}) since from it we get
\begin{equation} \label{4.8}
\frac{\la  {s}_{k+1,j},\sigma_k \widehat{\ra}}{\widehat{s}_{k+1,j}}
= \frac{|\la
 {s}_{k+1,j},\sigma_k  \ra |}{|{s}_{k+1,j}|} - \la \tau_{k+1,j},
 s_{k,j} \widehat{\ra}.
\end{equation}
The proof of the first equality is obtained, generally speaking, as
in proving (\ref{4.5}) except that we begin using once relation
(\ref{2.2}). In fact, when $k=1$ formula (\ref{4.6}) follows
directly from (\ref{2.2}) taking $\sigma_{\alpha}= \sigma_1 =
s_{1,1}$ and $\sigma_{\beta} = s_{2,j} $.

 Assume that $2 \leq k < j \leq m$. Using
(\ref{2.2}), it follows that
\[ \la \sigma_{k-1}, \sigma_k\widehat{\ra} = \la \frac{s_{k-1,j}}{\widehat{s}_{k,j}},
s_{k,k}\widehat{\ra} =   \la
\frac{\widehat{s}_{k,k}}{\widehat{s}_{k,j}}s_{k-1,j}\widehat{\ra} =
C_5 \widehat{s}_{k-1,j} + \la
s_{k-1,j},\frac{\tau_{k,j}}{\widehat{s}_{k+1,j}}, \la
s_{k+1,j},\sigma_{k} \ra \widehat{\ra}.
\]
Consequently,
\[ \frac{\widehat{s}_{1,2}}{\widehat{s}_{1,j}} =  C_5 + \frac{1}{\widehat{s}_{1,j}} \la
s_{1,j},\frac{\tau_{2,j}}{\widehat{s}_{3,j}}, \la s_{3,j},\sigma_2
\ra \widehat{\ra}, \qquad k =2,
\]
or
\[ \frac{\widehat{s}_{1,k}}{\widehat{s}_{1,j}} =  C_5 + \frac{1}{\widehat{s}_{1,j}} \la
\sigma_1,\ldots,\sigma_{k-2},
s_{k-1,j},\frac{\tau_{k,j}}{\widehat{s}_{k+1,j}}, \la
s_{k+1,j},\sigma_k \ra \widehat{\ra}, \qquad k \geq 3.
\]
From this point on we use (\ref{2.3}). From this formula, we obtain
\[ \frac{\la s_{k-1,j},\frac{\tau_{k,j}}{\widehat{s}_{k+1,j}},
\la s_{k+1,j},\sigma_k \ra \widehat{\ra}}{\widehat{s}_{k-1,j}} = C_6
- \la \tau_{k-1,j},\frac{\la
\tau_{k,j},s_{k-1,j}\ra}{\widehat{s}_{k+1,j}}, \la
s_{k+1,j},\sigma_k \ra \widehat{\ra},
\]
and (\ref{4.7}) readily follows if $k=2$. For $k \geq 3$ this
implies
\[ \frac{\widehat{s}_{1,3}}{\widehat{s}_{1,j}} =  C_7 - \frac{1}{\widehat{s}_{1,j}}
\la  s_{1,j},\tau_{2,j},\frac{\la
\tau_{3,j},s_{2,j}\ra}{\widehat{s}_{4,j}}, \la s_{4,j},\sigma_3 \ra
\widehat{\ra}, \qquad k =3,
\]
or
\[ \frac{\widehat{s}_{1,k}}{\widehat{s}_{1,j}} =  C_7 - \frac{1}{\widehat{s}_{1,j}}
\la \sigma_1,\ldots,\sigma_{k-3},s_{k-2,j},\tau_{k-1,j},\frac{\la
\tau_{k,j},s_{k-1,j}\ra}{\widehat{s}_{k+1,j}}, \la
s_{k+1,j},\sigma_k \ra \widehat{\ra}, \qquad k \geq 4.
\]
Continuing down, using (\ref{2.3}) on each step, we obtain
\[ \frac{\widehat{s}_{1,k}}{\widehat{s}_{1,j}} = C_8 + (-1)^{k-2} \frac{1}{\widehat{s}_{1,j}}
\la s_{1,j},\tau_{2,j},\la \tau_{3,j},s_{2,j}\ra, \ldots,  \la
\tau_{k-1,j},s_{k-2,j}\ra,\frac{\la
\tau_{k,j},s_{k-1,j}\ra}{\widehat{s}_{k+1,j}}, \la
s_{k+1,j},\sigma_k \ra \widehat{\ra}.
\]
One more use of (\ref{2.3}) with
\[ \sigma_{\alpha} = s_{1,j}, \qquad \sigma_{\beta} = \la \tau_{2,j},\la \tau_{3,j},s_{2,j}\ra, \ldots,  \la
\tau_{k-1,j},s_{k-2,j}\ra,\frac{\la
\tau_{k,j},s_{k-1,j}\ra}{\widehat{s}_{k+1,j}}, \la
s_{k+1,j},\sigma_k \ra  {\ra}
\]
gives the first equality of
(\ref{4.2}). With this we conclude the proof. \hfill $\Box$

\begin{rem} We wish to point out that formulas
(\ref{4.4})-(\ref{4.5}) and the second equalities in
(\ref{4.6})-(\ref{4.2}) are contained in \cite[Theorem
3.1.3]{DrSt2}, where they were deduced using the Stieltjes-Plemelj inversion formula. Our explicit
expressions of the right hand sides are necessary for the arguments
to follow. Additionally, the first equalities in
(\ref{4.6})-(\ref{4.2}) are of great value for the proof of the
general case.
\end{rem}

{\bf Proof of Lemma \ref{lem:4} when $j=1$.}
From (\ref{s22}) and (\ref{4.4}), we have
\[ \frac{ {\mathcal{L}}_{\bf n}}{\widehat{s}_{1,1} } =
\frac{ p_0}{ \widehat{s}_{1,1}}+p_1   + \sum_{k=2}^m p_k\frac{\widehat{s}_{1,k}}{\widehat{s}_{1,1}} =
\]
\[
 (\ell_{1,1}p_0  + p_1 + \sum_{k=2}^m \frac{|s_{1,k}|}{|s_{1,1}|}p_k) +  p_0\widehat{\tau}_{1,1}
- \sum_{k=2}^m p_k \la  {\tau}_{1,1},
s_{2,k},\sigma_1 \widehat{\ra} = {\mathcal{L}}_{\bf n}^*.
\]
We are done taking ${\bf n}^* = (n_1,n_0,n_2,\ldots,n_m)$ and
\[{\mathcal{N}}(\sigma_1^*,\ldots,\sigma_m^*) = {\mathcal{N}}(\tau_{1,1},\la \sigma_2,\sigma_1\ra,\sigma_3,\ldots,\sigma_m)
\]
since $\la s_{2,k},\sigma_1\ra = \la \la \sigma_2,\sigma_1\ra ,\sigma_3,\ldots,\sigma_k\ra$ when $k\geq 3$.  \hfill $\Box$ \vspace{0,2cm}

In the sequel $2 \leq j \leq m$. From (\ref{s22}),  (\ref{4.5}), and the first equalities in
(\ref{4.6})-(\ref{4.2}), one has
\[ \frac{ {\mathcal{L}}_{\bf n}}{\widehat{s}_{1,j}} =
\frac{p_0}{\widehat{s}_{1,j}}+
 p_j   + \sum_{k\neq j,k=1}^m p_k  \frac{\widehat{s}_{1,k}}{\widehat{s}_{1,j}} =
 (\ell_{1,j}p_0  +  p_j  + \sum_{k\neq j, k=1}^m \frac{|s_{1,k}|}{|s_{1,j}|}
p_k)  +
\]
\[p_0 \widehat{\tau}_{1,j}
+ p_1\la
\frac{\tau_{1,j}}{\widehat{s}_{2,j}}, \la s_{2,j},\sigma_1 \ra
\widehat{\ra} +
\]
\[ \sum_{k=2}^{j-1} (-1)^{k-1} p_k \la
\tau_{1,j},\la\tau_{2,j},s_{1,j} \ra,\ldots,
\la\tau_{k-1,j},s_{k-2,j}\ra, \frac{\la \tau_{k,j},
{s}_{k-1,j}\ra}{\widehat{s}_{k+1,j}},\la s_{k+1,j},\sigma_k\ra
\widehat{\ra} +
\]
\begin{equation}
\label{eq:L}
(-1)^j  \sum_{k=j+1}^{m} p_k \la
\tau_{1,j},\la\tau_{2,j},s_{1,j} \ra,\ldots,
\la\tau_{j,j},s_{j-1,j}\ra,\la s_{j+1,k},\sigma_j\ra \widehat{\ra}.
\end{equation}
Now, it is not so clear who the auxiliary Nikishin system should be because some annoying ratios of Cauchy transforms have appeared. We shall see that already for $j=2$ there are two candidates, and for general $j$ the number of candidates equals $2^{j-1}$.

We can use (\ref{4.8}) (see the second inequalities in (\ref{4.6})-(\ref{4.2})) to obtain
\[  \frac{ {\mathcal{L}}_{\bf n}}{\widehat{s}_{1,j}}  = (\ell_{1,j}p_0  +  p_j  + \sum_{k\neq j, k=1}^m \frac{|s_{1,k}|}{|s_{1,j}|}
p_k) + (p_0 + \frac{|\la
 {s}_{2,j},\sigma_1  \ra |}{|{s}_{2,j}|}p_1)\tau_{1,j} +
\]
\[ \sum_{k=2}^{j-1}
 (-1)^{k-1} (p_{k-1} + \frac{|\la
 {s}_{k+1,j},\sigma_k  \ra |}{|{s}_{k+1,j}|} p_k) \la \tau_{1,j},\la
 \tau_{2,j},s_{1,j}\ra,\ldots,\la \tau_{k,j},s_{k-1,j} \ra
 \widehat{\ra} +
\]
\[ (-1)^{j-1} p_{j-1} \la \tau_{1,j}, \la
 \tau_{2,j},s_{1,j}\ra,\ldots,\la \tau_{j,j},s_{j-1,j} \ra
 \widehat{\ra} +
\]
\begin{equation} \label{4.9}
(-1)^j \sum_{k=j+1}^{m} p_k \la \tau_{1,j},\la\tau_{2,j},s_{1,j} \ra,\ldots,
\la\tau_{j,j},s_{j-1,j}\ra,\la s_{j+1,k},\sigma_j\ra \widehat{\ra}.
\end{equation}
(The sum $\sum_{k=2}^{j-1}$ is empty if $j=2$.)

If we are in the class ${\mathbb{Z}}_+^{m+1}(*)$ of multi-indices, and we take $j$ to be the first component for which $n_j = \max\{n_0+1,n_1,\ldots,n_m\}$, then  $n_0 \geq \cdots \geq n_{j-1}$. It follows that
\[ \deg (\ell_{1,j}p_0  +  p_j  + \sum_{k\neq j, k=1}^m \frac{|s_{1,k}|}{|s_{1,j}|}
p_k) \leq n_j -1
\]
and
\[ \deg (p_{k-1} + \frac{|\la
 {s}_{k+1,j},\sigma_k  \ra |}{|{s}_{k+1,j}|} p_k) \leq n_{k-1} -1,\qquad k=1,\ldots,j-1.
\]
Thus ${\mathcal{L}}_{\bf n}^*$ is the right hand side of (\ref{4.9}), which is a linear form generated by the multi-index ${\bf n}^* = (n_j,n_0\ldots,n_{j-1},n_{j+1},\ldots,n_m) \in {\mathbb{Z}}_+^{m+1}$ and the Nikishin system
\[{\mathcal{N}}(\sigma_1^*,\ldots,\sigma_m^*) = {\mathcal{N}}(\tau_{1,j},\la\tau_{2,j},s_{1,j} \ra,\ldots,
\la\tau_{j,j},s_{j-1,j}\ra,\la \sigma_{j+1},\sigma_j\ra, \sigma_{j+2},\ldots,\sigma_m).
\]

This would be sufficient to prove the AT property within the class ${\mathbb{Z}}_+^{m+1}(*)$ because it is easy to observe that then $(n_0,\ldots,n_{j-1},n_{j+1},\ldots,n_m) \in {\mathbb{Z}}_+^m(*)$ (see the proof of Theorem \ref{teo:2}). This result was first obtained in \cite[Theorem 2]{LIF}.

Of course, (\ref{4.9}) is still valid in the general case but, if it is not true that $n_0 \geq \ldots \geq n_{j-1}$, some of the degrees of the polynomials in the linear form on the right hand blow up with respect to the bounds established by the components of ${\bf n}^*$. We must proceed with caution. For this, we need two more reduction formulas which are contained in the next lemma.  \vspace{0,2cm}

Let ${\tau}_{\alpha,\beta;\gamma,\gamma}$ denote the inverse measure of
$\langle \langle \sigma_{\alpha}, \sigma_{\beta} \rangle,
\sigma_{\gamma} \ra$. That is,
\[ {1}/{\langle \langle
\sigma_{\alpha}, \sigma_{\beta} \rangle,
\sigma_{\gamma}\widehat{\rangle}(z)} =
\ell_{\alpha,\beta;\gamma,\gamma}(z) +
\widehat{\tau}_{\alpha,\beta;\gamma,\gamma}(z)
\]
where $\ell_{\alpha,\beta;\gamma,\gamma}$ denotes a first degree
polynomial. This notation seems unnecessarily
complicated. It is consistent with the one used later for more general
inverse measures which will be needed.

\bl \label{freedom} Let $\Delta_{\gamma},$ $\Delta_{\alpha}$ and
$\Delta_{\beta}$ be three intervals such that $\Delta_{\gamma} \cap
\Delta_{\alpha} = \emptyset = \Delta_{\beta} \cap \Delta_{\alpha} .$
Let $\sigma_{\gamma} \in {\mathcal{M}}(\Delta_{\gamma}),$
$\sigma_{\alpha} \in {\mathcal{M}}(\Delta_{\alpha})$ and
$\sigma_{\beta} \in {\mathcal{M}}(\Delta_{\beta}).$ Then for any $f
\in L_1(\sigma_{\gamma})$
\begin{equation}\label{inversa2*}
\frac{\widehat{\sigma}_{\alpha}(z)}{\langle
\sigma_{\alpha},\sigma_{\beta}\widehat{\rangle}(z)} \langle \langle
\tau_{\alpha,\alpha} ,\sigma_{\beta},\sigma_{\alpha}\rangle, f
\sigma_{\gamma}, \sigma_{\alpha} \widehat{\rangle}(z)= \la \frac{\la
\sigma_{\beta},\sigma_{\alpha}
\widehat{\ra}}{\widehat{\sigma}_{\beta}} \tau_{\alpha,\beta},
f\sigma_{\gamma}, \sigma_{\alpha},\sigma_{\beta} \widehat{\ra}(z),
\end{equation}
\begin{equation}\label{shirenu}
\frac{\langle
\sigma_{\alpha},\sigma_{\beta}\widehat{\rangle}(z)}{\langle \langle
\sigma_{\alpha}, \sigma_{\beta} \rangle,
\sigma_{\gamma}\widehat{\rangle}(z)} \la \frac{\la
\sigma_{\beta},\sigma_{\alpha}
\widehat{\ra}}{\widehat{\sigma}_{\beta}} \tau_{\alpha,\beta},
\sigma_{\gamma}, \sigma_{\alpha},\sigma_{\beta} \widehat{\ra}(z)=
\la  \frac{\langle \sigma_{\beta}
,\sigma_{\alpha},\sigma_{\gamma}\widehat{\rangle}}{\widehat{\sigma}_{\beta}
} \frac{\langle
\sigma_{\gamma},\sigma_{\alpha},\sigma_{\beta}\widehat{\rangle} }{
\widehat{\sigma}_{\gamma}} {\tau}_{\alpha,\beta;\gamma,\gamma}
\widehat{\ra}(z).
\end{equation}
\el
{\bf Proof.} Let us prove (\ref{inversa2*}). Taking into account
(\ref{2.1}) and (\ref{2.3}), we have that
\[
\langle \langle \tau_{\alpha,\alpha}
,\sigma_{\beta},\sigma_{\alpha}\rangle,\langle  f \sigma_{\gamma},
\sigma_{\alpha}\rangle \widehat{\rangle}(z)=\langle f
\sigma_{\gamma}, \sigma_{\alpha} \widehat{\rangle}(z)\langle
\tau_{\alpha,\alpha} ,\sigma_{\beta},
{\sigma}_{\alpha}\widehat{\ra}(z) -\langle \langle f
\sigma_{\gamma}, \sigma_{\alpha}\rangle, \tau_{\alpha,\alpha}
,\sigma_{\beta},{\sigma}_{\alpha}\widehat{\ra}(z) =
\]
\[
\langle f   \sigma_{\gamma},\sigma_{\alpha}
\widehat{\rangle}(z)\left(\frac{|\langle
\sigma_{\alpha},\sigma_{\beta}\rangle|}{|\sigma_{\alpha}|}-
\frac{\langle\sigma_{\alpha},\sigma_{\beta}\widehat{\rangle}(z)}{
\widehat{\sigma}_{\alpha}(z)}\right) -\int\left(\frac{|\langle
\sigma_{\alpha},\sigma_{\beta}\rangle|}{|\sigma_{\alpha}|}-
\frac{\langle\sigma_{\alpha},\sigma_{\beta}\widehat{\rangle}(x_{\gamma})}{
\widehat{\sigma}_{\alpha}(x_{\gamma})}\right)
\frac{f(x_{\gamma})d\langle\sigma_{\gamma}, {\sigma}_{\alpha}\ra
(x_{\gamma})}{z-x_{\gamma}}=
\]
\[
\int
\langle\sigma_{\alpha},\sigma_{\beta}\widehat{\rangle}(x_{\gamma})\frac{
f(x_{\gamma}) d\sigma_{\gamma} (x_{\gamma})}{z-x_{\gamma}}-\langle f
\sigma_{\gamma}, \sigma_{\alpha}
\widehat{\rangle}(z)\frac{\langle\sigma_{\alpha},\sigma_{\beta}\widehat{\rangle}(z)}{
\widehat{\sigma}_{\alpha}(z)}.
\]
This and (\ref{2.2}), render
\[
\frac{ \widehat{\sigma}_{\alpha}(z)}{\langle
\sigma_{\alpha},\sigma_{\beta}\widehat{\rangle}(z)} \langle \langle
\tau_{\alpha,\alpha} ,\sigma_{\beta},\sigma_{\alpha}\rangle, f
 \sigma_{\gamma}, \sigma_{\alpha} \widehat{\rangle}(z)=
\]
\[
\frac{ \widehat{\sigma}_{\alpha}(z)}{\langle
\sigma_{\alpha},\sigma_{\beta}\widehat{\rangle}(z)}\int
\langle\sigma_{\alpha},\sigma_{\beta}\widehat{\rangle}(x_{\gamma})
\frac{f(x_{\gamma})d\sigma_{\gamma}(x_{\gamma})}{z-x_{\gamma}}-\langle
f \sigma_{\gamma},\sigma_{\alpha} \widehat{\rangle}(z) =
\]
\[ -\langle
f \sigma_{\gamma},\sigma_{\alpha} \widehat{\rangle}(z)
 + \frac{|{\sigma}_{\alpha} |}{|\langle
\sigma_{\alpha},\sigma_{\beta} {\rangle}|} \la f \sigma_{\gamma},
\sigma_{\alpha}, \sigma_{\beta} \widehat{\ra} (z)  + \la f
\sigma_{\gamma}, \sigma_{\alpha}, \sigma_{\beta} \widehat{\ra} (z)
\la \frac{\tau_{\alpha,\beta}}{\widehat{\sigma}_{\beta}},\la
\sigma_{\beta},\sigma_{\alpha} \ra \widehat{\ra}(z)
\]

Since
\[
\frac{\widehat{\sigma}_{\alpha}(z)}{\langle
\sigma_{\alpha},\sigma_{\beta}\widehat{\rangle}(z)} \langle \langle
\tau_{\alpha,\alpha} ,\sigma_{\beta},\sigma_{\alpha}\rangle, f
 \sigma_{\gamma}, \sigma_{\alpha}
\widehat{\rangle}(z)={\mathcal{O}} \left(\frac{1}{z}\right)\in
{\mathcal{H}}(\overline{{\mathbb{C}}}\setminus \Delta_{\alpha}), \qquad z  \to \infty,
\]
that $-\langle f  \sigma_{\gamma},\sigma_{\alpha} \widehat{\rangle}
 + \frac{|{\sigma}_{\alpha} |}{|\langle
\sigma_{\alpha},\sigma_{\beta} {\rangle}|} \la f \sigma_{\gamma},
\sigma_{\alpha}, \sigma_{\beta} \widehat{\ra}$ and $\la f
\sigma_{\gamma}, \sigma_{\alpha}, \sigma_{\beta} \widehat{\ra} $ are
analytic on a neighborhood of $\Delta_{\alpha},$ on account of
(\ref{eq:g}), we obtain (\ref{inversa2*}).

Now, we prove (\ref{shirenu}). From (\ref{2.3}) and (\ref{2.1})
\[ \la \frac{\la
\sigma_{\beta},\sigma_{\alpha}
\widehat{\ra}}{\widehat{\sigma}_{\beta}} \tau_{\alpha,\beta},
\sigma_{\gamma}, \sigma_{\alpha},\sigma_{\beta} \widehat{\ra}(z) =
\]
\[
\frac{|\langle \sigma_{\beta}, \sigma_{\alpha}\rangle|}{|\langle
\sigma_{\beta}\rangle|} \langle \tau_{\alpha,\beta},
\sigma_{\gamma},\sigma_{\alpha},\sigma_{\beta} \widehat{\rangle}(z)
-\langle  \langle \tau_{\alpha,\beta},
\sigma_{\gamma},\sigma_{\alpha},\sigma_{\beta} \rangle,
\tau_{\beta,\beta},\sigma_{\alpha},\sigma_{\beta}\widehat{\rangle}(z)=
\]
\[
\frac{|\langle \sigma_{\beta}, \sigma_{\alpha}\rangle|}{|\langle
\sigma_{\beta}\rangle|} \langle \tau_{\alpha,\beta},
\sigma_{\gamma},\sigma_{\alpha},\sigma_{\beta} \widehat{\rangle}(z)
- \langle \tau_{\alpha,\beta},
\sigma_{\gamma},\sigma_{\alpha},\sigma_{\beta} \widehat{\rangle}(z)
\langle
\tau_{\beta,\beta},\sigma_{\alpha},\sigma_{\beta}\widehat{\rangle}(z)+
\]
\[
\langle \langle
\tau_{\beta,\beta},\sigma_{\alpha},\sigma_{\beta}\rangle,\langle
\tau_{\alpha,\beta},
\sigma_{\gamma},\sigma_{\alpha},\sigma_{\beta}\ra
\widehat{\rangle}(z)=
\]
\[
\frac{\langle \sigma_{\beta}, {\sigma}_{\alpha}\widehat{\ra}(z)}{
\widehat{ \sigma}_{\beta} (z)}\langle \tau_{\alpha,\beta},
\sigma_{\gamma},\sigma_{\alpha},\sigma_{\beta}
\widehat{\rangle}(z)+\langle \langle
\tau_{\beta,\beta},\sigma_{\alpha},\sigma_{\beta}\rangle,\langle
\tau_{\alpha,\beta},
\sigma_{\gamma},\sigma_{\alpha},\sigma_{\beta}\ra
\widehat{\rangle}(z) =
\]
\[ \frac{\langle \sigma_{\beta},
 {\sigma}_{\alpha}\widehat{\ra}(z)}{ \widehat{ \sigma}_{\beta} (z)}
\left(\frac{|\langle \langle
\sigma_{\alpha},\sigma_{\beta}\rangle,\sigma_{\gamma}\rangle|}{|\langle\sigma_{\alpha},\sigma_{\beta}
\rangle|}-\frac{\langle\langle
\sigma_{\alpha},\sigma_{\beta}\rangle,\sigma_{\gamma}\widehat{\rangle}(z)}{\langle\sigma_{\alpha},\sigma_{\beta}\widehat{\rangle}(z)}\right)+
\]
\[
\int \left(\frac{|\langle \langle
\sigma_{\alpha},\sigma_{\beta}\rangle,\sigma_{\gamma}\rangle|}{|\langle\sigma_{\alpha},\sigma_{\beta}
\rangle|}-\frac{\langle\langle
\sigma_{\alpha},\sigma_{\beta}\rangle,\sigma_{\gamma}\widehat{\rangle}(x_{\beta})}{
\la
\sigma_{\alpha},\sigma_{\beta}\widehat{\rangle}(x_{\beta})}\right)
\frac{ d \langle
\tau_{\beta,\beta},\sigma_{\alpha},\sigma_{\beta} {\rangle}(x_{\beta})}{z-x_{\beta}}=
\]
\[
-\frac{\langle \sigma_{\beta},  {\sigma}_{\alpha}\widehat{\ra}(z)}{
\widehat{\sigma}_{\beta} (z)}\frac{\langle\langle
\sigma_{\alpha},\sigma_{\beta}\rangle,\sigma_{\gamma}\widehat{\rangle}(z)}{\langle\sigma_{\alpha},\sigma_{\beta}\widehat{\rangle}(z)}+
\]
\begin{equation}\label{final*}
\left(\frac{\langle \sigma_{\beta},
{\sigma}_{\alpha}\widehat{\ra}(z)}{ \widehat{\sigma}_{\beta}
(z)}+\langle
\tau_{\beta,\beta},\sigma_{\alpha},\sigma_{\beta}\widehat{\rangle}(z)\right)\frac{|\langle
\langle
\sigma_{\alpha},\sigma_{\beta}\rangle,\sigma_{\gamma}\rangle|}{|\langle\sigma_{\alpha},\sigma_{\beta}
\rangle|}-\langle\tau_{\beta,\beta} ,\langle
\sigma_{\alpha},\sigma_{\beta}
\rangle,\sigma_{\gamma}\widehat{\rangle}(z) =
\end{equation}
\[
-\frac{\langle \sigma_{\beta},  {\sigma}_{\alpha}\widehat{\ra}(z)}{
\widehat{\sigma}_{\beta}  (z)}\frac{\langle\langle
\sigma_{\alpha},\sigma_{\beta}\rangle,\sigma_{\gamma}\widehat{\rangle}(z)}{\langle\sigma_{\alpha},\sigma_{\beta}\widehat{\rangle}(z)}+\frac{|\langle
\sigma_{\beta}, \sigma_{\alpha}\rangle|}{|\langle
\sigma_{\beta}\rangle|}\frac{|\langle \langle
\sigma_{\alpha},\sigma_{\beta}\rangle,\sigma_{\gamma}\rangle|}{|\langle\sigma_{\alpha},\sigma_{\beta}
\rangle|}-\langle\tau_{\beta,\beta} ,\langle
\sigma_{\alpha},\sigma_{\beta}\rangle,
\sigma_{\gamma}\widehat{\rangle}(z) =
\]
\[
-\frac{\langle \sigma_{\beta},  {\sigma}_{\alpha}\widehat{\ra}(z)}{
\widehat{\sigma}_{\beta}  (z)}\frac{\langle\langle
\sigma_{\alpha},\sigma_{\beta}\rangle,\sigma_{\gamma}\widehat{\rangle}(z)}
{\langle\sigma_{\alpha},\sigma_{\beta}\widehat{\rangle}(z)}-
\frac{|\langle \langle
\sigma_{\alpha},\sigma_{\beta}\rangle,\sigma_{\gamma}\rangle|}{|\sigma_{\beta}
 |}- \langle\tau_{\beta,\beta} ,\langle
\sigma_{\alpha},\sigma_{\beta}\rangle,
\sigma_{\gamma}\widehat{\rangle}(z) =
\]
\[ -\frac{\langle \sigma_{\beta},  {\sigma}_{\alpha} \widehat{\ra}(z)}{
\widehat{\sigma}_{\beta}(z)}\frac{\langle\langle
\sigma_{\alpha},\sigma_{\beta}\rangle,\sigma_{\gamma}\widehat{\rangle}(z)}
{\langle\sigma_{\alpha},\sigma_{\beta}\widehat{\rangle}(z)} +
\frac{\la \sigma_{\beta},\sigma_{\alpha},\sigma_{\gamma}
\widehat{\ra}(z)}{\widehat{\sigma}_{\beta}(z)}.
\]
In the second last equality above, we employed that
\[
0 = \lim_{z \to \infty} z  \widehat{\sigma}_{\alpha}  (z)
\widehat{\sigma}_{\beta}(z)=\lim_{z \to \infty}z \langle
\sigma_{\alpha},\sigma_{\beta} \widehat{\rangle} (z) + \lim_{z \to
\infty}z \langle \sigma_{\beta},\sigma_{\alpha} \widehat{\rangle}
(z) = |\la \sigma_{\alpha},\sigma_{\beta} \ra| + |\la
\sigma_{\beta},\sigma_{\alpha} \ra|,
\]
which implies that  $|\la \sigma_{\alpha},\sigma_{\beta} \ra| = -
|\la \sigma_{\beta},\sigma_{\alpha} \ra|$. Analogously, $- |\la\la
\sigma_{\alpha},\sigma_{\beta}\ra,\sigma_{\gamma} \ra| = - |\la\la
\sigma_{\alpha},\sigma_{\gamma}\ra,\sigma_{\beta} \ra| = |\la
\sigma_{\beta},\la \sigma_{\alpha},\sigma_{\gamma} \ra \ra|  $ and
by (\ref{2.3})
\[ \frac{\la \sigma_{\beta},\la \sigma_{\alpha},\sigma_{\gamma}\ra
\widehat{\ra} }{\widehat{\sigma}_{\beta} } = \frac{|\la
\sigma_{\beta},\la \sigma_{\alpha},\sigma_{\gamma}\ra
 {\ra} |}{ |{\sigma}_{\beta}|} - \la \tau_{\beta,\beta}, \la \la
 \sigma_{\alpha},\sigma_{\gamma} \ra, \sigma_{\beta} \ra
 \widehat{\ra} = -
\frac{|\langle \langle
\sigma_{\alpha},\sigma_{\beta}\rangle,\sigma_{\gamma}\rangle|}{|\sigma_{\beta}
 |}- \langle\tau_{\beta,\beta} ,\langle
\sigma_{\alpha},\sigma_{\beta}\rangle,
\sigma_{\gamma}\widehat{\rangle},
\]
which is used in the last equality. Therefore, from the chain of
equations and (\ref{2.2}), we derive that
\[ \frac{\langle
\sigma_{\alpha},\sigma_{\beta}\widehat{\rangle}(z)}{\langle \langle
\sigma_{\alpha}, \sigma_{\beta} \rangle,
\sigma_{\gamma}\widehat{\rangle}(z)} \la \frac{\la
\sigma_{\beta},\sigma_{\alpha}
\widehat{\ra}}{\widehat{\sigma}_{\beta}} \tau_{\alpha,\beta},
\sigma_{\gamma}, \sigma_{\alpha},\sigma_{\beta} \widehat{\ra}(z) =
-\frac{\langle \sigma_{\beta}, {\sigma}_{\alpha}\widehat{\ra}(z)}{
\widehat{\sigma}_{\beta}(z)}+ \frac{\la
\sigma_{\beta},\sigma_{\alpha},\sigma_{\gamma}
\widehat{\ra}(z)}{\widehat{\sigma}_{\beta}(z)}\frac{\langle\sigma_{\alpha},\sigma_{\beta}\widehat{\rangle}(z)}
{\langle\langle\sigma_{\alpha},\sigma_{\beta}\rangle,\sigma_{\gamma}\widehat{\rangle}(z)}
=
\]
\[
-\frac{\langle \sigma_{\beta}, {\sigma}_{\alpha}\widehat{\ra}(z)}{
\widehat{\sigma}_{\beta}(z)}+ \frac{\la
\sigma_{\beta},\sigma_{\alpha},\sigma_{\gamma}
\widehat{\ra}(z)}{\widehat{\sigma}_{\beta}(z)}\frac{|\langle\sigma_{\alpha},\sigma_{\beta}
\rangle|}
{|\langle\langle\sigma_{\alpha},\sigma_{\beta}\rangle,\sigma_{\gamma}
{\rangle}|} + \frac{\la
\sigma_{\beta},\sigma_{\alpha},\sigma_{\gamma}
\widehat{\ra}(z)}{\widehat{\sigma}_{\beta}(z)} \la  \frac{
\tau_{\alpha,\beta;\gamma,\gamma} }{\widehat{\sigma}_{\gamma}},
\sigma_{\gamma},\sigma_{\alpha},\sigma_{\beta} \widehat{\ra}(z).
\]

Taking into consideration that
\[ \frac{\langle
\sigma_{\alpha},\sigma_{\beta}\widehat{\rangle}(z)}{\langle \langle
\sigma_{\alpha}, \sigma_{\beta} \rangle,
\sigma_{\gamma}\widehat{\rangle}(z)} \la \frac{\la
\sigma_{\beta},\sigma_{\alpha}
\widehat{\ra}}{\widehat{\sigma}_{\beta}} \tau_{\alpha,\beta},
\sigma_{\gamma}, \sigma_{\alpha},\sigma_{\beta} \widehat{\ra}(z)
={\mathcal{O}} \left(\frac{1}{z}\right) \in
{\mathcal{H}}(\overline{{\mathbb{C}}}\setminus \Delta_{\alpha}), \qquad z \to \infty,
\]
whereas $-\frac{\langle \sigma_{\beta},
{\sigma}_{\alpha}\widehat{\ra} }{ \widehat{\sigma}_{\beta} }+
\frac{\la \sigma_{\beta},\sigma_{\alpha},\sigma_{\gamma}
\widehat{\ra} }{\widehat{\sigma}_{\beta}
}\frac{|\langle\sigma_{\alpha},\sigma_{\beta} \rangle|}
{|\langle\langle\sigma_{\alpha},\sigma_{\beta}\rangle,\sigma_{\gamma}
{\rangle}|}$ and $\frac{\la
\sigma_{\beta},\sigma_{\alpha},\sigma_{\gamma} \widehat{\ra}
}{\widehat{\sigma}_{\beta} }$ are analytic on a neighborhood of
$\Delta_{\alpha}$, using (\ref{eq:g}) we obtain (\ref{shirenu}). \fp
\vspace{0,2cm}

{\bf Proof of Lemma \ref{lem:4} when $2 \leq j \leq m$.}  Set
\[ {\mathcal{L}}_{\bf n}^* = p_0^* + p_0 \widehat{\tau}_{1,j}
+ p_1\la
\frac{\tau_{1,j}}{\widehat{s}_{2,j}}, \la s_{2,j},\sigma_1 \ra
\widehat{\ra}  +
\]
\[ \sum_{k=2}^{j-1} (-1)^{k-1} p_k \la
\tau_{1,j},\la\tau_{2,j},s_{1,j} \ra,\ldots,
\la\tau_{k-1,j},s_{k-2,j}\ra, \frac{\la \tau_{k,j},
{s}_{k-1,j}\ra}{\widehat{s}_{k+1,j}},\la s_{k+1,j},\sigma_k\ra
\widehat{\ra} +
\]
\begin{equation}
\label{eq:A}
(-1)^j  \sum_{k=j+1}^{m} p_k \la
\tau_{1,j},\la\tau_{2,j},s_{1,j} \ra,\ldots,
\la\tau_{j,j},s_{j-1,j}\ra,\la s_{j+1,k},\sigma_j\ra \widehat{\ra},
\end{equation}
where
\[ p_0^* = \ell_{1,j}p_0  +  p_j  + \sum_{k\neq j, k=1}^m \frac{|s_{1,k}|}{|s_{1,j}|}
p_k, \quad \deg p_0^* \leq n_j -1 = n_0^* -1.
\]
We took this function from the right hand side of (\ref{eq:L}). We must show that there exist a multi-index ${\bf n}^* \in {\mathbb{Z}}_+^{m+1}$, which is a permutation of ${\bf n}$,
and a Nikishin system ${\mathcal{N}}(\sigma_0^*,\ldots,\sigma_m^*)$ which allow to express ${\mathcal{L}}_{\bf n}^*$ as a linear form generated by them with polynomials with real coefficients. So far, $n_0^*$ defined above serves the purpose of being the first component of ${\bf n}^*$ and $p_0^*$ of being the polynomial part of the linear form.

{\bf First step.} Is $n_0 \geq n_1$ or $n_0 \leq n_1$? (When $n_0 =
n_1$ we can proceed either ways.)

\noindent {\bf A1)} If $n_0 \geq n_1$, take ${n}_1^* = n_0$  and $\sigma_1^* = \tau_{1,j}$.
Decompose $\frac{ \la s_{2,j},\sigma_1
\widehat{\ra}}{\widehat{s}_{2,j}}$ using (\ref{2.3}). Then, the
first three terms of ${\mathcal{L}}_{\bf n}^*$ are
\[ p_0^* + p_0 \widehat{\sigma}_{1}^*
+ p_1\la
\frac{\sigma_{1}^*}{\widehat{s}_{2,j}}, \la s_{2,j},\sigma_1 \ra
\widehat{\ra}  = p_0^* + (p_0
+ \frac{|\la
 {s}_{2,j},\sigma_1  \ra |}{|{s}_{2,j}|}p_1) \widehat{\sigma}_{1}^* - p_1 \la {\sigma}_{1}^*, \la \tau_{2,j},
s_{1,j} \ra \widehat{\ra}.
\]
Consequently, taking  $ {p}_1^* =  p_0 + \frac{|\la
 {s}_{2,j},\sigma_1  \ra |}{|{s}_{2,j}|}p_1$, we have that $\deg  {p}_1^* \leq  {n}_1^* -1 (= n_0-1)$.

In case that $j=2$, we obtain
\[ {\mathcal{L}}_{\bf n}^* = p_0^* +  {p}_1^* \widehat{\tau}_{1,2} - p_1 \la \tau_{1,2},
 \la \tau_{2,2},s_{1,2} \ra \widehat{\ra} + \sum_{k=3}^{m} p_k
\la \tau_{1,2}, \la\tau_{2,2},s_{1,2} \ra, \la s_{3,k},\sigma_2 \ra
\widehat{\ra}
\]
(compare with (\ref{4.9})). Then, the proof would be complete taking
$ {\bf n}^* = (n_2, n_0,n_1,n_3\ldots,n_m)$ and the Nikishin
system
\[{\mathcal{N}}(\sigma_1^*,\ldots,\sigma_m^*) = {\mathcal{N}}(\tau_{1,2},\la \tau_{2,2},s_{1,2}\ra, \la
\sigma_3,\sigma_2 \ra, \sigma_4,\ldots,\sigma_m ).
\]
(If $m=2$, then
$ {\bf n}^* = (n_2, n_0,n_1)$ and the Nikishin system is
${\mathcal{N}}(\tau_{1,2}, \la\tau_{2,2},s_{1,2}\ra)$.)

If $j \geq 3$, we obtain
\[ {\mathcal{L}}_{\bf n}^* = p_0^* + {p}_1^* \widehat{\sigma}_{1}^* - p_1 \la
 \sigma_{1}^*, \la\tau_{2,j},s_{1,j} \ra \widehat{\ra} - p_2
\la \sigma_{1}^*, \frac{\la s_{3,j},\sigma_2  \widehat{\ra}}{ \widehat{s}_{3,j} }
\la \tau_{2,j}, s_{1,j} \ra \widehat{\ra} +
\]
\[\sum_{k=3}^{j-1} (-1)^{k-1} p_k \la \sigma_{1}^*, \la\tau_{2,j},s_{1,j}
\ra,\ldots, \la\tau_{k-1,j},s_{k-2,j}\ra, \frac{\la \tau_{k,j},
{s}_{k-1,j}\ra}{\widehat{s}_{k+1,j}},\la s_{k+1,j},\sigma_k\ra
\widehat{\ra}+   \]
\[  (-1)^j \sum_{k=j+1}^{m} p_k
\la \sigma_{1}^*,\la\tau_{2,j},s_{1,j} \ra,\ldots, \la\tau_{j,j},s_{j-1,j}\ra,\la
s_{j+1,k},\sigma_j\ra \widehat{\ra}.
\]

\noindent {\bf B1)} If $n_0 \leq n_1$, take $ {n}_1^* = n_1$ and $\sigma_1^* =  \frac{\la s_{2,j},\sigma_1\widehat{\ra}}{\widehat{s}_{2,j}} \tau_{1,j}$. We can rewrite (\ref{eq:A}) as follows
\[ {\mathcal{L}}_{\bf n}^* = p_0^* + p_1\widehat{\sigma}_1^*  + p_0 \la \frac{ \widehat{s}_{2,j}}{\la s_{2,j},\sigma_1
\widehat{\ra}} {\sigma}_1^*\widehat{\ra}
+
\]
\[ \sum_{k=2}^{j-1} (-1)^{k-1} p_k \la
\frac{ \widehat{s}_{2,j}}{\la s_{2,j},\sigma_1
\widehat{\ra}} {\sigma}_1^*,\la\tau_{2,j},s_{1,j} \ra,\ldots,
\la\tau_{k-1,j},s_{k-2,j}\ra, \frac{\la \tau_{k,j},
{s}_{k-1,j}\ra}{\widehat{s}_{k+1,j}},\la s_{k+1,j},\sigma_k\ra
\widehat{\ra} +
\]
\begin{equation}
\label{eq:B}
(-1)^j  \sum_{k=j+1}^{m} p_k \la
\frac{ \widehat{s}_{2,j}}{\la s_{2,j},\sigma_1
\widehat{\ra}} {\sigma}_1^*,\la\tau_{2,j},s_{1,j} \ra,\ldots,
\la\tau_{j,j},s_{j-1,j}\ra,\la s_{j+1,k},\sigma_j\ra \widehat{\ra},
\end{equation}

Decompose $\frac{ \widehat{s}_{2,j}}{\la s_{2,j},\sigma_1
\widehat{\ra}}$ using (\ref{2.2}). Then, the
first three terms of ${\mathcal{L}}_{\bf n}^*$ in (\ref{eq:B}) can be expressed as
\[ p_0^* + p_1 \widehat{\sigma}_{1}^*
+ p_0 \la \frac{ \widehat{s}_{2,j}}{\la s_{2,j},\sigma_1
 \widehat{\ra}}  {\sigma}_{1}^* \widehat{\ra} = p_0^* + (p_1
+ \frac{|{s}_{2,j}|}{|\la
 {s}_{2,j},\sigma_1  \ra |}p_0) \widehat{\sigma}_{1}^* + p_0 \la \sigma_1^*, \frac{\widehat{s}_{1,j}}{\widehat{\sigma}_1}  \tau_{2,j;1,1} \widehat{\ra}.
\]
Taking  $ {p}_1^* = p_1 + \frac{|{s}_{2,j}|}{|\la
 {s}_{2,j},\sigma_1  \ra |}p_0  $, we have that $\deg  {p}_1^* \leq  {n}_1^* -1 (= n_1-1)$.

If $j=2$, due to (\ref{inversa2*}) (in the next formula
$\widehat{s}_{4,k} \equiv 1$ if $k=3$)
\[ \frac{\widehat{s}_{2,2}}{\la s_{2,2},\sigma_1 \widehat{\ra}}
\la \la\tau_{2,2},s_{1,2} \ra,\la s_{3,k},\sigma_2\ra \widehat{\ra}
= \frac{\widehat{\sigma}_{2}}{\la \sigma_{2},\sigma_1 \widehat{\ra}}
\la \la\tau_{2,2},s_{1,2} \ra,\widehat{s}_{4,k} \sigma_{3},\sigma_2
\widehat{\ra} =
\]
\[ \la \frac{\la \sigma_{1},\sigma_{2}
\widehat{\ra}}{\widehat{\sigma}_{1}} \tau_{2,2;1,1},
\widehat{s}_{4,k} \sigma_{3}, \sigma_{2},\sigma_{1}  \widehat{\ra} =
\la \frac{\la \sigma_{1},\sigma_{2}
\widehat{\ra}}{\widehat{\sigma}_{1}} \tau_{2,2;1,1}, \la \sigma_{3},
\sigma_{2},\sigma_{1} \ra, s_{4,k}\widehat{\ra}.
\]
Consequently,
\[ {\mathcal{L}}_{\bf n}^* = p_0^* + p_1^* \widehat{\sigma}_{1}^* + p_0 \la \sigma_1^*, \frac{\widehat{s}_{1,2}}{\widehat{\sigma}_1}  \tau_{2,2;1,1} \widehat{\ra} + \sum_{k=3}^m p_k \la \sigma_1^*, \frac{\widehat{s}_{1,2}}{\widehat{\sigma}_1}  \tau_{2,2;1,1} , \la \sigma_{3},
\sigma_{2},\sigma_{1} \ra, s_{4,k}\widehat{\ra}.
\]
In this situation, we would be done considering $ {\bf n}^* =
(n_2, n_1,n_0,n_3,\ldots,n_m)$ and the system
\[{\mathcal{N}}(\sigma_1^*,\ldots,\sigma_m^*) = {\mathcal{N}}(\frac{\la \sigma_{2},\sigma_1\widehat{\ra}}{\widehat{\sigma}_{2}} \tau_{1,2},\frac{\la \sigma_{1},\sigma_{2}
\widehat{\ra}}{\widehat{\sigma}_{1}} \tau_{2,2;1,1}, \la \sigma_{3},
\sigma_{2},\sigma_{1} \ra, \sigma_{4},\ldots,\sigma_m).
\]
(Should
$m=2$, then $ {\bf n}^* = (n_2,n_1,n_0)$ and the Nikishin system
is ${\mathcal{N}}(\frac{\la \sigma_{2},\sigma_1\widehat{\ra}}{\widehat{\sigma}_{2}} \tau_{1,j},\frac{\la \sigma_{1},\sigma_{2}
\widehat{\ra}}{\widehat{\sigma}_{1}} \tau_{2,2;1,1})$.) This result
already includes new multi-indices for which it is possible to prove
normality. The only sub-case studied previously (see \cite{LF1}) was
for $m=2$. Therefore, if $j=2$ we are done.

Let us assume that $j \geq 3$. (The algorithm ends after $j-1$
steps.)  So far, we have used the notation $s_{k,l}$ only with $k \leq l$. We will extend its meaning to $k > l$ in which case
\[s_{k,l} = \la \sigma_k, \sigma_{k-1}, \ldots, \sigma_l \ra, \qquad k > l.
\]
Notice that if $l < k < j$, then
\[ \la s_{k,j}, s_{k-1, l} \ra = \la s_{k,l}, s_{k+1,j} \ra.
\]
The inverse measure of $\la s_{k,j}, s_{k-1, l} \ra$ we denote by $\tau_{k,j;k-1,l}$; that is,
\[ {1}/{\la s_{k,j}, s_{k-1, l} \widehat{\ra}} = \ell_{k,j;k-1,l} +  \widehat{\tau}_{k,j;k-1,l}
\]
In particular, $\tau_{2,j;1,1}$ denotes the inverse measure of $\la s_{2,j},\sigma_1 \ra$.

Let us transform the measures in $\sum_{k=2}^{j-1}$ of (\ref{eq:B}). Regarding the term with $p_2$, using (\ref{shirenu}) with
$\sigma_{\alpha} = \sigma_2, \sigma_{\beta} = s_{3,j}$ and
$\sigma_{\gamma} = \sigma_1$, we obtain
\[ \frac{\widehat{s}_{2,j}}{\la s_{2,j},\sigma_1
\widehat{\ra}}\la \frac{\la \tau_{2,j}, s_{1,j}\ra}{
\widehat{s}_{3,j} }, \la s_{3,j},\sigma_2 \ra \widehat{\ra} = \la
\frac{\la s_{3,j},\sigma_2,\sigma_1\widehat{\ra}}{\widehat{s}_{3,j}}
\frac{\widehat{s}_{1,j}}{\widehat{\sigma}_1}\tau_{2,j;1,1}
\widehat{\ra}.
\]
For $j=3$, $\sum_{k=3}^{j-1}$ is empty, so here the formulas make
sense when $j\geq 4$. Using (\ref{inversa2*}), with $\sigma_{\alpha}
= s_{2,j},$ $\sigma_{\beta} = \sigma_1,$ $\sigma_{\gamma} =
\tau_{3,j}$, and
\[ f = f_{1,j,k} = \left\{
\begin{array}{ll}
\displaystyle{\frac{\la s_{4,j}, \sigma_3 \widehat{\ra}}{\widehat{s}_{4,j}}}, &  3 = k < j \leq m, \\
\la \displaystyle{\frac{\la
\tau_{4,j},s_{3,j}\ra}{\widehat{s}_{5,j}}},\la s_{5,j},\sigma_4 \ra
\widehat{\ra}, &  4 = k < j \leq m, \\
\la \la \tau_{4,j},s_{3,j}\ra,\ldots,\la \tau_{k-1,j},s_{k-2,j} \ra,
\displaystyle{\frac{\la \tau_{k,j},s_{k-1,j}
\ra}{\widehat{s}_{k+1,j}}} , \la s_{k+1,j},\sigma_k \ra
\widehat{\ra}, &  5 \leq k < j \leq m,
\end{array}
\right.
\]
we obtain
\[ \frac{\widehat{s}_{2,j}}{\la s_{2,j},\sigma_1 \widehat{\ra}} \la \la\tau_{2,j},s_{1,j}
\ra,\ldots, \la\tau_{k-1,j},s_{k-2,j}\ra, \frac{\la \tau_{k,j},
{s}_{k-1,j}\ra}{\widehat{s}_{k+1,j}},\la s_{k+1,j},\sigma_k\ra
\widehat{\ra} =
\]
\[ \la  {\frac{\widehat{s}_{1,j}}{\widehat{\sigma}_1}} \tau_{2,j;1,1} ,
f_{1,j,k} \tau_{3,j},s_{2,j},\sigma_1  \widehat{\ra} =
\]
\[
\left\{
\begin{array}{ll}
\la \displaystyle{\frac{\widehat{s}_{1,j}}{\widehat{\sigma}_1}}
\tau_{2,j;1,1} , \displaystyle{\frac{\la \tau_{3,j},s_{2,j},\sigma_1
\ra }{\widehat{s}_{4,j}}}, \la s_{4,j},\sigma_3 \ra
\widehat{\ra}, &  3 = k < j \leq m, \\
\la \displaystyle{\frac{\widehat{s}_{1,j}}{\widehat{\sigma}_1}}
\tau_{2,j;1,1} , \la \tau_{3,j},s_{2,j},\sigma_1 \ra,\displaystyle{
\frac{\la \tau_{4,j}, s_{3,j} \ra}{\widehat{s}_{5,j}} },\la s_{5,j},
\sigma_4 \ra \widehat{\ra}, &  4 = k < j \leq m,\\
\la  {\frac{\widehat{s}_{1,j}}{\widehat{\sigma}_1}} \tau_{2,j;1,1} ,
\la \widehat{\sigma}_1 s_{2,j} \widehat{\ra}\tau_{3,j}, \widehat{s}_{3,j} \tau_{4,j},
\ldots, \widehat{s}_{k-2,j}\tau_{k-1,j} ,  {\frac{\la
\widehat{\sigma}_k s_{k+1,j} \widehat{\ra}\widehat{s}_{k-1,j} }{\widehat{s}_{k+1,j}}}
\tau_{k,j}  \widehat{\ra}, & 5 \leq k < j \leq m.
\end{array}
\right.
\]
In the last row, we used a more compact notation to fit the line. Notice that it is the same as
\[ \la  {\frac{\widehat{s}_{1,j}}{\widehat{\sigma}_1}} \tau_{2,j;1,1} ,
\la \tau_{3,j}, s_{2,j},\sigma_1  {\ra}, \la \tau_{4,j},  {s}_{3,j} \ra,
\ldots, \la \tau_{k-1,j},s_{k-2,j}\ra,  {\frac{\la
\tau_{k,j},s_{k-1,j} \ra}{\widehat{s}_{k+1,j}}}, \la
s_{k+1,j},\sigma_k \ra \widehat{\ra}.
\]

As for the terms in  $\sum_{k = j+1}^m$, applying
(\ref{inversa2*}) with $\sigma_{\alpha} = s_{2,j},$ $\sigma_{\beta}
= \sigma_1,$ $\sigma_{\gamma} = \tau_{3,j}$, and
\[ f  = f_{1,j,k}= \left\{
\begin{array}{ll}
\la s_{4,k}, \sigma_3 \widehat{\ra},  & 3 = j < k \leq m, \\
\la \la \tau_{4,j},s_{3,j}\ra,\ldots,\la \tau_{j,j},s_{j-1,j}
\ra,\la s_{j+1,k},\sigma_j \ra  \widehat{\ra}, &  4 \leq j < k \leq
m,
\end{array}
\right.
\]
it follows that
\[ \frac{\widehat{s}_{2,j}}{\la s_{2,j},\sigma_1 \widehat{\ra}}
\la \la\tau_{2,j},s_{1,j} \ra,\ldots, \la\tau_{j,j},s_{j-1,j}\ra,\la
s_{j+1,k},\sigma_j\ra \widehat{\ra} =
\]
\[ = \la \displaystyle{\frac{\widehat{s}_{1,j}}{\widehat{\sigma}_1}}
\tau_{2,j;1,1} , f_{1,j,k} \tau_{3,j},s_{2,j},\sigma_1 \widehat{\ra}
=
\]
\[
\left\{
\begin{array}{ll}
\la \displaystyle{\frac{\widehat{s}_{1,3}}{\widehat{\sigma}_1}}
\tau_{2,3;1,1} , \la \tau_{3,3},s_{2,3},\sigma_1 \ra, \la
s_{4,k},\sigma_3 \ra
\widehat{\ra},  & 3 = j < k \leq m, \\
\la \displaystyle{\frac{\widehat{s}_{1,j}}{\widehat{\sigma}_1}}
\tau_{2,j;1,1} , \la \tau_{3,j},s_{2,j},\sigma_1 \ra,\la \tau_{4,j},
s_{3,j} \ra, \ldots, \la \tau_{j,j},s_{j-1,j}\ra, \la s_{j+1,k},
\sigma_j \ra \widehat{\ra},  & 4 \leq j < k \leq m.
\end{array}
\right.
\]
When $j=m$ no such terms exist.

Therefore,
\[
 {\mathcal{L}}_{\bf n}^* = p_0^* + p_1^* \widehat{\sigma}_{1}^* + p_0 \la \sigma_1^*, \frac{\widehat{s}_{1,j}}{\widehat{\sigma}_1}  \tau_{2,j;1,1} \widehat{\ra}
  - p_2 \la \sigma_1^*, \frac{ \la
s_{3,j},\sigma_2,\sigma_1\widehat{\ra}}{\widehat{s}_{3,j}}
\frac{\widehat{s}_{1,j}}{\widehat{\sigma}_1}\tau_{2,j;1,1}\widehat{\ra}
+
\]
\[ \sum_{k=3}^{j-1} (-1)^{k-1} p_k \la \sigma_1^*,\displaystyle{\frac{\widehat{s}_{1,j}}{\widehat{\sigma}_1}}
\tau_{2,j;1,1} , f_{1,j,k} \tau_{3,j},s_{2,j},\sigma_1 \widehat{\ra}
+ (-1)^j \sum_{k= j+1}^m p_k \la \sigma_1^*,
\displaystyle{\frac{\widehat{s}_{1,j}}{\widehat{\sigma}_1}}
\tau_{2,j;1,1} , f_{1,j,k} \tau_{3,j},s_{2,j},\sigma_1
\widehat{\ra}.
\]
Using the notation for $f_{1,j,k}$ defined previously, in {\bf A1)} we ended up with
\[ {\mathcal{L}}_{\bf n}^* = p_0^* + {p}_1^* \widehat{\sigma}_{1}^* - p_1 \la
 \sigma_{1}^*, \la\tau_{2,j},s_{1,j} \ra \widehat{\ra} - p_2
\la \sigma_{1}^*, \frac{\la s_{3,j},\sigma_2  \widehat{\ra}}{ \widehat{s}_{3,j} }
\la \tau_{2,j}, s_{1,j} \ra \widehat{\ra} +
\]
\[\sum_{k=3}^{j-1} (-1)^{k-1} p_k \la \sigma_{1}^*, \la\tau_{2,j},s_{1,j}
\ra,f_{1,j,k} \tau_{3,j}, s_{2,j} \widehat{\ra} +
  (-1)^j \sum_{k=j+1}^{m} p_k
\la \sigma_{1}^*,\la\tau_{2,j},s_{1,j} \ra,f_{1,j,k} \tau_{3,j}, s_{2,j} \widehat{\ra}.
\]
Denote
\[  \sigma_2^{(l_1)} = \left\{
\begin{array}{ll}
\widehat{s}_{1,j}\tau_{2,j},  & l_1=1, \\
\frac{\widehat{s}_{1,j}}{\widehat{\sigma}_1}\tau_{2,j;1,1}, & l_1=0.
\end{array}
\right.
\]
The two formulas for ${\mathcal{L}}_{\bf n}^*$ may be put in one writing
\begin{equation}
\label{eq:C}
{\mathcal{L}}_{\bf n}^* = p_0^* + {p}_1^* \widehat{\sigma}_{1}^* + (-1)^{l_1} p_{l_1} \la
 \sigma_{1}^*, \sigma_2^{(l_1)} \widehat{\ra} - p_2
\la \sigma_{1}^*, \frac{\la s_{3,j},s_{2,l_1+1}  \widehat{\ra}}{ \widehat{s}_{3,j} }
\sigma_2^{(l_1)}\widehat{\ra} +
\end{equation}
\[\sum_{k=3, k\neq j}^{m} \delta_{k,j} p_k \la \sigma_{1}^*, \sigma_2^{(l_1)},f_{1,j,k} \tau_{3,j}, s_{2,l_1+1},s_{3,j} \widehat{\ra},  \qquad n_{l_1} = \min\{n_0,n_1\},
\]
where
\[ \delta_{k,j} = \left\{
\begin{array}{cc}
(-1)^{k-1}, & k < j, \\
(-1)^{j}, & k > j.
\end{array}
\right.
\]
We also have
\[ \deg p_0^* \leq n_j -1 = n_0^* -1, \qquad \deg p_1^* \leq \max\{n_0,n_1\} -1 = n_1^* -1.
\]

For $j=2$ we found the following solutions
\[
\begin{array}{lll}
\underline{\max\{n_0,n_1\}} & \underline{(n_0^*,\ldots,n_m^*)} & \underline{{\mathcal{N}}(\sigma_1^*,\ldots,\sigma_m^*)} \\
n_0 & (n_2,n_0,n_1,n_3,\ldots,n_m) & {\mathcal{N}}(\tau_{1,2},\la \tau_{2,2},s_{1,2}\ra, s_{3,2}, \sigma_4,\ldots,\sigma_m ) \\
n_1 & (n_2,n_1,n_0,n_3,\ldots,n_m) & {\mathcal{N}}(\frac{\widehat{s}_{2,1}}{\widehat{s}_{2,2}} \tau_{1,2},\frac{\widehat{s}_{1,2}}{\widehat{s}_{1,1}} \tau_{2,2;1,1}, s_{3,1}, \sigma_{4},\ldots,\sigma_m)
\end{array}
\]
(If $m=2$, then $ {\bf n}^*$ has only the first three components and
the Nikishin system has only the first two measures indicated.)

We are ready for the induction hypothesis but, for the sake of clarity, let us take one more step.

{\bf Second step.} The case $j=2$ has been solved;
therefore, $j\geq 3$. We ask whether $\min\{n_0,n_1\} \geq n_2$ or
$\min\{n_0,n_1\} \leq n_2$? (When $\min\{n_0,n_1\} = n_2$ we can
proceed either ways.)

\noindent {\bf A2)} If $\min\{n_0,n_1\} \geq n_2$, take
$ {n}_2^* = \min\{n_0,n_1\}$ and $\sigma_2^* = \sigma_2^{(l_1)}, l_1 \in \{0,1\}$. Decompose $\frac{\la
s_{3,j},s_{2,l_1+1} \widehat{\ra}}{ \widehat{s}_{3,j} }$  using
(\ref{2.3}).   Then, the first four terms of ${\mathcal{L}}_{\bf n}^*$ reduce to
\[ p_0^* + {p}_1^* \widehat{\sigma}_{1}^* + (-1)^{l_1}p_{l_1} \la
 \sigma_{1}^*, \sigma_2^* \widehat{\ra} - p_2
\la \sigma_{1}^*, \frac{\la s_{3,j},s_{2,l_1+1}  \widehat{\ra}}{ \widehat{s}_{3,j} }
\sigma_2^*\widehat{\ra} =
\]
\[ p_0^* + {p}_1^* \widehat{\sigma}_{1}^* + p_2^*\la
 \sigma_{1}^*, \sigma_2^* \widehat{\ra} + p_2
\la \sigma_{1}^*,
\sigma_2^*, \tau_{3,j}, s_{2,l_1+1}, s_{3,j}\widehat{\ra}, \qquad n_{l_1} = \min\{n_0,n_1\},
\]
with $  {p}_2^*  = (-1)^{l_1}p_{l_1} - \frac{|\la s_{3,j},s_{2,l_1+1} {\ra}|}{ |{s}_{3,j}| } p_2$, $\deg  {p}_2^* \leq
 {n}_2^* -1.$

In case that $j=3$, we have
\[ {\mathcal{L}}_n^*  =  p_0^* + {p}_1^* \widehat{\sigma}_{1}^* + p_2^*\la
 \sigma_{1}^*, \sigma_2^* \widehat{\ra} + p_2
\la \sigma_{1}^*,
\sigma_2^*, \tau_{3,3}, s_{2,l_1+1}, s_{3,3}\widehat{\ra} - \sum_{k=4}^{m} p_k
\la \sigma_{1}^*,\sigma_2^*,f_{1,3,k} \tau_{3,3}, s_{2,l_1+1},s_{3,3} \widehat{\ra},
\]
This subcase produces the two solutions (in both ${\min\{n_0,n_1\}} \geq n_2$)
\[
\begin{array}{lll}
\underline{\max\{n_0,n_1\}} &  \underline{(n_0^*,\ldots,n_m^*)} & \underline{{\mathcal{N}}(\sigma_1^*,\ldots,\sigma_m^*)} \\
n_0  & (n_3,n_0,n_1,n_2,\ldots)  & {\mathcal{N}}(\tau_{1,3},\widehat{s}_{1,3}\tau_{2,3}, \widehat{s}_{2,3}\tau_{3,3}, s_{4,3},\ldots ) \\
n_1 &  (n_3,n_1,n_0,n_2,\ldots)  & {\mathcal{N}}(\frac{\la {s}_{2,3},\sigma_1 \widehat{\ra}}{\widehat{s}_{2,3}} \tau_{1,3},\frac{\widehat{s}_{1,3}}{\widehat{s}_{1,1}} \tau_{2,3;1,1},\la \tau_{3,3},s_{2,3},\sigma_1 \ra, s_{4,3},\ldots)
\end{array}
\]
If $m=3$, then $ {\bf n}^*$ has only the first four components and
the Nikishin system has only the first three measures indicated. When $m \geq 4$
\[ (n_4^*,\ldots,n_m^*) = (n_4,\ldots,n_m), \qquad (\sigma_5^*,\ldots,\sigma_m^*) = (\sigma_5 ,\ldots,\sigma_m ),\,\, (\mbox{if}\,\, m \geq 5).
\]

Let $j \geq 4$. Define
\[ f  = f_{2,j,k}= \left\{
\begin{array}{ll}
\displaystyle{\frac{\la s_{5,j}, \sigma_4 \widehat{\ra}}{\widehat{s}_{5,j}}}, &  4 = k < j \leq m, \\
\la \displaystyle{\frac{\la
\tau_{5,j},s_{4,j}\ra}{\widehat{s}_{6,j}}},\la s_{6,j},\sigma_5 \ra
\widehat{\ra}, &  5 = k < j \leq m, \\
\la \la \tau_{5,j},s_{4,j}\ra,\ldots,\la \tau_{k-1,j},s_{k-2,j} \ra,
\displaystyle{\frac{\la \tau_{k,j},s_{k-1,j}
\ra}{\widehat{s}_{k+1,j}}} , \la s_{k+1,j},\sigma_k \ra
\widehat{\ra}, &  6 \leq k < j \leq m, \\
\la s_{5,k}, \sigma_4 \widehat{\ra},  & 4 = j < k \leq m, \\
\la \la \tau_{5,j},s_{4,j}\ra,\ldots,\la \tau_{j,j},s_{j-1,j}
\ra,\la s_{j+1,k},\sigma_j \ra  \widehat{\ra}, &  5 \leq j < k \leq
m,
\end{array}
\right.
\]
From (\ref{eq:C}), we obtain
\[ {\mathcal{L}}_n^*  = p_0^* + {p}_1^* \widehat{\sigma}_{1}^* + p_2^*\la
 \sigma_{1}^*, \sigma_2^* \widehat{\ra} + p_2
\la \sigma_{1}^*,
\sigma_2^*, \tau_{3,j}, s_{2,l_1+1}, s_{3,j}\widehat{\ra} + p_3 \la \sigma_{1}^*, \sigma_2^*,\frac{\la {s}_{4,j}, \sigma_3 \widehat{\ra}}{\widehat{s}_{4,j}}\tau_{3,j}, s_{2,l_1+1},s_{3,j} \widehat{\ra} +
\]
\begin{equation}
\label{eq:D}
\sum_{k=4, k\neq j}^{m} \delta_{k,j} p_k \la \sigma_{1}^*, \sigma_2^*,\la \tau_{3,j}, s_{2,l_1+1},s_{3,j} \ra, f_{2,j,k}  \tau_{4,j}, s_{3,j}  \widehat{\ra}, \qquad n_{l_1} = \min\{n_0,n_1\} .
\end{equation}

\noindent {\bf B2)} If $\min\{n_0,n_1\} \leq n_2$, take
$ {n}_2^* = n_2$ and $\sigma_2^* = \frac{\la s_{3,j},s_{2,l_1+1}  \widehat{\ra}}{ \widehat{s}_{3,j} }
\sigma_2^{(l_1)}$. Using   (\ref{2.2}) to decompose $\frac{\widehat{s}_{3,j}}{\la s_{3,j},s_{2,l_1+1}  \widehat{\ra}  }$, the first four terms of ${\mathcal{L}}_{\bf n}^*$ in (\ref{eq:C}) become
\[ p_0^* + {p}_1^* \widehat{\sigma}_{1}^* - p_2
\la \sigma_{1}^*,
\sigma_2^*\widehat{\ra} + (-1)^{l_1}p_{l_1} \la
 \sigma_{1}^*, \frac{\widehat{s}_{3,j}}{\la s_{3,j},s_{2,l_1+1} \widehat{\ra} }\sigma_2^* \widehat{\ra} =
\]
\[p_0^* + {p}_1^* \widehat{\sigma}_{1}^* + p_2^*
\la \sigma_{1}^*,
\sigma_2^*\widehat{\ra} + (-1)^{l_1}p_{l_1} \la
 \sigma_{1}^*, \sigma_2^*, \frac{\la {s}_{2,l_1+1},s_{3,j}\widehat{\ra}}{\widehat{s}_{2,l_1+1} }\tau_{3,j;2,l_1+1} \widehat{\ra}, \qquad n_{l_1} = \min\{n_0,n_1\}
\]
where $p_2^* = - p_2 + (-1)^{l_1}\frac{|\widehat{s}_{3,j}|}{|\la s_{3,j},s_{2,l_1+1} \widehat{\ra} |}p_{l_1}, \deg p_2^* \leq n_2^* -1$.

If $j=m=3$ we are done. Should $m \geq 4$ and $j=3$,
using (\ref{inversa2*}) with $\sigma_{\alpha} = \sigma_3 = s_{3,3},
\sigma_{\beta}=  s_{2,l_1+1}, \sigma_{\gamma}= s_{4,k}$, and $f \equiv 1,$
for the terms in $\sum_{k=4}^m$, we obtain (see definition of
$f_{1,3,k}$)
\[ \la
\sigma_2^{(l_1)},
 f_{1,3,k} \tau_{3,3},s_{2,l_1+1}, s_{3,3} \widehat{\ra}
= \la \frac{\widehat{s}_{3,3}}{\la
s_{3,3},s_{2,l_1+1}\widehat{\ra}}
\sigma_2^*, \la \tau_{3,3},s_{2,l_1+1},s_{3,3} \ra,   s_{4,k},
s_{3,3} \widehat{\ra} =
\]
\[ \la  \sigma_2^*, \frac{\la s_{2,l_1+1}, s_{3,3} \widehat{\ra}}{\widehat{s}_{2,l_1+1}}\tau_{3,3;2,l_1+1}, s_{4,k},
s_{3,l_1+1} \widehat{\ra}.
\]
Therefore,
\[ {\mathcal{L}}_{\bf n}^* = p_0^* + {p}_1^* \widehat{\sigma}_{1}^* + p_2^*
\la \sigma_{1}^*,
\sigma_2^*\widehat{\ra} + (-1)^{l_1}p_{l_1} \la
 \sigma_{1}^*, \sigma_2^*, \frac{\la {s}_{2,l_1+1},s_{3,j}\widehat{\ra}}{\widehat{s}_{2,l_1+1} }\tau_{3,j;2,l_1+1} \widehat{\ra} -
\]
\[  \sum_{k=4}^m p_k \la \sigma_1^*,  \sigma_2^*, \frac{\la s_{2,l_1+1}, s_{3,3} \widehat{\ra}}{\widehat{s}_{2,l_1+1}}\tau_{3,3;2,l_1+1}, s_{4,k},
s_{3,l_1+1} \widehat{\ra}
\]
The case $j=3$ is finished with two additional solutions (in both ${\min\{n_0,n_1\}} \leq n_2$)
\[
\begin{array}{lll}
\underline{\max\{n_0,n_1\}} &  \underline{(n_0^*,\ldots,n_m^*)} & \underline{{\mathcal{N}}(\sigma_1^*,\ldots,\sigma_m^*)} \\
n_0  & (n_3,n_0,n_2,n_1,\ldots)  & {\mathcal{N}}(\tau_{1,3},\frac{\widehat{s}_{3,2}}{\widehat{s}_{3,3}}\widehat{s}_{1,3}\tau_{2,3}, \frac{\widehat{s}_{2,3}}{\widehat{s}_{2,2}}\tau_{3,3;2,2}, s_{4,2},\ldots) \\
n_1 &  (n_3,n_1,n_2,n_0,\ldots)  & {\mathcal{N}}(\frac{\la {s}_{2,3},\sigma_1 \widehat{\ra}}{\widehat{s}_{2,3}} \tau_{1,3},\frac{\widehat{s}_{3,1}}{\widehat{s}_{3,3}}\frac{\widehat{s}_{1,3}}{\widehat{s}_{1,1}} \tau_{2,3;1,1}, \frac{\la s_{2,1},s_{3,3} \widehat{\ra}}{\widehat{s}_{2,1}}\tau_{3,3;2,1},  s_{4,1},\ldots)
\end{array}
\]
If $m=3$, then $ {\bf n}^*$ has only the first four components and
the Nikishin system has only the first three measures indicated. When $m \geq 4$
\[ (n_4^*,\ldots,n_m^*) = (n_4,\ldots,n_m), \qquad (\sigma_5^*,\ldots,\sigma_m^*) = (\sigma_5 ,\ldots,\sigma_m ),\,\, (\mbox{if}\,\, m \geq 5).
\]

Let us transform the terms in $\sum_{k=3, k\neq j}^m$ of (\ref{eq:C}).
Assume that $j \geq 4$. Consider the function multiplying
$p_3$; that is, when $3=k < j \leq m$. Using (\ref{shirenu}) with
$\sigma_{\alpha} = \sigma_3, \sigma_{\beta} = s_{4,j},$ and
$\sigma_{\gamma}  = s_{2,l_1+1} $, we obtain
\[ \la \sigma_2^{(l_1)},f_{1,j,3} \tau_{3,j}, s_{2,l_1+1},s_{3,j} \widehat{\ra} =
 \la \frac{\widehat{s}_{3,j}}{\la {s}_{3,j}, s_{2,l_1+1}\widehat{\ra}} \sigma_2^*, \frac{\la {s}_{4,j}, \sigma_3 \widehat{\ra}}{\widehat{s}_{4,j}} \tau_{3,j}, s_{2,l_1+1}, s_{3,j} \widehat{\ra} =
\]
\[\la \sigma_2^*, \frac{\la s_{4,j},s_{3,l_1+1}
\widehat{\ra}}{\widehat{s}_{4,j}}\frac{\la {s}_{2,l_1+1},s_{3,j} \widehat{\ra}}{\widehat{s}_{2,l_1+1}}\tau_{3,j;2,l_1+1}
\widehat{\ra},\qquad l_1 \in \{0,1\}.
\]

To reduce the terms in $\sum_{k=4}^{j-1}$ when $5 \leq j \leq
m,$ and $\sum_{k=j+1}^{m}$ for $4 \leq j < m$, apply
(\ref{inversa2*}) with $\sigma_{\alpha} = s_{3,j},$ $\sigma_{\beta}
= \sigma_2$ or $\sigma_{\beta} = \la \sigma_2,\sigma_1\ra,$
$\sigma_{\gamma} = \tau_{4,j}$, and $f= f_{2,j,k}$. It follows that
\[ \la \sigma_2^{(l_1)},f_{1,j,k} \tau_{3,j}, s_{2,l_1+1},s_{3,j} \widehat{\ra} = \la \frac{\widehat{s}_{3,j} }{ \la s_{3,j},s_{2,l_1+1}  \widehat{\ra}}\sigma_2^*, \la\tau_{3,j},s_{2,l_1+1},s_{3,j} \ra, f_{2,j,k} \tau_{4,j}, s_{3,j} \widehat{\ra} =
\]
\[ \la \sigma_2^*, \frac{\la {s}_{2,l_1+1}, s_{3,j}\widehat{\ra}}{\widehat{s}_{2,l_1+1}} \tau_{3,j;2,l_1+1}, f_{2,j,k} \tau_{4,j},s_{3,j}, s_{2,l_1+1} \widehat{\ra}, \qquad l_1 \in \{0,1\}.
\]

Therefore,
\begin{equation}
\label{eq:E}
{\mathcal{L}}_{\bf n}^* = p_0^* + {p}_1^* \widehat{\sigma}_{1}^* + p_2^*
\la \sigma_{1}^*,
\sigma_2^*\widehat{\ra} + (-1)^{l_1}p_{l_1} \la
 \sigma_{1}^*, \sigma_2^*, \frac{\la {s}_{2,l_1+1},s_{3,j}\widehat{\ra}}{\widehat{s}_{2,l_1+1} }\tau_{3,j;2,l_1+1} \widehat{\ra} +
\end{equation}
\[ p_3 \la \sigma_1^*, \sigma_2^*, \frac{\la s_{4,j},s_{3,l_1+1}
\widehat{\ra}}{\widehat{s}_{4,j}}\frac{\la {s}_{2,l_1+1},s_{3,j} \widehat{\ra}}{\widehat{s}_{2,l_1+1}}\tau_{3,j;2,l_1+1}
\widehat{\ra},
\]
\[ \sum_{k=4, k\neq j} \delta_{k,j} p_k \la \sigma_1^*, \sigma_2^*,  \frac{\la {s}_{2,l_1+1}, s_{3,j}\widehat{\ra}}{\widehat{s}_{2,l_1+1}} \tau_{3,j;2,l_1+1}, f_{2,j,k} \tau_{4,j},s_{3,j}, s_{2,l_1+1} \widehat{\ra}, \qquad l_1 \in \{0,1\}.
\]

Let us write down (\ref{eq:D})-(\ref{eq:E}) in one expression. Recall that $j \geq 4$. Take
\[  {\sigma}_3^{(l_2)} =
\left\{
\begin{array}{ll}
\widehat{s}_{2,j} \tau_{3,j}, &  l_1=1 \,\,\mbox{in}\,\, (\ref{eq:D}),\\
\la s_{2,j},\sigma_1 \widehat{\ra} \tau_{3,j}, &  l_1=0 \,\,\mbox{in} \,\, (\ref{eq:D}),\\
\frac{\widehat{s}_{2,j}}{\widehat{s}_{2,2}} \tau_{3,j;2,2}, &  l_1=1 \,\,\mbox{in}\,\, (\ref{eq:E}), \\
\frac{\la s_{2,1},s_{3,j}\widehat{\ra}}{\widehat{s}_{2,1}} \tau_{3,j;2,1},  & l_1=0\,\,\mbox{in} \,\, (\ref{eq:E}) .
\end{array}
\right.
\]
Then
\[ {\mathcal{L}}_n^*  = p_0^* + {p}_1^* \widehat{\sigma}_{1}^* + p_2^*\la
 \sigma_{1}^*, \sigma_2^* \widehat{\ra} + (-1)^{l_2}p_{l_2}
\la \sigma_{1}^*,
\sigma_2^*, \sigma_3^{(l_2)}\widehat{\ra} + p_3 \la \sigma_{1}^*, \sigma_2^*,\frac{\la {s}_{4,j}, s_{3,l_2+1} \widehat{\ra}}{\widehat{s}_{4,j}}\sigma_3^{(l_2)}\widehat{\ra} +
\]
\begin{equation}
\label{eq:F}
\sum_{k=4, k\neq j}^{m} \delta_{k,j} p_k \la \sigma_{1}^*, \sigma_2^*,\sigma_3^{(l_2)}, f_{2,j,k}  \tau_{4,j}, s_{3,l_2+1},s_{4,j}  \widehat{\ra}, \qquad n_{l_2} = \min\{n_0,n_1,n_2\} .
\end{equation}
We also have that $\deg p_k^* \leq n_k^*-1, k=0,1,2,$ where
\begin{equation} \label{eq:G}
{n}_k^* =
\left\{
\begin{array}{ll}
n_j = \max\{n_0+1,n_1,\ldots,n_m\}, & k=0, \\
\max\{n_0,n_1\}, & k =1,  \\
\max\{\min\{n_0,n_1\},n_2\}, & k= 2,
\end{array}
\right.
\end{equation}
and the measures $\sigma_1^*,\sigma_2^*$ have also been determined. The polynomials $p_0^*,{p}_1^*,{p}_2^*$, the indices $n_0^*,{n}_1^*,  {n}_2^*$, and the measures ${\sigma}_1^*, {\sigma}_2^*$ will not change in subsequent reductions.   The structure of a Nikishin systems brakes down in the Cauchy transform multiplying $p_3$ due to an annoying ratio of Cauchy transforms modifying the third measure in the corresponding product.

With this unified formula, you can repeat the line of reasonings employed in step 2 (or 1) and so on. Let us write down the eight solutions corresponding to $j=4$. In the table, besides (\ref{eq:G}), we have that $n_3^* = \max\{\min\{n_0,n_1,n_2\},n_3\}$,  $n_4^* = \min\{n_0,\ldots,n_3\}$, and of course $n_0^* = n_4$.
\[
\begin{array}{ll}
 \underline{(n_1^*,\ldots,n_4^*)} & \underline{{\mathcal{N}}(\sigma_1^*,\ldots,\sigma_m^*)} \\
 ( n_0, n_1, n_2, n_3 )  & {\mathcal{N}}(\tau_{1,4},\widehat{s}_{1,4} \tau_{2,4}, \widehat{s}_{2,4}  \tau_{3,4},  \widehat{s}_{3,4}  \tau_{4,4}, s_{5,4},\ldots) \\
  ( n_1, n_0, n_2, n_3 )  & {\mathcal{N}}(\frac{\la s_{2,4},\sigma_1\widehat{\ra}}{\widehat{s}_{2,4}} \tau_{1,4},\frac{\widehat{s}_{1,4}}{\widehat{\sigma}_1}\tau_{2,4;1,1}, \la \tau_{3,4}, {s}_{2,4},\sigma_1 \ra,  \widehat{s}_{3,4}\tau_{4,4}, s_{5,4},\ldots) \\
 ( n_0, n_2, n_1, n_3 ) & {\mathcal{N}}(\tau_{1,4}, \frac{\la  {s}_{3,4},\sigma_2 \widehat{\ra}}{\widehat{s}_{3,4}} \widehat{s}_{1,4} \tau_{2,4},
\frac{\widehat{s}_{2,4}}{\widehat{\sigma}_2} \tau_{3,4;2,2},  \la \tau_{4,4}, s_{3,2}, \sigma_4\ra, s_{5,4},\ldots) \\
( n_1, n_2, n_0, n_3) & {\mathcal{N}}(\frac{\la s_{2,4},\sigma_1\widehat{\ra}}{\widehat{s}_{2,4}} \tau_{1,4},\frac{\la  {s}_{3,4},s_{2,1} \widehat{\ra}}{\widehat{s}_{3,4}} \frac{\widehat{s}_{1,4}}{\widehat{\sigma}_1} \tau_{2,4;1,1},
\frac{\la s_{2,1},{s}_{3,4}\widehat{\ra}}{ \widehat{s}_{2,1}}  \tau_{3,4;2,1}, \la \tau_{4,4}, s_{3,1}, \sigma_4\ra ,  s_{5,4},\ldots) \\
( n_0, n_1, n_3, n_2 )  & {\mathcal{N}}(\tau_{1,4} ,\widehat{s}_{1,4} \tau_{2,4}, \frac{\widehat{s}_{4,3} }{\widehat{\sigma}_{4}}\widehat{s}_{2,4}  \tau_{3,4},   \frac{\widehat{s}_{3,4}}{\widehat{\sigma}_3} \tau_{4,4;3,3}, s_{5,3},\ldots) \\
( n_1, n_0, n_3, n_2 )& {\mathcal{N}}(\frac{\la s_{2,4},\sigma_1\widehat{\ra}}{\widehat{s}_{2,4}} \tau_{1,4},\frac{\widehat{s}_{1,4}}{\widehat{\sigma}_1}\tau_{2,4;1,1}, \frac{\widehat{s}_{4,3} }{\widehat{\sigma}_{4}}\la \tau_{3,4}, {s}_{2,4},\sigma_1 \ra,  \frac{\widehat{s}_{3,4}}{\widehat{\sigma}_3} \tau_{4,4;3,3}, s_{5,3},\ldots) \\
( n_0, n_2, n_3, n_1 ) & {\mathcal{N}}(\tau_{1,4},\frac{\la  {s}_{3,4},\sigma_2 \widehat{\ra}}{\widehat{s}_{3,4}} \widehat{s}_{1,4} \tau_{2,4},
\frac{\widehat{s}_{4,2}}{\widehat{\sigma}_{4}} \frac{\widehat{s}_{2,4}}{\widehat{\sigma}_2} \tau_{3,4;2,2}, \frac{\la s_{3,2},\sigma_4\widehat{\ra}}{\widehat{s}_{3,2}}\tau_{4,4;3,2}  , s_{5,2},\ldots) \\
( n_1, n_2, n_3, n_0 )& {\mathcal{N}}(\frac{\la s_{2,4},\sigma_1\widehat{\ra}}{\widehat{s}_{2,4}} \tau_{1,4},\frac{\la  {s}_{3,4},s_{2,1} \widehat{\ra}}{\widehat{s}_{3,4}} \frac{\widehat{s}_{1,4}}{\widehat{\sigma}_1} \tau_{2,4;1,1}, \frac{\widehat{s}_{4,1}}{\widehat{\sigma}_{4}}
\frac{\la s_{2,1},{s}_{3,4}\widehat{\ra}}{ \widehat{s}_{2,1}}  \tau_{3,4;2,1}, \frac{\la s_{3,1},\sigma_4\widehat{\ra}}{\widehat{s}_{3,1}}\tau_{4,4;3,1} ,  s_{5,1},\ldots)
\end{array}
\]
If $m=4$, then $ {\bf n}^*$ has only the first five components and
the Nikishin system has only the first four measures indicated. When $m \geq 5$
\[ (n_5^*,\ldots,n_m^*) = (n_5,\ldots,n_m), \qquad (\sigma_6^*,\ldots,\sigma_m^*) = (\sigma_6 ,\ldots,\sigma_m ),\,\, (\mbox{if}\,\, m \geq 6).
\]

Summarizing, in step 1 we proved the statement of Lemma \ref{lem:4} when $j=2$, showing that there are two solutions,  and  for $j\geq 3$ we obtained formula  (\ref{eq:C}) which allowed us to prove in step 2 that the lemma is true when $j=3$, with 4 solutions, and for $j \geq 4$ to obtain formula (\ref{eq:F}), similar to (\ref{eq:C}), which provides the instruments to carry out step 3 and so on. The case $j=1$ was treated separately. The induction will be on the number of successful steps we have been able to carry out. The counter of the steps will be denoted $j^*$. Fix $j, 2 \leq j \leq m$. Assume that we have been able to carry out $j^*$ steps. Let us describe what we have obtained in

{\bf Step $j^*$. Induction hypothesis.} We have proved that the statement of Lemma $\ref{lem:4}$ holds when $j = j^*+1 $, with $2^{j-1} = 2^{j^*}$ solutions, and when $j \geq j^* +2,$ we have defined:
\begin{enumerate}
\item indices
\[ {n}_k^* =
\left\{
\begin{array}{ll}
n_j = \max\{n_0 +1,n_1,\ldots,n_m\}, & k=0, \\
\max\{\min\{n_0,\ldots,n_{k-1}\},n_k\}, & k =1,\ldots,j^*.
\end{array}
\right.
\]
\item integers $l_k, k=0,\ldots,j^*$, inductively as follows: $l_0 =0$, $l_1$ is the subindex between those of $n_0,n_1$ not employed in defining $n_1^*$, and so forth until $l_{j^*}$ which is the subindex between those of $n_0,\ldots,n_{j^*}$ which was not employed in defining $n_1^*,\ldots,n_{j^*}^*$. In particular,
\[ n_{l_k} = \min\{n_0,\ldots,n_k\}.
\]
\item polynomials
\[ p_k^* =
\left\{
\begin{array}{lll}
\ell_{1,j}p_0  +  p_j  + \sum_{i=1, i\neq j}^m \frac{|s_{1,i}|}{|s_{1,j}|} p_i, & k=0, & \mbox{} \\
(-1)^{l_{k-1}}p_{l_{k-1}} + C_{1,j,k}p_k, & k=1,\ldots,j^*, &  \min\{n_0,\ldots,n_{k-1}\} \geq n_k, \\
(-1)^{k-1}p_{k} + C_{2,j,k}p_{l_{k-1}}, & k=1,\ldots,j^*, &   \min\{n_0,\ldots,n_{k-1}\} \leq n_k.
\end{array}
\right.
\]
where $C_{1,j,k}, C_{2,j,k},$ are real constants different from zero.
\item functions
\[ f_{j^*,j,k}= \left\{
\begin{array}{ll}
\la s_{j^*+3,k}, \sigma_{j^*+2} \widehat{\ra},  & j^* +2 = j < k \leq m, \\
\la \la \tau_{j^*+3,j},s_{j^*+2,j}\ra,\ldots,\la \tau_{j,j},s_{j-1,j}
\ra,\la s_{j+1,k},\sigma_j \ra  \widehat{\ra}, & j^*+3 \leq j < k \leq
m, \\
\displaystyle{\frac{\la s_{j^*+3,j}, \sigma_{j^*+2} \widehat{\ra}}{\widehat{s}_{j^*+3,j}}}, & j^*+2 = k < j \leq m, \\
\la \displaystyle{\frac{\la
\tau_{j^*+3,j},s_{j^*+2,j}\ra}{\widehat{s}_{j^*+4,j}}},\la s_{j^*+4,j},\sigma_{j^*+3} \ra
\widehat{\ra}, &  j^*+3 = k < j \leq m, \\
\la \widehat{s}_{j^*+2,j} \tau_{j^*+3,j},\ldots,\widehat{s}_{k-2,j}\tau_{k-1,j},
\displaystyle{\frac{\la \tau_{k,j},s_{k-1,j}
\ra}{\widehat{s}_{k+1,j}}} , \la s_{k+1,j},\sigma_k \ra
\widehat{\ra}, &  j^*+4 \leq k < j \leq m.
\end{array}
\right.
\]
\item measures $\sigma_1^*,\ldots,\sigma_{j^*}^*$ and $\sigma_{j^*+1}^{(l_{j^*})}$ whose supports are contained in the same intervals $\Delta_1,\ldots,$ $\Delta_{j^*+1}$ as $\sigma_1,\ldots,\sigma_{j^* +1}$, respectively, where
    \[ \sigma_{j^*+1}^{(l_{j^*})} = \left\{
\begin{array}{ll}
\tau_{j^*+1,j}, s_{j^*,l_{j^*-1}+1}, s_{j^*+1,j}, & \mbox{if}\quad \max\{\min\{n_0,\ldots,n_{j^*-1}\}, n_{j^*} \} = n_{l_{j^*-1}}, \\
\displaystyle{\frac{\la s_{j^*,l_{j^*-1}+1},s_{j^*+1,j}\widehat{\ra}}{\widehat{s}_{j^*,l_{j^*-1}+1}}} \tau_{j^*+1,j;j^*,l_{j^*-1}+1}, & \mbox{if}\quad \max\{\min\{n_0,\ldots,n_{j^*-1}\}, n_{j^*} \} = n_{j^*}.
\end{array}
\right.
    \]
\end{enumerate}
With these elements, we have proved the formula (analogous to those obtained in steps 1 and 2)
\begin{equation}
\label{eq:H}
{\mathcal{L}}_{\bf n}^* = p^*_0 + \sum_{k=1}^{j^*} p^*_k \la \sigma_1^* ,\ldots,\sigma^*_k \widehat{\ra} + (-1)^{l_{j^*}}p_{l_{j^*}} \la\sigma_1^*, \ldots, \sigma^*_{j^*} , \sigma_{j^*+1}^{(l_{j^*})}  \widehat{\ra} +
\end{equation}
\[ (-1)^{j^*} p_{j^*+1} \la\sigma_1^*, \ldots, \sigma^*_{j^*} , \frac{\la s_{j^* +2,j}, s_{j^* +  1,l_{j^*}+1} \widehat{\ra}}{\widehat{s}_{j^* +2,j}} \sigma_{j^*+1}^{(l_{j^*})}  \widehat{\ra} +
\]
\[ \sum_{k = j^* +2 , k \neq j}^m \delta_{k,j} p_k \la \sigma_1^*, \ldots, \sigma^*_{j^*} , \sigma_{j^*+1}^{(l_{j^*})} ,f_{j^*,j,k}\tau_{j^*+2,j}, s_{j^*+1,l_{j^*}+1}, s_{j^*+2,j} \widehat{\ra}.
\]

{\bf Step $j^* +1$. Induction proof.} To complete the induction we must prove in this step that Lemma \ref{lem:4} is satisfied when $j=j^* +2$, with $2^{j^* +1}$ solutions, and when $j \geq j^* +3$ to produce a formula which extends (\ref{eq:H}) one more step.

\noindent {\bf Case A)} Suppose that $\max\{\min\{n_0,\ldots,n_{j^*}\}, n_{j^*+1} \} = \min\{n_0,\ldots,n_{j^*}\} = n_{l_{j^*}}$. Take
\[ {n}_{j^* +1}^*  = n_{l_{j^*}} \quad \mbox{and} \quad\sigma^*_{j^*+1} = \sigma_{j^*+1}^{(l_{j^*})}.
\]
Using (\ref{2.3}) on $\frac{\la s_{j^* +2,j}, s_{j^* +  1,l_{j^*}+1} \widehat{\ra}}{\widehat{s}_{j^* +2,j}}$, from (\ref{eq:H}) it follows that
\[
{\mathcal{L}}_{\bf n}^* = p^*_0 + \sum_{k=1}^{j^*+1} p^*_k \la \sigma_1^* ,\ldots,\sigma^*_k \widehat{\ra} +
(-1)^{j^*+1} p_{j^*+1} \la\sigma_1^*, \ldots, \sigma^*_{j^*} , \sigma^*_{j^* +1} , \tau_{j^*+2,j}, s_{j^*+1,l_{j^*}+1}, s_{j^*+2,j} \widehat{\ra} +
\]
\[
\sum_{k = j^* +2 , k \neq j}^m \delta_{k,j} p_k \la \sigma_1^*, \ldots, \sigma^*_{j^*} , \sigma_{j^*+1}^{(l_{j^*})} ,f_{j^*,j,k}\tau_{j^*+2,j}, s_{j^*+1,l_{j^*}+1}, s_{j^*+2,j} \widehat{\ra} ,
\]
where
\[p^*_{j^*+1} = (-1)^{l_{j^*}}p_{l_{j^*}} + (-1)^{j^*} \frac{|\la s_{j^*+2,j}, s_{j^*+1,l_{j^*}+1}\ra|}{\widehat{s}_{j^*+2,j}} p_{j^*+1}, \qquad \deg p^*_{j^*+1} \leq n_{j^*+1} -1.\]

Since $f_{j^*,j^*+2,k} = \la s_{j^*+3,k}, \sigma_{j^*+2} \widehat{\ra}$, if $j=j^*+2$ we have that $ {\mathcal{L}}_{\bf n}^*$ is a linear form generated by
\[ {\bf n}^* = ( {n}_0^*,\ldots, {n}_{j^*+1}^*,n_{j^*+1},n_{j^*+3},\ldots,n_m)
\]
and the Nikishin system
\[ {\mathcal{N}}(\sigma_1^* , \ldots, \sigma^*_{j^*+1} ,   \la \tau_{j^*+2,j^*+2}, s_{j^*+1,l_{j^*}+1}, s_{j^*+2,j^*+2}\ra,  s_{j^*+3,j^*+2}, \sigma_{j^*+4} , \ldots,\sigma_m).
\]
When $m = j (= j^*+2)$, the system ends with the measure $\la \tau_{j^*+2,j^*+2}, s_{j^*+1,l+1}, s_{j^*+2,j^*+2}\ra$.

If $j \geq j^*+3$, (\ref{eq:H}) may be expressed as
\[
{\mathcal{L}}_{\bf n}^* = p^*_0 + \sum_{k=1}^{j^*+1} p^*_k \la \sigma_1^* ,\ldots,\sigma^*_k \widehat{\ra} +
(-1)^{j^*+1} p_{j^*+1} \la\sigma_1^*, \ldots,   \sigma^*_{j^* +1} , \tau_{j^*+2,j}, s_{j^*+1,l_{j^*}+1}, s_{j^*+2,j} \widehat{\ra} +
\]
\begin{equation}
\label{eq:I}
(-1)^{j^*+1} p_{j^*+2} \la \sigma_1^*, \ldots,  \sigma^*_{j^* +1}  ,\displaystyle{\frac{\la s_{j^*+3,j}, \sigma_{j^*+2} \widehat{\ra}}{\widehat{s}_{j^*+3,j}}}\tau_{j^*+2,j}, s_{j^*+1,l_{j^*}+1}, s_{j^*+2,j} \widehat{\ra} +
\end{equation}
\[
\sum_{k = j^* +3 , k \neq j}^m \delta_{k,j} p_k \la \sigma_1^*, \ldots,  \sigma_{j^*+1}^*, \la \tau_{j^*+2,j}, s_{j^*+1,l_{j^*}+1}, s_{j^*+2,j} \ra,f_{j^*+1,j,k} \tau_{j^*+3,j}, s_{j^*+2,j} \widehat{\ra} ,
\]
where $f_{j^*+1,j,k}$ is defined as $f_{j^*,j,k}$ substituting $j^*$ by $j^*+1$.

\noindent
{\bf Case B)}
If $\max\{\min\{n_0,\ldots,n_{j^*}\}, n_{j^*+1}\} = n_{j^*+1} $, take
\[  {n}_{j^*+1}^* = n_{j^*+1} \quad \mbox{and} \quad \sigma^*_{j^*+1} = \frac{\la s_{j^* +2,j}, s_{j^* +  1,l_{j^*}+1} \widehat{\ra}}{\widehat{s}_{j^* +2,j}} \sigma_{j^*+1}^{(l_{j^*})}.   \]
Rewrite (\ref{eq:H}) as follows
\[
{\mathcal{L}}_{\bf n}^* = p^*_0 + \sum_{k=1}^{j^*} p^*_k \la \sigma_1^* ,\ldots,\sigma^*_k \widehat{\ra} + (-1)^{j^*} p_{j^*+1} \la\sigma_1^*, \ldots, \sigma^*_{j^*}, \sigma^*_{j^*+1}  \widehat{\ra} +
\]
\[ (-1)^{l_{j^*}}p_{l_{j^*}} \la\sigma_1^*, \ldots, \sigma^*_{j^*} , \frac{\widehat{s}_{j^* +2,j}}{\la s_{j^* +2,j}, s_{j^* + 1,l_{j^*}+1} \widehat{\ra}} \sigma_{j^*+1}^*  \widehat{\ra} +
\]
\[ \sum_{k = j^* +2 , k \neq j}^m \delta_{k,j} p_k \la \sigma_1^*, \ldots, \sigma^*_{j^*} , \frac{\widehat{s}_{j^* +2,j}}{\la s_{j^* +2,j}, s_{j^* + 1,l_{j^*}+1} \widehat{\ra}} \sigma_{j^*+1}^*  ,f_{j^*,j,k}\tau_{j^*+2,j}, s_{j^*+1,l_{j^*}+1}, s_{j^*+2,j} \widehat{\ra}.
\]
Using (\ref{2.2}), reduce $\frac{\widehat{s}_{j^* +2,j}}{\la s_{j^* +2,j}, s_{j^* +  1,l_{j^*}+1} \widehat{\ra}}$ in the term with $p_{l_{j^*}}$. This formula transforms into
\begin{equation}
\label{eq:J}
{\mathcal{L}}_{\bf n}^* = p^*_0 + \sum_{k=1}^{j^*+1} p^*_k \la \sigma_1^* ,\ldots,\sigma^*_k \widehat{\ra} +
\end{equation}
\[ (-1)^{l_{j^*}}p_{l_{j^*}} \la\sigma_1^*, \ldots, \sigma^*_{j^*+1} ,  \frac{\la s_{j^*+1,l_{j^*}+1},s_{j^*+2,j}\widehat{\ra}}{\widehat{s}_{j^*+1,l_{j^*}+1}} \tau_{j^*+2,j;j^*+1,l_{j^*}+1}  \widehat{\ra} +
\]
\[ \sum_{k = j^* +2 , k \neq j}^m \delta_{k,j} p_k \la \sigma_1^*, \ldots, \sigma^*_{j^*} , \frac{\widehat{s}_{j^* +2,j}}{\la s_{j^* +2,j}, s_{j^* + 1,l_{j^*}+1} \widehat{\ra}} \sigma_{j^*+1}^*  ,f_{j^*,j,k}\tau_{j^*+2,j}, s_{j^*+1,l_{j^*}+1}, s_{j^*+2,j} \widehat{\ra}.
\]
where
\[ p^*_{j^*+1} = (-1)^{j^*}p_{j^*+1} + (-1)^{l_{j^*}} \frac{|{s}_{j^* +2,j}|}{|\la s_{j^* +2,j}, s_{j^* +  1,l_{j^*}+1}{\ra}|}  p_{l_{j^*}}, \qquad \deg p^*_{j^*+1} \leq  {n}_{j^*+1}-1.
\]
If $m=j^*+2$ we are done, since the sum $\sum_{k = j^*+2, k \neq j}^m$ is empty (recall that $j^* +2 \leq j\leq m$).

Suppose that $j^* +2  = j < m$. Using (\ref{inversa2*}) with $\sigma_{\alpha} = \sigma_{j^*+2},
\sigma_{\beta}= s_{j^*+1,l_{j^*}+1}, \sigma_{\gamma}= s_{j^*+3,k}$, and $f \equiv 1,$
for the terms in $\sum_{k=j^*+3}^m$ (see definition of
$f_{j^*,j^*+2,k}$), (\ref{eq:J}) becomes
\[
{\mathcal{L}}_{\bf n}^* = p^*_0 + \sum_{k=1}^{j^*+1} p^*_k \la \sigma_1^* ,\ldots,\sigma^*_k \widehat{\ra} +
\]
\[ (-1)^{l_{j^*}}p_{l_{j^*}} \la\sigma_1^*, \ldots, \sigma^*_{j^*+1} ,  \frac{\la s_{j^*+1,l_{j^*}+1},\sigma_{j^*+2}\widehat{\ra}}{\widehat{s}_{j^*+1,l_{j^*}+1}} \tau_{j^*+2,j^*+2;j^*+1,l_{j^*}+1}  \widehat{\ra} +
\]
\[ \sum_{k = j^* +3 }^m \delta_{k,j} p_k \la \sigma_1^*, \ldots,   \sigma^*_{j^*+1}, \frac{\la s_{j^*+1,l_{j^*}+1},\sigma_{j^*+2}\widehat{\ra}}{\widehat{s}_{j^*+1,l_{j^*}+1}}
\tau_{j^*+2,j^*+2;j^*+1,l_{j^*}+1}, s_{j^*+3,k},s_{j^*+2,l_{j^*}+1} \widehat{\ra} .
\]
Thus, we conclude with $j=j^*+2$ taking
\[ {\bf n}^* =
( {n}_0^*,\ldots, {n}_{j^*+1}^*,n_{l_{j^*}},n_{j^*+3},\ldots,n_m)
\]
and the Nikishin system
\[
 {\mathcal{N}} (\sigma_1^*, \ldots,   \sigma^*_{j^*+1}, \frac{\la s_{j^*+1,l_{j^*}+1},\sigma_{j^*+2}\widehat{\ra}}{\widehat{s}_{j^*+1,l_{j^*}+1}}
\tau_{j^*+2,j^*+2;j^*+1,l_{j^*}+1}, s_{j^*+3,l_{j^*}+1},\sigma_{j^*+4}, \ldots,\sigma_m).
\]
If $m = j (= j^*+2)$, the system ends with the measure $\frac{\la s_{j^*+1,l_{j^*}+1},\sigma_{j^*+2}\widehat{\ra}}{\widehat{s}_{j^*+1,l_{J^*}+1}}
\tau_{j^*+2,j^*+2;j^*+1,l_{j^*}+1}$.

Let $j^*+3 \leq j \leq m$. In (\ref{eq:J}), the measure multiplying $p_{j^* +2}$  is transformed by means of (\ref{shirenu}), with $\sigma_{\alpha} = \sigma_{j^* +2}, \sigma_{\beta} = s_{j^*+3,j}$ and $\sigma_{\gamma} = s_{j^*+1,l_{j^*}+1}$. All other terms of $\sum_{k=j^*+2, k\neq j}^m$  are reduced employing (\ref{inversa2*}) taking $\sigma_{\alpha} = s_{j^*+2,j}, \sigma_{\beta} = s_{j^*+1,l_{j^*}+1}, \sigma_{\gamma} = \tau_{j^*+2,j}$ and
$f = f_{j^*+1,j,k}$ (remember that $f_{j^*+1,j,k}$ is defined substituting $j^*$ by $j^* +1$ in the definition of $f_{j^*,j,k}$ as was already used at the end of case A)). It is easy to verify that (\ref{eq:J}) adopts the expression
\begin{equation}
\label{eq:K}
{\mathcal{L}}_{\bf n}^* = p^*_0 + \sum_{k=1}^{j^*+1} p^*_k \la \sigma_1^* ,\ldots,\sigma^*_k \widehat{\ra} +
\end{equation}
\[ (-1)^{l_{j^*}}p_{l_{j^*}} \la\sigma_1^*, \ldots, \sigma^*_{j^*+1} ,  \frac{\la s_{j^*+1,l_{j^*}+1},s_{j^*+2,j}\widehat{\ra}}{\widehat{s}_{j^*+1,l_{j^*}+1}} \tau_{j^*+2,j;j^*+1,l_{j^*}+1}  \widehat{\ra} +
\]
\[(-1)^{j^*+1} p_{j^*+2} \la \sigma_1^*, \ldots, \sigma^*_{j^*+1} , \frac{\la {s}_{j^*+3,j}, s_{j^*+2,l_{j^*}+1} \widehat{\ra}}{\widehat{s}_{j^*+3,j}} \frac{\la s_{j^*+1,l_{j^*}+1},s_{j^*+2,j}\widehat{\ra}}{\widehat{s}_{j^*+1,l_{j^*}+1}} \tau_{j^*+2,j;j^*+1,l_{j^*}+1}  \widehat{\ra} +
\]
\[  \sum' \delta_{k,j} p_k \la \sigma_1^*, \ldots, \sigma^*_{j^*+1} ,  \frac{\la s_{j^*+1,l_{j^*}+1},s_{j^*+2,j}\widehat{\ra}}{\widehat{s}_{j^*+1,l_{j^*}+1}} \tau_{j^*+2,j;j^*+1,l_{j^*}+1}  ,f_{j^*+1,j,k}\tau_{j^*+3,j}, s_{j^*+2,j}, s_{j^*+1,l_{j^*}+1} \widehat{\ra}
\]
($\sum' = \sum_{k = j^* +3 , k \neq j}^m$), which has the same structure as (\ref{eq:I}).

Let $l_{j^*+1}$  denote the subindex between those of $n_0,\ldots,n_{j^*+1}$ which was not employed in defining $n_1^*,\ldots,n_{j^*+1}^*$, then
\[ n_{l_{j^* +1}} = \min\{n_0,\ldots,n_{j^* +1}\}.
\]
Notice that in the situation of case A), $l_{j^*+1} = l_{j^*}$, and in case B), $l_{j^*+1} = j^*+1.$ Define
\[ \sigma_{j^*+2}^{(l_{j^*+1})} = \left\{
\begin{array}{ll}
\tau_{j^*+2,j}, s_{j^*+1,l_{j^*}+1}, s_{j^*+2,j}, & \mbox{if}\,\, \max\{\min\{n_0,\ldots,n_{j^*}\}, n_{j^*+1} \} = n_{l_{j^*}}, \\
\displaystyle{\frac{\la s_{j^*+1,l_{j^*}+1},s_{j^*+2,j}\widehat{\ra}}{\widehat{s}_{j^*+1,l_{j^*}+1}}} \tau_{j^*+2,j;j^*+1,l_{j^*}+1}, & \mbox{if}\,\, \max\{\min\{n_0,\ldots,n_{j^*}\}, n_{j^*+1} \} = n_{j^*+1}.
\end{array}
\right.
\]
With these notations, formulas (\ref{eq:I}) and (\ref{eq:K}) may be unified in
\[
{\mathcal{L}}_{\bf n}^* = p^*_0 + \sum_{k=1}^{j^*+1} p^*_k \la \sigma_1^* ,\ldots,\sigma^*_k \widehat{\ra} +
(-1)^{l_{j^*+1}}p_{l_{j^*+1}} \la\sigma_1^*, \ldots, \sigma^*_{j^*+1} ,  \sigma_{j^*+2}^{(l_{j^*+1})} \widehat{\ra} +
\]
\[(-1)^{j^*+1} p_{j^*+2} \la \sigma_1^*, \ldots, \sigma^*_{j^*+1} , \frac{\la {s}_{j^*+3,j}, s_{j^*+2,l_{j^*+1}+1} \widehat{\ra}}{\widehat{s}_{j^*+3,j}} \sigma_{j^*+2}^{(l_{j^*+1})}  \widehat{\ra} +
\]
\[ \sum_{k = j^* +3 , k \neq j}^m \delta_{k,j} p_k \la \sigma_1^*, \ldots, \sigma^*_{j^*+1} ,  \sigma_{j^*+2}^{(l_{j^*+1})} ,f_{j^*+1,j,k}\tau_{j^*+3,j}, s_{j^*+2,l_{j^*+1}+1}, s_{j^*+3,j} \widehat{\ra}.
\]
With this we conclude the induction and Lemma \ref{lem:4} has been proved. \hfill $\Box$ \vspace{0,2cm}

Before moving on, let us write down the expressions of the $p_k^*$ after carrying out $j-1$ steps. We will need their structure  for further developments. We have
\[
{\mathcal{L}}_{\bf n}^* = p^*_0 + \sum_{k=1}^{m} p^*_k \la \sigma_1^* ,\ldots,\sigma^*_k \widehat{\ra}
\]
where
\begin{equation}
\label{eq:M}
p_k^* = \left\{
\begin{array}{lll}
\ell_{1,j}p_0  +  p_j  + \sum_{i=1, i\neq j}^m \frac{|s_{1,i}|}{|s_{1,j}|} p_i, & k=0, & \mbox{} \\
(-1)^{l_{k-1}}p_{l_{k-1}} + C_{1,j,k}p_k, & k=1,\ldots,j-1, &  \max\{{n_{l_{k-1}},n_k}\}  = n_{l_{k-1}}, \\
(-1)^{k-1}p_{k} + C_{2,j,k}p_{l_{k-1}}, & k=1,\ldots,j-1, &   \max\{{n_{l_{k-1}},n_k}\} = n_k, \\
(-1)^{j-1} p_{j-1} & k=j, & \min\{n_{l_{j-2}},n_{j-1}\} = n_{j-1}, \\
(-1)^{l_{j-1}}p_{l_{j-2}}, & k = j, & \min\{n_{l_{j-2}},n_{j-1}\} = n_{l_{j-2}}, \\
(-1)^{j}p_k, & k=j+1,\ldots,m. & \mbox{}
\end{array}
\right.
\end{equation}

Lemma \ref{lem:4} has an immediate consequence in terms of the orthogonality conditions satisfied by the linear form ${\mathcal{A}}_{{\bf n}}$ (see (\ref{orto})). We state it as a lemma which is useful to prove Corollary \ref{cor:2} and the results on the asymptotic behavior of sequences of these linear forms.

\bl \label{lem:ortox} Let  $\widehat{\mathbb{S}}$ and ${\bf n} = ({\bf n}_1;{\bf n}_2) \in {\mathbb{Z}}_+^{m_1 +1}\times {\mathbb{Z}}_+^{m_2 +1}, |{\bf n}_1| = |{\bf n}_2| +1,$ be given. Suppose that $n_{2,j} = \max\{n_{2,0}+1,n_{2,1},\ldots,n_{2,m_2}\}$. Let ${\bf n}_2^* =(n_{2,0}^*,\ldots,n_{2,m_2}^*)$ and $ {\mathcal{N}}(\sigma_1^{2*},\ldots,\sigma_{m_2}^{2*})$ be a multi-index and a Nikishin system associated with  ${\bf n}_2$ and ${\mathcal{N}}(\sigma_1^{2},\ldots,\sigma_{m_2}^{2})$ through Lemma $\ref{lem:4}$. Set $d\sigma_{0}^{2*} = \widehat{s}_{1,j}^{2}(x)d\sigma_0^2 $ and $(s_{0,0}^{2*},s_{0,1}^{2*},\ldots,s_{0,m_2}^{2*}) = {\mathcal{N}}(\sigma_{0}^{2*},\sigma_{1}^{2*},\ldots,\sigma_{m_2}^{2*})$. Then,
\begin{equation} \label{ortox}
\int x^{\nu} {\mathcal{A}}_{{\bf n}}(x) d {s}_{0,k}^{2*}(x) = 0, \qquad \nu = 0,\ldots,n_{2,k}^* -1,\qquad k=0,\ldots, m_2.
\end{equation}
\el

{\bf Proof.} For $k=0$, since $d s_{0,0}^{2*} = d \sigma_{0}^{2*} = \widehat{s}_{1,j}^{2}(x)d\sigma_0^2$, (\ref{orto}) reduces to (\ref{ortox}) when ${\mathcal{L}}_{{\bf n}_2}(x) = x^{\nu} \widehat{s}^2_{1,j}(x), \nu =0,\ldots,n_{2,0}^* -1$, taking into consideration that $n_{2,0}^* = n_{2,j}$.

Let $j +1 \leq k \leq m_2$. On account of Lemma \ref{lem:4} and (\ref{eq:M})
\[ {\mathcal{L}}_{{\bf n}_2}(x) = x^{\nu} \widehat{s}^2_{1,k} (x) = (\frac{|s_{1,k}^2|}{|s_{1,j}^2|}x^{\nu} + (-1)^{j} x^{\nu} \widehat{s}_{1,k}^{2*})\widehat{s}_{1,j}^2.
\]
Consequently, from (\ref{orto}) and the orthogonality for $k=0$, we obtain
\[
0 = \int x^{\nu}\widehat{s}^2_{1,k} (x) {\mathcal{A}}_{{\bf n}}(x) d\sigma_0^2(x) = \frac{|s_{1,k}^2|}{|s_{1,j}^2|}\int x^{\nu} {\mathcal{A}}_{{\bf n}}(x) ds_{0,0}^{2*}(x)+ (-1)^{j} \int x^{\nu} {\mathcal{A}}_{{\bf n}}(x)  ds_{0,k}^{2*}(x) =
\]
\[
(-1)^{j} \int x^{\nu} {\mathcal{A}}_{{\bf n}}(x)  ds_{0,k}^{2*}(x), \qquad \nu =0,\ldots,n_{2,k}^* -1,
\]
since $n_{2,k}^* = n_{2,k} \leq n_{2,j} = n_{2,0}^*, k=j+1,\ldots,m_2.$ Thus, for these values of $k$ the assertion  also holds.

If $k \in \{1,\ldots,j-1\}$ and $\max\{{n_{2,l_{k-1}},n_{2,k}}\} = n_{2,k},$ from Lemma \ref{lem:4} and (\ref{eq:M}), it follows that
\[
{\mathcal{L}}_{{\bf n}_2}(x) = x^{\nu} \widehat{s}^2_{1,k} (x) = (\frac{|s_{1,k}^2|}{|s_{1,j}^2|}x^{\nu} + (-1)^{k-1}x^{\nu} \widehat{s}_{1,k}^{2*})\widehat{s}_{1,j}^2.
\]
Using the same arguments as in the previous case, we obtain what is needed. Similarly, when $k =j$ and $\min\{n_{2,l_{j-2}},n_{2,j-1}\} = n_{2,j-1},$ then
\[
{\mathcal{L}}_{{\bf n}_2}(x) = x^{\nu} \widehat{s}^2_{1,j-1} (x) = (\frac{|s_{1,j-1}^2|}{|s_{1,j}^2|}x^{\nu} + (-1)^{j-1}x^{\nu} \widehat{s}_{1,j}^{2*})\widehat{s}_{1,j}^2
\]
and integrating the statement follows.

Assume that (\ref{ortox}) holds for all values of the parameter less than $k, 1\leq k \leq j,$ and let us show that it is also true for $k$. When $\max\{{n_{2,l_{k-1}},n_{2,k}}\} = n_{2,k}, 1 \leq k \leq j-1,$ or $\min\{n_{2,l_{j-2}},n_{2,j-1}\} = n_{2,j-1}, k=j,$ we just proved that (\ref{eq:M}) is satisfied, so we must only consider the cases when $\max\{n_{2,l_{k-1}},n_{2,k}\} = n_{2,l_{k-1}}, 1 \leq k \leq j-1,$ and $\min\{n_{2,l_{j-2}},n_{2,j-1}\} = n_{2,l_{j-2}}, k=j.$ In this situation, if $1\leq k \leq j-1$
we have that
\[
{\mathcal{L}}_{{\bf n}_2}(x) = x^{\nu} \widehat{s}^2_{1,l_{k-1}} (x) = (C_0 x^{\nu}  + \sum_{i=1}^{k-1} C_1 x^{\nu}\widehat{s}_{1,i}^{2*}  + (-1)^{l_{k-1}} x^{\nu}\widehat{s}_{1,k}^{2*})\widehat{s}_{1,j}^2, \qquad (\widehat{s}^2_{1,0} \equiv 0).
\]
where $C_i, i=1,\ldots,k-1,$ are constants, $C_0$ is also a constant if $l_{k-1} \neq 0$ and it is a first degree polynomial when $l_{k-1} = 0$. From (\ref{orto}) and the induction hypothesis, it follows that
\[
0 = \int x^{\nu}\widehat{s}^2_{1,l_{k-1}} (x) {\mathcal{A}}_{{\bf n}}(x) d\sigma_0^2(x) = \sum_{i=0}^{k-1}  \int C_i x^{\nu} {\mathcal{A}}_{{\bf n}}(x) ds_{0,i}^{2*}(x)  +
\]
\[ (-1)^{l_{k-1}} \int x^{\nu} {\mathcal{A}}_{{\bf n}}(x) ds_{0,k}^{2*}(x) = (-1)^{l_{k-1}} \int x^{\nu} {\mathcal{A}}_{{\bf n}}(x) ds_{0,k}^{2*}(x), \qquad \nu = 0,\ldots,n_{2,k}^* -1,
\]
since $n_{2,k}^* = n_{2,l_{k-1}} = \min\{n_{2,0},\ldots,n_{2,k-1}\} \leq n_{2,i}^*, i=0,\ldots,k-1$, and when $l_{k-1} =0$ then $n_{2,0} < n_{2,j} = n_{2,0}^*$. Notice that we have already proved (\ref{ortox}) for all values of the parameter up to $j-1$. When $k=j$ and $\min\{n_{2,l_{j-2}},n_{2,j-1}\} = n_{2,l_{j-2}}$, we proceed analogously taking
${\mathcal{L}}_{{\bf n}_2} = x^{\nu}\widehat{s}_{1,l_{j-2}}^2$.  Again,   $n_{2,j}^* \leq n_{2,i}^*, i=0,\ldots,n_{2,j-1}^*$ and we can complete the induction. \hfill $\Box$ \vspace{0,2cm}

\begin{rem}
\label{rem:2}
Fix $k \in \{1,\ldots,m_2\}$. Taking $p_k \equiv 1, p_i \equiv 0, i \neq k, i=0,\ldots,m_2$, one obtains the formula that links $s_{0,k}^2$ with the measures $s_{0,0}^{2*},\ldots,s_{0,k}^{2*}$ (not all have to appear).
\end{rem}

{\bf Proof of Corollary \ref{cor:2}.} From Lemma \ref{lem:orto}, we have that
\begin{equation}
\label{eq:O}
 \int x^{\nu} Q_{\bf n}(x) ds_{0,k}(x) =0, \qquad \nu=0,\ldots,n_k -1,\qquad k=0,\ldots,m.
\end{equation}
Using the definition of type II Pad\'{e} approximation, we have that for any polynomial $Q, \deg Q \leq n_k$
\[ Q(z)(Q_{\bf n}(z) \widehat{s}_{0,k}(z) -  P_{{\bf n},k}(z)) = {\mathcal{O}}(1/z), \qquad z \to \infty.
\]
On account of (\ref{eq:g}), it follows that
\begin{equation}
 \label{eq:N}
Q(z)(Q_{\bf n}(z) \widehat{s}_{0,k}(z) -  P_{{\bf n},k}(z)) = \int \frac{Q(x)Q_{\bf n}(x)ds_{0,k}(x)}{z-x}
\end{equation}
Taking $Q \equiv 1$, we obtain that
\[ P_{{\bf n},k}(z) = \int \frac{Q_{\bf n}(z) - Q_{\bf n}(x)}{z-x} ds_{0,k}(x),\qquad \frac{P_{{\bf n},k}(z)}{Q_{\bf n}(z)} = \frac{1}{Q_{\bf n}(z)} \int \frac{Q_{\bf n}(z) - Q_{\bf n}(x)}{z-x} ds_{0,k}(x).
\]
Consequently, since the zeros of $Q_{\bf n}(z) = \prod_{i=1}^{|{\bf n}|} (z- x_{{\bf n},i})$ are simple and lie in the interior of $\mbox{Co}(\supp \sigma_{0})$
\begin{equation}
\label{eq:P}
\frac{P_{{\bf n},k}(z)}{Q_{\bf n}(z)} = \sum_{i=1}^{|{\bf n}|} \frac{\lambda_{{\bf n},k,i}}{z- x_{{\bf n},i}}, \,\,\,\, \lambda_{{\bf n},k,i} = \lim_{z \to x_{{\bf n},i}} (z- x_{{\bf n},i})\frac{P_{{\bf n},k}(z)}{Q_{\bf n}(z)} =  \int \frac{ Q_{\bf n}(x)}{Q_{\bf n}'(x_{{\bf n},i})} \frac{ds_{0,k}(x)}{x- x_{{\bf n},i}}.
\end{equation}

Let $p$ be an arbitrary polynomials of degree $\leq |{\bf n}| + n_k -1$ and $\ell_{\bf n}(z) = \sum_{i=0}^{|{\bf n}|}\frac{Q_{\bf n}(z)p(x_{{\bf n},i})}{Q_{\bf n}'(x_{{\bf n},i})(z- x_{{\bf n},i})}$ be the Lagrange polynomial of degree $|{\bf n}|-1$ which interpolates $p$ at the zeros of $Q_{\bf n}$. Then
\[ (p - \ell_{\bf n})(z) = q(z)Q_{\bf n}(z), \qquad \deg q \leq n_k -1.
\]
From (\ref{eq:O}),
\[ \int (p - \ell_{\bf n})(x) d s_{0,k}(x) = 0.
\]
Consequently, using (\ref{eq:P})
\[ \int p(x) d s_{0,k}(x) = \int   \ell_{\bf n}(x) d s_{0,k}(x) = \sum_{i=1}^{|{\bf n}|} p(x_{{\bf n},i})\int \frac{Q_{\bf n}(x)d s_{0,k}(x)}{Q_{\bf n}(x_{{\bf n},i})(x - x_{{\bf n},i})} =  \sum_{i=1}^{|{\bf n}|} \lambda_{{\bf n},k,i}p(x_{{\bf n},i}),
\]
which gives the first statement of the corollary.

Let us consider the case when $n_0 = \max\{n_0,n_1-1,\ldots,n_m-1\}$. According to (\ref{eq:N}) with $Q \equiv 1$, Fubini's theorem, and (\ref{eq:O})
\[ \int t^{\nu}(Q_{\bf n}(t) \widehat{s}_{0,0}(t) -  P_{{\bf n},0}(t))ds_{1,k}(t) = \int t^{\nu} \int \frac{Q_{\bf n}(x)ds_{0,0}(x)}{t-x} ds_{1,k}(t) =
\]
\[ \int Q_{\bf n}(x)\int \frac{t^{\nu} \mp x^{\nu}}{t-x} ds_{1,k}(t)ds_{0,0}(x) = \int q_{\nu}(x) Q_{\bf n}(x)d s_{0,0}(x) - \int x^{\nu} Q_{\bf n}(x)d s_{0,k}(x) =0,
\]
for all $\nu = 1,\ldots,n_k-1$ and $k=1,\ldots,m$, since $q_{\nu}(x) = \int \frac{t^{\nu} - x^{\nu}}{t-x} ds_{1,k}(t)$ is a polynomial such that $\deg q_{\nu} \leq n_k -2 \leq n_0 -1$. This implies that
\[ \int  (p_1(t) + \sum_{k=2}^{m}p_k(t)\widehat{s}_{2,k}(t))(Q_{\bf n}(t) \widehat{s}_{0,0}(t) -  P_{{\bf n},0}(t))ds_{1,1}(t) =0,
\]
for arbitrary polynomials $p_k, \deg p_k \leq n_k -1, k=1,\ldots,m$. Since, by Theorem \ref{teo:1}, $(1,\widehat{s}_{2,2},\ldots,\widehat{s}_{2,m})$ is an AT-system, it follows that $Q_{\bf n} \widehat{s}_{0,0} -  P_{{\bf n},0}$, has at least $|{\bf n}| - n_0$ sign changes in the interior of $\mbox{Co}(\supp \sigma_{1})$.

Let $Q_{{\bf n},1}$ be a monic polynomial of degree $|{\bf n}| - n_0$ whose simple zeros are points where $Q_{\bf n} \widehat{s}_{0,0} -  P_{{\bf n},0}$ changes sign on $\mbox{Co}(\supp \sigma_{1})$. It is easy to verify that
\[ \frac{z^{\nu}(Q_{\bf n} \widehat{s}_{0,0} -  P_{{\bf n},0})(z)}{Q_{{\bf n},1}(z)} = {\mathcal{O}}(1/z^2), \qquad z \to \infty, \qquad \nu =0,\ldots,|{\bf n}| -1.
\]
Using (\ref{eq:f}), we obtain
\[ \int x^{\nu} Q_{\bf n}(x) \frac{d s_{0,0}(x)}{Q_{{\bf n},1}(x)}
=0 , \qquad \nu =0,\ldots,|{\bf n}| -1. \]
Then, for any $\widetilde{Q}, \deg \widetilde{Q} \leq |{\bf n}|$,
\[ \int Q_{\bf n}(x) \frac{\widetilde{Q}(z) - \widetilde{Q}(x)}{z-x} \frac{ds_{0,0}(x)}{Q_{{\bf n},1}(x)} = 0
\]
which implies that
\[ \widetilde{Q}(z) \int \frac{Q_{\bf n}(x)}{z-x} \frac{ds_{0,0}(x)}{Q_{{\bf n},1}(x)} =
\int  \frac{\widetilde{Q}(z)Q_{\bf n}(x)}{z-x} \frac{ds_{0,0}(x)}{Q_{{\bf n},1}(x)}.
\]
In particular, with $z = x_{{\bf n},i}$, taking $\widetilde{Q}(x) = Q_{{\bf n},1}(x)$ and then $\widetilde{Q}(x) = \frac{Q_{\bf n}(x)}{Q_{\bf n}'(x_{{\bf n},i})(x- x_{{\bf n},i})}$, we have
\[ \lambda_{{\bf n},0,i} =   \int \frac{ Q_{\bf n}(x)}{Q_{\bf n}'(x_{{\bf n},i})} \frac{ds_{0,0}(x)}{x- x_{{\bf n},i}} = \int \frac{Q_{{\bf n},1}(x) Q_{\bf n}(x)}{Q_{\bf n}'(x_{{\bf n},i})(x- x_{{\bf n},i})} \frac{ds_{0,0}(x)}{Q_{{\bf n},1}(x) } =
\]
\[\frac{Q_{\bf n}'(x_{{\bf n},i})}{Q_{\bf n}'(x_{{\bf n},i})}Q_{{\bf n},1}(x_{{\bf n},i} )\int \frac{ Q_{\bf n}(x)}{Q_{\bf n}'(x_{{\bf n},i})(x- x_{{\bf n},i})} \frac{ds_{0,0}(x)}{Q_{{\bf n},1}(x) } =
\]
\[\int \left(\frac{Q_{\bf n}(x)}{Q_{\bf n}'(x_{{\bf n},i})(x- x_{{\bf n},i})} \right)^2\frac{Q_{{\bf n},1}(x_{{\bf n},i})}{Q_{{\bf n},1}(x)} ds_{0,0}(x), \qquad i =1,\ldots,|{\bf n}|.
\]
Since $\left(\frac{Q_{\bf n}(x)}{Q_{\bf n}'(x_{{\bf n},i})(x- x_{{\bf n},i})} \right)^2\frac{Q_{{\bf n},1}(x_{{\bf n},i})}{Q_{{\bf n},1}(x)}$ is positive for all $x \in \mbox{Co}(\supp s_{0,0})$, the second statement of Corollary \ref{cor:2} follows for $k=0$. Using standard arguments of one-sided polynomial approximation of Riemann-Stieltjes integrable functions (see e.g. \cite[Theorem 15.2.2]{Sze} and \cite[Lemma 2]{LIF}), the third statement is a consequence of the first two for any sequence of multi-indices $\Lambda \subset {\mathbb{Z}}_+^{m+1}$ such that for all ${\bf n} \in \Lambda, n_0 = \max\{n_0,n_1-1,\ldots,n_m -1\}$. In particular, the second and third statements are valid when ${\bf n} = (n,n+1,\ldots,n+1)$.

In the rest of the proof, we restrict our attention to multi-indices of the form ${\bf n} = (n,n+1,\ldots,n+1)$. Fix $k \in \{1,\ldots,m\}.$ Since we have that $n_k = n+1 = \max\{n_0 +1,n_1,\ldots,n_m\},$ we can apply Lemma \ref{lem:ortox} with $j=k$ and we obtain that $Q_{\bf n}$ is multiple orthogonal with respect to ${\bf n}^*$, which has $n+1$ in the first component, and a Nikishin system ${\mathcal{N}}(\sigma_0^*,\ldots,\sigma_m^*)$  whose first measure is $s_{0,0}^*= s_{0,k}$. Consequently, the coefficients
\[ \lambda_{{\bf n},k,i} =  \int \frac{ Q_{\bf n}(x)}{Q_{\bf n}'(x_{{\bf n},i})} \frac{ds_{0,k}(x)}{x- x_{{\bf n},i}} = \int \frac{ Q_{\bf n}(x)}{Q_{\bf n}'(x_{{\bf n},i})} \frac{ds_{0,0}^*(x)}{x- x_{{\bf n},i}}
\]
must all have the same sign as $s_{0,k}$. The convergence of the quadratures is obtained as before. \hfill $\Box$

\section{Proof of Theorems \ref{teo:3}-\ref{teo:4}}

{\bf Proof of Theorem \ref{teo:3}.} For $m=0$ the result is trivially true since ${\mathcal{L}}_{\bf n} = p_0$. When $m=1$ it is easy to deduce. Indeed, if $n_0 \geq n_1$, take $\lambda$ the identity and $S(\lambda) = {\mathcal{N}}(\sigma_{1})$; otherwise, $n_0 < n_1$ and by Lemma \ref{lem:4}
\[ p_0 + p_1\widehat{s}_{1,1} = (p_0^* + p_1^* \widehat{\tau}_{1,1}) \widehat{s}_{1,1}, \quad \deg p_0^* \leq n_1-1, \quad \deg p_1^* \leq n_0-1.
\]
Hence, the solution is $\lambda$ such that $\lambda(0) =1, \lambda(1) =0,$ and $S(\lambda) = {\mathcal{N}}(\tau_{1,1})$. In the following $m \geq 2$.

Next, let us consider the case when $n_0 = \max\{n_0,\ldots,n_m\}$. If $n_0 \geq \cdots \geq n_m$, the result is trivial taking $\lambda$ the identity and $S(\lambda) = {\mathcal{N}}(\sigma_1,\ldots,\sigma_m)$. Otherwise,   there exists $\overline{m}, 0 \leq \overline{m} \leq m-2,$ such that $n_0 \geq \cdots \geq n_{\overline{m}}$, $n_{\overline{m}} = \max\{n_{\overline{m}},\ldots,n_m\}$, and  $n_{\overline{m} +1} < \max\{n_{\overline{m}+2},\ldots,n_m\}$. (Consequently, $n_{\overline{m}+1} < n_{\overline{m}}$.)

We have
\[ {\mathcal{L}}_{\bf n} =  p_0 + \sum_{k=1}^m p_k \widehat{s}_{1,k} = p_0 + \sum_{k=0}^{\overline{m}} p_k \widehat{s}_{1,k}
+ \sum_{k=\overline{m}+1}^m p_k \la \sigma_1,\ldots,\sigma_{\overline{m}},s_{\overline{m}+1,k}\widehat{\ra}.
\]
It is easy to check (see, for example, Lemma 2.1 in \cite{LL}) that for each $k \in \{ \overline{m}+1,\ldots,m\},$
\begin{equation}
\label{eq:Q}
p_k \la \sigma_1,\ldots,\sigma_{\overline{m}},s_{\overline{m}+1,k}\widehat{\ra} = \ell_{k,0} +
\sum_{i=1}^{\overline{m}} \ell_{k,i} \widehat{s}_{1,i} + \la \sigma_1,\ldots,\sigma_{\overline{m}},p_k s_{\overline{m}+1,k} \widehat{\ra},
\end{equation}
$\deg \ell_{k,i} \leq \deg p_k -1,$ and these polynomials have real coefficients.  Since $n_k \leq n_{\overline{m}} \leq n_{\overline{m}-1} \leq n_0 $ whenever $k \in \{\overline{m}+1,\ldots,n_m\}$, the polynomials  $\ell_{k,i}, i=0,\ldots,\overline{m},$ $ k=\overline{m}+1,\ldots,m$, are absorbed by the polynomials $p_k, k=0\ldots,\overline{m},$ without altering the bound on the degrees of the second. Therefore, there exist polynomials with real coefficients $\widetilde{p}_k,  \deg \widetilde{p}_k \leq n_k -1,k=0\ldots,\overline{m},$ such that
\[ {p}_0 + \sum_{k=1}^m  {p}_k \widehat{s}_{1,k} = \widetilde{p}_0 + \sum_{k=0}^{\overline{m}} \widetilde{p}_k \widehat{s}_{1,k} + \la \sigma_1,\ldots,\sigma_{\overline{m}},\sum_{k=\overline{m}+1}^m p_k s_{\overline{m}+1,k}\widehat{\ra} =
\]
\begin{equation}
\label{eq:R}
\widetilde{p}_0 + \sum_{k=0}^{\overline{m}} \widetilde{p}_k \widehat{s}_{1,k} + \la \sigma_1,\ldots,\sigma_{\overline{m}},(p_{\overline{m}+1} + \sum_{k=\overline{m}+2}^m p_k \widehat{s}_{\overline{m}+2,k})  \sigma_{\overline{m}+1}\widehat{\ra}.
\end{equation}

By assumption $n_{\overline{m}  +1} < \max\{n_{\overline{m}+2},\ldots,n_m\}.$ So, we can apply Lemma \ref{lem:4} on the linear form $p_{\overline{m}+1} + \sum_{k=\overline{m}+2}^m p_k \widehat{s}_{\overline{m}+2,k}$. Thus, there exist a Nikishin system ${\mathcal{N}}(\sigma_{\overline{m}+2}^*,\ldots,\sigma_{m}^*)$, a multi-index $(n_{\overline{m}+1}^*,\ldots,n_m^*) \in {\mathbb{Z}}_+^{m - \overline{m}}$, which is a permutation of $(n_{\overline{m}+1},\ldots,n_m)$, and polynomials with real coefficients $p_k^*, \deg p_k^* \leq n_k^* -1, k= \overline{m} +1,\ldots,m,$ such that
\[
p_{\overline{m}+1} + \sum_{k=\overline{m}+2}^m p_k \widehat{s}_{\overline{m}+2,k} = ( {p}_{\overline{m}+1}^* + \sum_{k=\overline{m}+2}^m  {p}_k^* \widehat{s}_{\overline{m}+2,k}^*)\widehat{s}_{\overline{m}+2,j},
\]
where $j$ is such that $n_j = \max\{n_{\overline{m}+1},\ldots,n_m\}$ and $n_{\overline{m}+1}^* = n_j$. Substitute this formula in (\ref{eq:R}) and reverse the application of (\ref{eq:Q}) to pull out the polynomials $\widetilde{p}_k$ of the product of measures. The new polynomials $l_{k,i}, k \in \{\overline{m}+1,\ldots,n_m\}$ which arise from this second application of (\ref{eq:Q}) are also absorbed by the polynomials $\widetilde{p}_k, k=0,\ldots,\overline{m}$ in (\ref{eq:R}) without changing the bound on their degrees.  Therefore, we get that for certain polynomials with real coefficients $p_k^*, k=0,\ldots,m$
\begin{equation}
 \label{eq:S}
 {\mathcal{L}_{\bf n}} = p_0^* + \sum_{k=1}^{\overline{m}}p_k^* \widehat{s}_{1,k} + p_{\overline{m}+1}^* \widehat{s}_{1,j} + \sum_{k=\overline{m}+2}^m p_k^* \la \sigma_1,\ldots,\sigma_{\overline{m}},s_{\overline{m}+1,j},\sigma_{\overline{m}+2}^*,\ldots,\sigma_{m}^*\widehat{\ra},
\end{equation}
which is a linear form generated by the multi-index ${\bf n}^* = (n_0,\ldots,n_{\overline{m}},n_j,n_{\overline{m}+2}^*,\ldots,n_m^*),$ and the Nikishin system ${\mathcal{N}}(\sigma_1,\ldots,\sigma_{\overline{m}},s_{\overline{m}+1,j},\sigma_{\overline{m}+2}^*,\ldots,\sigma_{m}^*)$.

Now, we have that $n_0 \geq \cdots \geq n_{\overline{m}} \geq n_j$ and $n_j = \max \{n_j,n_{\overline{m}+2}^*,\ldots,n_m^*\}$. If $n_{\overline{m}+2} \geq \cdots \geq n_m,$ we are done taking $\lambda$ such that $n_{\lambda(k)} = n_k^*$, $\widehat{s}_{1,\lambda(0)} \equiv 1,$ and
\[S(\lambda) = {\mathcal{N}}(\sigma_1,\ldots,\sigma_{\overline{m}},s_{\overline{m}+1,j},\sigma_{\overline{m}+2}^*,\ldots,\sigma_{m}^*).
\]
Otherwise, we repeat the process with the linear form on the right hand of (\ref{eq:S}). The new $\overline{m}$ will certainly be larger that the previous one and in a finite number of iterations we reorganize the entries of ${\bf n}$ in decreasing order obtaining with it $\lambda$ and $S(\lambda)$.

If $n_0 < \max\{n_0,\ldots,n_m\}$, we apply first Lemma \ref{lem:4} and  then proceed as have done before but with the form ${\mathcal{L}}_{\bf n}^*= p_0^* + \sum_{k=1}^m p_k^*\widehat{s}_{1,k}^*$ and the multi-index ${\bf n}^*$ arising from that lemma. Thus, we find a permutation $\overline{\lambda}$ of $(0,\ldots,m)$ such that $n_{\overline{\lambda}(0)}^*\geq \cdots \geq n_{\overline{\lambda}(m)}^*$, $S(\overline{\lambda}) = {\mathcal{N}}(\rho_1,\ldots,\rho_m),$ and polynomials with real coefficients $q_k, k=0,\ldots,m,$ such that
\[ {\mathcal{L}}_{\bf n}^* = q_0 + \sum_{k=1}^m q_k \widehat{r}_{1,k}, \qquad \deg q_k \leq n_{\overline{\lambda}(k)}^* -1, \qquad k=0,\ldots,m.
\]
If $\lambda^*$ is the permutation due to Lemma \ref{lem:4} which transports ${\bf n}$ into ${\bf n}^*$,
taking $\lambda = \overline{\lambda} {\rm o} \lambda^*$ and $S(\lambda) = S(\overline{\lambda})$, applying the formula of Lemma \ref{lem:4} the assertion of Theorem \ref{teo:3} again follows. \hfill $\Box$

\vspace{0,2cm}

Using induction we could have reduced a bit the proof of Theorem \ref{teo:3}. Nevertheless, for the proof of Theorem \ref{teo:4} it was convenient to underline the fact that Theorem \ref{teo:3} is a consequence of iterating the use of Lemma \ref{lem:4} and the trick exhibited in (\ref{eq:Q})-(\ref{eq:S}). \vspace{0,2cm}

{\bf Proof of Theorem \ref{teo:4}.} If $m_2 =0$ or $m_2 =1$ and $n_{2,0} \geq n_{2,1}$ the result is trivial. For $m_2 = 1$ and $n_{2,0} < n_{2,1}$ the statement is contained in Lemma \ref{lem:ortox}. So, we restrict our attention to $m_2 \geq 2$.

First, we consider the case when $n_{2,0} = \max\{n_{2,0},\ldots,n_{2,m_2}\}$.
In this situation, the result is trivial again if $n_{2,0} \geq \cdots \geq n_{2,m_2}$. If this is not the case, there exists $\overline{m}, 0 \leq \overline{m} \leq m_2 -2,$ such that
$n_{2,0} \geq \cdots \geq n_{2,\overline{m}}$, $n_{2,\overline{m}} = \max\{n_{2,\overline{m}},\ldots,n_{2,m_2}\}$, and  $n_{2,\overline{m} +1} < \max\{n_{2,\overline{m}+2},\ldots,n_{2,m_2}\}$.

According to (\ref{eq:S}), for any polynomials $p_k, \deg p_k \leq n_{2,k}-1,$ (\ref{eq:S}) takes place; that is,
\[
 p_0 + \sum_{k=1}^{m_2}p_k \widehat{s}_{1,k}^2 = p_0^* + \sum_{k=1}^{\overline{m}}p_k^* \widehat{s}_{1,k}^2 + p_{\overline{m}+1}^* \widehat{s}_{1,j}^2 + \sum_{k=\overline{m}+2}^m p_k^* \la \sigma_1^2,\ldots,\sigma_{\overline{m}}^2,s_{\overline{m}+1,j}^2,\sigma_{\overline{m}+2}^{2*},\ldots,\sigma_{m}^{2*}\widehat{\ra},
\]
where $j$ is such that $n_{2,j} = \max\{n_{2,\overline{m}+1},\ldots,n_{2,m_2}\}.$ This linear form is generated
\[
{\bf n}_2^* = (n_{2,0}^*,\ldots,n_{2,m_2}^*) = (n_{2,0},\ldots,n_{2,\overline{m}},n_{2,j},n_{2,\overline{m}+2}^*,\ldots,n_{2,m_2}^*),
 \]
 and the Nikishin system \[{\mathcal{N}}(\sigma_1^{2*},\ldots,\sigma_{m_2}^{2*}) = {\mathcal{N}}(\sigma_1^2,\ldots,\sigma_{\overline{m}}^2,s_{\overline{m}+1,j}^2,\sigma_{\overline{m}+2}^{2*},\ldots
,\sigma_{m_2}^{2*}).\]
Consider the extended Nikishin system ${\mathcal{N}}(\sigma_0^{2*},\sigma_1^{2*},\ldots,\sigma_{m_2}^{2*})= {\mathcal{N}}(\sigma_0^{2},\sigma_1^{2*},\ldots,\sigma_{m_2}^{2*}).$ The form ${\mathcal{A}}_{\bf n}$ is of multiple orthogonality with respect to ${\bf n}_2^*$ and the extended Nikishin system.

In fact, by definition,  ${\mathcal{A}}_{\bf n}$ satisfies
\[ \int x^{\nu}{\mathcal{A}}_{\bf n}(x) ds_{0,k}^{2*}(x) =0 , \qquad \nu =0,\ldots,n_{2,k}^* -1, \qquad k=0,\ldots,\overline{m} +1,
\]
since $s_{0,k}^{2*} = s_{0,k}^{2}, n_{2,k}^* = n_{2,k}, k=0,\ldots,\overline{m}, s_{0,\overline{m}+1}^{2*} = s_{0,j}^{2}$, and $n_{2,\overline{m}+1}^* = n_{2,j}$. To prove
\[ \int x^{\nu}{\mathcal{A}}_{\bf n}(x) ds_{0,k}^{2*}(x) =0 , \qquad \nu =0,\ldots,n_{2,k}^* -1, \qquad k=\overline{m} +2,\ldots,m_2,
\]
one follows  arguments similar to those employed in proving Lemma \ref{lem:ortox}, choosing particular expressions for ${\mathcal{L}}_{{\bf n}_2}$ of the form $x^{\nu}\widehat{s}^2_{1,\overline{k}}$ ($\overline{k}$ is not always equal to $k$), and taking into consideration  (\ref{eq:Q})-(\ref{eq:S})  as well as (\ref{eq:M}). The details are left to the reader.

Once we have proved that ${\mathcal{A}}_{\bf n}$ is of multiple orthogonality with respect to ${\bf n}_2^*$ and the extended Nikishin system, one repeats the process finding a new $\overline{m}$, which is obviously larger than the previous one, and in a finite number of iterations the statement follows.

If $n_{2,0} < \max\{n_{2,0},\ldots,n_{2,m_2}\}$, the proof is reduced to the previous case by Lemma \ref{lem:ortox}. \hfill $\Box$

\section{Proof of Theorems \ref{teo:5}-\ref{teo:6} and Corollary \ref{cor:3}}

If we apply Theorem \ref{teo:3} to the form ${\mathcal{A}}_{\bf n}$, we obtain that there exists a permutation $\lambda_1$ of $(0,\ldots,m_1)$ and an associated Nikishin system $S(\lambda_1) = (r_{1,1}^1,\ldots,r_{1,m_1}^1) = {\mathcal{N}}(\rho_{1}^1,\ldots,\rho_{m_1}^1)$ such that
\[ {\mathcal{A}}_{\bf n} = a_{{\bf n},0} + \sum_{k=1}^{m_1} a_{{\bf n},k} \widehat{s}_{1,k}^1 = (b_{{\bf n},0} + \sum_{k=1}^{m_1} b_{{\bf n},k} \widehat{r}_{1,k})\widehat{s}_{1,\lambda_1(0)} = {\mathcal{B}}_{\bf n}\widehat{s}_{1,\lambda_1(0)},
\]
where $\widehat{s}_{1,\lambda_1(0)} \equiv 1$ if $\lambda_1(0) =0$, and $\deg b_{{\bf n},k} \leq n_{1,\lambda_1(k)} -1, k=0,\ldots,m_1.$ On the other hand, from Theorem  \ref{teo:4}, we know that there exists a permutation $\lambda_2$   of $(0,\ldots,m_2)$ and a  Nikishin system ${\mathcal{N}}(\rho_0^2,\ldots,\rho_{m_2}^2)$, where  $\rho^2_0 = \widehat{s}^2_{1,\lambda_2 (0)}\sigma_0^2$, such that for each $k=0,\ldots,m_2,$
\[ \int x^{\nu} {\mathcal{A}}_{{\bf n}}(x) d r_{0,k}^2(x) = \int x^{\nu} {\mathcal{B}}_{{\bf n}}(x) \widehat{s}_{1,\lambda_1(0)}(x) d r_{0,k}^2(x) =0,\qquad \nu =0,\ldots,n_{2,\lambda_2 (k)} -1.
\]
Therefore, ${\mathcal{B}}_{{\bf n}}$ is a linear form, generated by the multi-index $(n_{1,\lambda_1(0)},\ldots,n_{1,\lambda_1(m_1)})$ and $S(\lambda_1)$, which is of multiple orthogonality with respect to the multi-index  $(n_{2,\lambda_2(0)},\ldots,n_{2,\lambda_2(m_2)})$ and the Nikishin system $(\widehat{s}_{1,\lambda_1(0)}r_0^2,r_1^2,\ldots,r_{m_2}^2)={\mathcal{N}}(\widehat{s}_{1,\lambda_1(0)}
\rho_0^2,\rho_1^2,\ldots,\rho_{m_2}^2)$. In other words,
\[ {\mathbb{B}}_{\bf n} = (b_{{\bf n},0},\ldots,b_{{\bf n},m_1})
\]
is the mixed type multiple orthogonal polynomial relative to the pair of Nikishin systems $(\widetilde{S}^1,\widetilde{S}^2)$ and the multi-index $\widetilde{\bf n} = (\widetilde{\bf n}_1;\widetilde{\bf n}_2) \in {\mathbb{Z}}_+^{m_1+1} \times {\mathbb{Z}}_+^{m_2 +1}$, where
\[ (\widetilde{S}^1,\widetilde{S}^2) = ({\mathcal{N}}(\widehat{s}_{1,\lambda_1(0)}
\rho_0^2,\rho_{1}^1,\ldots,\rho_{m_1}^1),{\mathcal{N}}(\widehat{s}_{1,\lambda_1(0)}
\rho_0^2,\rho_1^2,\ldots,\rho_{m_2}^2)),
\]
and
\[ \widetilde{\bf n}_i = (n_{i,\lambda_i(0)},\ldots,n_{i,\lambda_i(m_1)}), \qquad i=1,2.
\]
Both $\widetilde{\bf n}_1$ and $\widetilde{\bf n}_2$ have decreasing components. Therefore, to derive Theorems \ref{teo:5} and \ref{teo:6} we can apply the results of \cite{FLLS}.

\begin{lemma}
\label{lem:cond}
If $(S^1,S^2)$ satisfies the hypotheses of Theorem 5 $($respectively 6$)$ the same is true for $(\widetilde{S}^1,\widetilde{S}^2)$.
\end{lemma}

{\bf Proof.} The systems  $ \widetilde{S}^1$ and $\widetilde{S}^2$ are obtained transforming the generating measures of ${S}^1$ and ${S}^2$ through inversion of measures and multiplication by Cauchy transforms of measures supported on disjoint intervals. We have to check that these operations preserve the quasi-regularity of supports, the regularity of measures, and the property concerning the Radon-Nikodym derivative of the measure.

Let $\sigma \in {\mathcal{M}}(\Delta), \Delta = \mbox{\rm Co}(\supp
\sigma)$, $\tau$ is the inverse measure of $\sigma$, and $g$ is a continuous function on $\Delta$ with constant sign and different from zero on $\Delta$.

It is trivial that the supports of $\sigma$ and $g\sigma$ coincide and that $\sigma' > 0$ if and only if $g\sigma' >0$. It is well known and easy to verify (using, for example, the minimality property of monic orthogonal polynomials) that $\sigma \in \mbox{\bf Reg}$ if and only if $g \sigma \in \mbox{\bf Reg}$ as well.

Regarding the inversion of measures, the Stieltjes-Plemelj inversion formula implies that the continuous parts of the supports of $\sigma$ and $\tau$ coincide. From the formula relating  $\widehat{\sigma}$ and $\widehat{\tau}$ it is obvious that isolated mass points of $\sigma$ outside its continuous support become zeros of $\widehat{\tau}$ (thus are no longer in the support of $\tau$). On the other hand, in each connected component of $\Delta \setminus \supp \sigma$, $\widehat{\sigma}$ may have at most one
zero (counting multiplicity), because $\widehat{\sigma}$ is strictly monotonic
when restricted to any one of those components. Such zeros of
$\widehat{\sigma}$ become mass points of $\tau$. They are isolated, so they can only accumulate on $\supp \sigma$.  Therefore, if $\supp \sigma = \widetilde{E} \cup e$, where $\widetilde{E}$ is regular with respect to the Dirichlet problem and $e$ is at most a denumerable set of points which may only accumulate on $\widetilde{E}$, then $\supp \tau = \widetilde{E} \cup \widetilde{e}$ where  $\widetilde{e}$ is at most a denumerable set of points which may only accumulate on $\widetilde{E}$. In particular, the same holds when $\widetilde{E}$ is an interval (case of Theorem \ref{teo:6}).

If $\sigma \in \mbox{\bf Reg}$ then $\tau \in \mbox{\bf Reg}$. Indeed, the denominator $Q_n$ of the $n$-th diagonal Pad\'{e} approximant of $\widehat{\sigma}$, taken with leading coefficient equal to $1$, is the $n$-th monic orthogonal polynomial with respect to $\sigma$. The numerator $P_{n-1}$ is an $(n-1)$-th orthogonal polynomial with respect to $\tau$.  By Markov's theorem
\[ \lim_{n} \frac{P_{n-1}(z)}{Q_n(z)} = \widehat{\sigma}(z), \qquad {\mathcal{K}} \subset \overline{\mathbb{C}} \setminus \Delta,
\]
In particular, the leading coefficient $c_{n-1}$ of $P_{n-1}$ satisfies
\[ \lim_{n} c_{n-1} = \lim_n \lim_{z \to \infty} \frac{zP_{n-1}(z)}{Q_n(z)} = \lim_{z\to \infty} z \widehat{\sigma}(z)  = \sigma (\Delta) \neq 0.
\]
Therefore,
\[ \lim_n |Q_n(z)|^{1/n} = \mbox{cap}(\supp \sigma) e^{g_{\Omega}(z;\infty)} \,\,\Leftrightarrow\,\, \lim_n \left|\frac{P_n(z)}{c_n}\right|^{1/n} = \mbox{cap}(\supp \sigma) e^{ g_{\Omega}(z;\infty)},
\]
uniformly on compact subsets of ${\mathbb{C}} \setminus \Delta$, where $g_{\Omega}(z;\infty)$ denotes Green's function of the region $\Omega = \overline{\mathbb{C}} \setminus \supp \sigma$ with singularity at $\infty$. But $\supp \sigma$ and $\supp \tau$ differ on a set of capacity zero so their capacities coincide as well as the Green's function of the complement of their supports. The limits above are equivalent to regularity (see \cite[Theorem 3.1.1]{stto}).

That $\sigma' >0$ a.e. on an interval is equivalent to $\tau' >0$ a.e. on the same interval follows from the Stieltjes--Plemelj inversion  formula. \hfill $\Box$ \vspace{0,2cm}

In order to prove Theorem 5 there is still one thing to be considered. The corresponding result  \cite[Theorem 1.3]{FLLS}  for the case of decreasing components in ${\bf n}_1,{\bf n}_2,$ was proved assuming that the supports of the measures were regular. We have to extend its applicability  to the case of quasi-regular supports  because, as follows from the proof of the previous lemma, the regularity of the supports of the measures generating $(S^1,S^2)$ does not guarantee regularity of the supports of the measures which generate  $(\widetilde{S}^1,\widetilde{S}^2)$ since isolated mass points may arise.

We need some notation. ${\mathcal{M}}_1(E)$ denotes the class of
probability measures supported on $E$,   and
\[
V^{\mu}(z) = \int \log \frac{1}{|z - \zeta|} d \mu (z)
\]
the logarithmic potential of the measures $\mu$. If $q_l$ is
a polynomials of degree $l$,
\[ \mu_{q_l} = \frac{1}{l} \sum_{q_l(x) =0} \delta_x
\]
is the associated normalized zero counting measure, where $\delta_x$
is the Dirac measure with mass $1$ at $x$. In \cite[Theorem
I.1.3]{ST} the authors prove \vspace{0,2cm}

\begin{lemma}\label{lemextremal}
Let $E \subset \mathbb{C}$ be a compact subset of the complex plane
and $\phi$ a continuous function on $E$. Then, there exists a unique
$\overline{\mu} \in {\mathcal{M}}_1(E)$ and a constant $w$ such that
\[
V^{\overline{\mu}}(z)+\phi(z) \left\{ \begin{array}{l} \leq w,\quad
z
\in \supp \overline{\mu} \,, \\
\geq w, \quad z \in  E \setminus A,\,\, \mbox{\rm cap}(A) = 0.
\end{array} \right.
\]
\end{lemma}
\vspace{0,2cm}

$\overline{\mu}$ and $w$ are called the equilibrium measure and the
equilibrium constant, respectively, in presence of the
external field $\phi$ on the compact $E$.
\vspace{0,2cm}

We are especially grateful to H. Stahl who gave us the clue for the
following improvement of \cite[Lemma 4.2]{FLLS}.

\begin{lemma} \label{lem:asint}
Let $\sigma \in \mbox{\bf Reg}$, $\supp{\sigma} \subset \mathbb{R}$,
where $\supp{\sigma}$ is quasi-regular. Let $\{\phi_l\}, l \in
\Lambda \subset \mathbb{Z}_+,$ be a sequence of positive continuous
functions on $\supp{\sigma}$ such that
\[
\lim_{l\in \Lambda}\frac{1}{2l}\log\frac{1}{|\phi_l(x)|}= \phi(x)
> -\infty ,
\]
uniformly on $\supp \sigma$. By  $\{q_l\}, l \in \Lambda,$ denote a
sequence of monic polynomials, $\deg q_l = l,$ and
\[
\int x^k q_l(x)\phi_l(x)d\sigma(x)=0,\qquad k=0,\ldots, l-1.
\]
Then
\[
*\lim_{l \in \Lambda}\mu_{q_l} = \overline{\mu},
\]
in the weak star topology of measures, and
\[
\lim_{l\in \Lambda}\left|\int |q_l(x)|^2\phi_l(x)
d\sigma(x)\right|^{1/{2l}}= \exp{(-w)},
\]
where $\overline{\mu}$ and $w$ are the  equilibrium measure and
equilibrium constant in the presence of the external field $\phi$ on
$\supp{\sigma} =: E$. We also have that
\[
\lim_{l \in \Lambda} \left(\frac{|q_l(z)|}{\|q_l
\phi_l^{1/2}\|_E}\right)^{1/l} = \exp{(w - V^{\overline{\mu}}(z))},
\qquad  \mathcal{K} \subset {\mathbb{C}} \setminus \mbox{\rm
Co}(\supp(\sigma)),
 \]
 where $\|\cdot\|_E$ denotes the sup norm on $E$.
\end{lemma}

{\bf Proof.}  Proceeding as in the proof of \cite[Lemma 4.2]{FLLS} one
shows that  for any sequence of monic polynomials  $\{p_l\}, l \in
\Lambda,$  such that $\deg p_l =l$,
\begin{equation} \label{eq:20}
\limsup_{l \in \Lambda} \left(\frac{|p_l(z)|}{\|p_l
\phi_l^{1/2}\|_E}\right)^{1/l} \leq \exp{(w -
V^{\overline{\mu}}(z))}, \qquad  \mathcal{K} \subset {\mathbb{C}},
\end{equation}
and
\begin{equation} \label{eq:21}
\liminf_{l \in \Lambda} \|p_l \phi_l^{1/2}\|_E^{1/l} \geq
\exp{(-w)}.
\end{equation}
In particular, these relations hold for $\{q_l\}, l \in \Lambda$. In \cite[Lemma 4.2]{FLLS},
it is also proved that
\begin{equation} \label{eq:22} \limsup_{l \in \Lambda}
\|q_l\phi_l^{1/2}\|_2^{1/l} \leq \exp({-w}),
\end{equation}
where $\|q_l\phi_l^{1/2}\|_2$ is the $L^2$ norm of $q_l\phi_l^{1/2}$
with respect to $\sigma$. We may assume, without loss of generality, that $\sigma$ is positive. In deducing (\ref{eq:20})-(\ref{eq:22}), the regularity of
$\supp \sigma$ is not required.

Combining (\ref{eq:21})-(\ref{eq:22}), it follows that
\[ \liminf_{l \in \Lambda}
\left(\frac{\|q_l\phi_l^{1/2}\|_E}{\|q_l\phi_l^{1/2}\|_2}\right)^{1/l}
\geq 1.
\]
Should
\begin{equation} \label{eq:23}
\limsup_{l \in \Lambda} \left(\frac{\|q_l
\phi_l^{1/2}\|_E}{\|q_l\phi_l^{1/2}\|_2}\right)^{1/l} \leq 1,
\end{equation}
then
\[
\lim_{l \in \Lambda} \left(\frac{\|q_l
\phi_l^{1/2}\|_E}{\|q_l\phi_l^{1/2}\|_2}\right)^{1/l} = 1.
\]
and due to (\ref{eq:21})-(\ref{eq:22}), we would have
\begin{equation} \label{eq:24}
\limsup_{l \in \Lambda} \|q_l\phi_l^{1/2}\|_E^{1/l} = \limsup_{l \in
\Lambda} \|q_l\phi_l^{1/2}\|_2^{1/l} = \exp({-w}).
\end{equation}

Once (\ref{eq:24}) is attained, with the help of (\ref{eq:20}), one
can conclude the proof as  in \cite[Lemma 7]{FLLS}. So, it remains
to show that (\ref{eq:23}) takes place when we relax the regularity
of $\supp \sigma$ to quasi-regularity.

In  \cite[Theorem 3.2.3]{stto} it is proved (see (v)
$\Rightarrow$ (vi)) that (\ref{eq:23}) holds for any sequence of
polynomials  $\{p_l\}, l \in \Lambda,$  such that $\deg p_l =l$, if
the same property is satisfied when $\phi_l \equiv 1, l \in
\Lambda$. Though the hypothesis of that theorem also contains the
assumption that $\supp \sigma$ be regular, the proof of this
assertion is independent of the regularity condition. Therefore, let
us show that (\ref{eq:23}) holds true when $\phi_l \equiv 1, l \in
\Lambda,$ and $\supp \sigma$ is quasi-regular.

In fact, according to \cite[Theorem 3.2.1 ii)]{stto} , we have that
\begin{equation} \label{eq:25}
\limsup_{l \in \Lambda} \left(\frac{|p_l(z)|}{\|p_l
\|_2}\right)^{1/l} \leq \exp{(g_{\Omega}(z;\infty))}, \qquad
\mathcal{K} \subset {\mathbb{C}},
\end{equation}
where $g_{\Omega}(z;\infty)$ is the Green's function of the region
$\Omega = \overline{\mathbb{C}} \setminus E$ with singularity at
$\infty$. Since $E = \widetilde{E} \cup e$, where $\widetilde{E}$ is
regular with respect to the Dirichlet problem and $\mbox{cap}(e)
=0$, we have that $g_{\Omega}(z;\infty) =
g_{\widetilde{\Omega}}(z;\infty), \widetilde{\Omega} =
\overline{\mathbb{C}} \setminus \widetilde{E}$, and
$g_{\widetilde{\Omega}}(z;\infty)$ extends continuously to all
$\mathbb{C}$.

Fix $\varepsilon > 0$ and let $U_{\varepsilon} = \{z \in
{\mathbb{C}}: g_{\widetilde{\Omega}}(z;\infty) < \varepsilon\}$.
This is an open set which contains $\widetilde{E}$ where
$g_{\widetilde{\Omega}}(z;\infty) =0$ identically. Since the set $e$
is at most denumerable and all its accumulation points are contained
in $\widetilde{E}$, it follows that $e \setminus U_{\varepsilon}$
has at most a finite number of points (or may be empty). Let
$\{z_1,\ldots,z_N\}$ be the set of such points (should there be any). For each fixed $k =
1,\ldots,N,$
\[ \frac{|p_l(z_k)|^2}{\|p_l\|_2^2}\sigma(z_k) \leq \int
\frac{|p_l(x)|^2}{\|p_l\|_2^2} d\sigma(x) =  1.
\]
Consequently,
\[ \limsup_{l \in \Lambda}
\left(\frac{|p_l(z_k)|}{\|p_l\|_2}\right)^{1/l} \leq 1, \qquad
k=1,\ldots,N.
\]
since $\sigma(z_k) > 0$. Because of (\ref{eq:25})
\[ \limsup_{l \in \Lambda}
\left(\frac{\|p_l\|_{E \cap
\overline{U}_{\varepsilon}}}{\|p_l\|_2}\right)^{1/l} \leq
\exp({\varepsilon})
\]
This, together with the previous inequality for $z_k,
k=1,\ldots,N$, immediately imply that
\[ \limsup_{l \in \Lambda}
\left(\frac{\|p_l\|_{E}}{\|p_l\|_2}\right)^{1/l} \leq
\exp({\varepsilon}).
\]
The arbitrariness of $\varepsilon > 0$ renders what we set out
to prove. \hfill $\Box$ \vspace{0,2cm}

The assumption that the points in $e$ only accumulate on
$\widetilde{E}$ is essential. If this was not the case one can
construct examples where $(\ref{eq:23})$ does not hold.

\vspace{0,2cm}

{\bf Proof of Theorem \ref{teo:5}.} We will prove $|{\bf n}_1|$-th root asymptotics for the sequence $\{{\mathcal{B}}_{\bf n}\}, {\bf n} \in \Lambda.$ Since $\widehat{s}_{1,\lambda_1(0)}(z) \neq 0, z \in {{\mathbb{C}}}\setminus \Delta^1_1,$ the statement of the theorem readily follows with the same limit.

For definiteness, in reordering the components of a given ${\bf n}$, let us take that unique pair of permutations $(\lambda_1,\lambda_2)$ such that for each $i=1,2,$ whenever $n_{i,\lambda_i(j)} = n_{i,\lambda_i(k)}$ for some $0 \leq j < k \leq m_i,$ then $\lambda_i(j) < \lambda_i(k).$ By $\Lambda(\lambda_1,\lambda_2)$, we denote the set of all multi-indices in $\Lambda$ whose components ${\bf n}_1,{\bf n}_2,$ are reordered decreasingly with $\lambda_1$ and $\lambda_2$ respectively. We are only interested in those $\Lambda(\lambda_1,\lambda_2)$ containing an infinite number of elements of $\Lambda.$ Fix $(\lambda_1,\lambda_2)$ and let
\[ (\widetilde{S}^1,\widetilde{S}^2) = ({\mathcal{N}}(\widehat{s}_{1,\lambda_1(0)}
\rho_0^2,\rho_{1}^1,\ldots,\rho_{m_1}^1),{\mathcal{N}}(\widehat{s}_{1,\lambda_1(0)}
\rho_0^2,\rho_1^2,\ldots,\rho_{m_2}^2))
\]
be the pair of Nikishin systems associated with ${\mathcal{A}}_{\bf n}$ by Theorems \ref{teo:3}-\ref{teo:4} with respect to which ${\mathcal{B}}_{\bf n}$ is a multiple orthogonal linear form.

Set
\[
P_j=\sum_{k=j}^{m_1} p_{1,\lambda_1(k)},\,\, j=0,\ldots,m_1, \qquad P_{-j} =
\sum_{k=j}^{m_2} p_{2,\lambda_2(k)},\,\, j=0,\ldots,m_2.
\]
Define the  tri-diagonal matrix
\begin{equation}\label{matriz}
{\mathcal{C}}=
\begin{pmatrix}
P_{-m_2}^2 & -\frac{P_{-m_2}P_{-m_2+1}}{2} & 0 &  \cdots  & 0\\
-\frac{P_{-m_2}P_{-m_2+1}}{2} & P_{-m_2+1}^2 &
-\frac{P_{-m_2+1}P_{-m_2+2}}{2} & \cdots  & 0\\
0 & -\frac{P_{-m_2+1}P_{-m_2+2}}{2} & P_{-m_2+2}^2
& \cdots  & 0 \\
\vdots &\vdots&\vdots&\ddots &\vdots\\
0 & 0& 0 &\cdots  & P_{m_1}^2
\end{pmatrix}.
\end{equation}
The sub-indices of the entries $c_{j,k}$ of $\mathcal{C}$ run from $-m_2 -1$ to $m_1+1$.

Let ${\mathcal{M}}_1(E_k) $ be the subclass of probability measures of
${\mathcal{M}}(E_k),$
\[ E_k = \left\{
\begin{array}{ll}
\supp \rho_k^1, & k=1,\ldots,m_1, \\
\supp \rho_{-k}^2, & k=-m_2,\ldots,0.
\end{array}
\right.
\]
Denote
\[{\mathcal{M}}_1= {\mathcal{M}}_1(E_{-m_2})
\times \cdots \times {\mathcal{M}}_1(E_{m_1})  \,.
\]
Given a vector measure $\mu=(\mu_{-m_2},\ldots,\,\mu_{m_1}) \in
{\mathcal{M}}_1$ and $j \in \{ -m_2,\ldots,m_1\},$ we define the combined
potential
\[W^{\mu}_j(x) = \sum_{k=-m_2}^{m_1} c_{j,k}
V^{\mu_k}(x) , \qquad  V^{\mu_k}(x) = \int \log \frac{1}{|x-y|} \,d\mu_k(y).
\]
Set
\[ J(\mu) = \sum_{k,j = -m_2}^{m_1} c_{j,k} \int \int \log \frac{1}{|x-y|} d\mu_j(x) d\mu_k(y) =     \int W^{\mu}_j(x) d\mu_j(x)  .
\]

From Propositions 4.1-4.5 in \cite[Chapter 5]{NiSo} it follows that there exists a unique vector measure
$\overline{\mu} = (\overline{\mu}_{-m_2},\ldots,\overline{\mu}_{m_1}) \in
{\mathcal{M}}_1$ such that
\begin{equation} \label{eq:minima}
  J(\overline{\mu}) = \inf \{J(\mu): \mu \in {\mathcal{M}}_1\},
\end{equation}
and that exist constants $w_j^{\overline{\mu}}, j=-m_2,\ldots,m_1,$ for which
\begin{equation}
\label{eq:equilibrio}
W_j^{\overline{\mu}}(x) \left\{
\begin{array}{ll}
\leq w_j^{\overline{\mu}}, & x \in \supp \overline{\mu}_j, \\
\geq w_j^{\overline{\mu}}, & x \in E_j \setminus  A_j,\, \mbox{\rm
cap} (A_j) = 0.
\end{array}
\right.
\end{equation}
for certain Borel sets $A_j$. For any two vector measures $ {\mu}^1,
{\mu}^2 \in {\mathcal{M}}_1$ such that $J( {\mu}^1) < \infty, J(
{\mu}^2) < \infty,$ straightforward calculations yield
\[ J(\mu^2) - J(\mu^1) = J(\mu^2 - \mu^1) + 2 \sum_{j=-m_2}^{m_1} \int W_j^{\mu^1} (x) d(\mu_j^2 - \mu_j^1)(x).
\]
Since $J(\mu^2 - \mu^1) \geq 0$ for all $\mu^1,\mu^2 \in {\mathcal{M}}_1$ (see \cite[Proposition 4.2]{NiSo}), and sets of capacity zero are negligible for measures with finite energy, if $\mu^1$ satisfies (\ref{eq:equilibrio}) it also satisfies (\ref{eq:minima}). Thus, (\ref{eq:minima})-(\ref{eq:equilibrio}) are equivalent and a measure verifying any one of the two is unique (and so are the constants in (\ref{eq:equilibrio})).
$\overline{\mu}$ is called the equilibrium vector measure, and $w^{\overline{\mu}} = (w_{-m_2}^{\overline{\mu}},\ldots,w_{m_1}^{\overline{\mu}})$ the equilibrium vector constant, for the logarithmic potential governed by the interaction matrix ${\mathcal{C}}$ on the system of compact sets $E_j, j=-m_2,\ldots,m_1.$

From Lemma \ref{lem:cond} we have that $(\widetilde{S}^1,\widetilde{S}^2) \in \mbox{\bf Reg}$ and the supports of the generating measures are quasi-regular. If ${\mathcal{A}}_{\bf n}$ is monic (see Definition \ref{monico}), due to the way in which ${\mathcal{B}}_{\bf n}$ is constructed (in particular, see (\ref{eq:M}) in the proof of Lemma $\ref{lem:4}$ and the proof of Theorem \ref{teo:4}) it follows that $b_{{\bf n},{m_1}}$ is either plus or minus $ a_{{\bf n},{\lambda_1^{-1}(m_1)}}$. Thus, its leading coefficient is either $1$ or $-1$; that is, except for a sign change, ${\mathcal{B}}_{\bf n}$ is monic with the normalization imposed in \cite[Theorem 5.1]{FLLS}. Following the proof of \cite[Theorem 5.1]{FLLS}, but using Lemma \ref{lem:asint} instead of \cite[Lemma 5.1]{FLLS}, one finds that
\[
\lim_{{\bf n}\in \Lambda(\lambda_1,\lambda_2)} |{\mathcal{B}}_{\bf n}(z)|^{1/|{\bf n}_1} = \exp\left(P_{1} V^{\overline{\mu}_{1}}(z)-P_0
V^{\overline{\mu}_0}(z) - 2\sum_{k=1}^{m_1}
\frac{\omega_k^{\overline{\mu}}}{P_k}\right), \qquad {\mathcal{K}} \subset {\mathbb{C}} \setminus (\Delta_0^1 \cup \Delta_1^1),
\]
where $\overline{\mu}= \overline{\mu}({\mathcal{C}})=
(\overline{\mu}_{-m_2},\ldots,\overline{\mu}_{m_1})$ is  the
equilibrium vector measure   and
$(\omega_{-m_2}^{\overline{\mu}},\ldots,\omega_{m_1}^{\overline{\mu}})$
is the system of equilibrium constants for the vector potential
problem determined by the interaction matrix $\mathcal{C}$ defined
in $(\ref{matriz})$  on the system of compact sets $E_j,  j=-m_2,\ldots,m_1.$

It is easy to see that the interaction matrix ${\mathcal{C}}$ does not depend on $(\lambda_1,\lambda_2)$ and
that the compact sets $E_j, j=-m_2,\ldots,m_1,$ for each fixed $j$,  may differ only on a (denumerable) set of capacity zero depending on $\lambda_1,\lambda_2$ (see proof of Lemma \ref{lem:cond}). Therefore, the equilibrium measure and the equilibrium constant are uniquely determined for any $\Lambda(\lambda_1,\lambda_2)$ containing infinitely many terms of $\Lambda$.  Consequently,
\begin{equation}
 \label{eq:lim5}
 \lim_{{\bf n}\in \Lambda} |{\mathcal{A}}_{\bf n}(z)|^{1/|{\bf n}_1} =  \exp\left(P_{1} V^{\overline{\mu}_{1}}(z)-P_0
V^{\overline{\mu}_0}(z) - 2\sum_{k=1}^{m_1}
\frac{\omega_k^{\overline{\mu}}}{P_k}\right), \quad {\mathcal{K}} \subset {\mathbb{C}} \setminus (\Delta_0^1 \cup \Delta_1^1).
\end{equation}
With this we conclude the proof. \hfill $\Box$

\begin{rem}
\label{debil}
If we denote by $Q_{{\bf n},0}$ the monic polynomial whose zeros are those of ${\mathcal{A}}_{\bf n}$, under the assumptions of Theorem \ref{teo:5}, we have (see \cite[Theorem 4.2]{FLLS})
\[
*\lim_{{\bf n} \in \Lambda} \mu_{Q_{{\bf n},0}} = \overline{\mu}_0.
\]
There are other linear forms related with ${\mathcal{A}}_{\bf n}$ whose asymptotic zero distribution and logarithmic asymptotic is described in terms of the other components of $\overline{\mu}$ and the vector equilibrium constant. This allows to give the logarithmic asymptotics of the polynomials $a_{{\bf n},k}$ as well. For a view of what can be expected, see \cite[Section 5]{FLLS}. These results can be used to give the exact rate of convergence of mixed type Hermite Pad\'e
approximants, see \cite[Section 7]{FLLS}  and \cite[Theorem 7]{LF2}. For example, in case of type II approximation, under regularity of the generating measures and quasi--regularity of their supports, the following limit exists
\[ \lim_{{\bf n}\in \Lambda} \|\widehat{s}_{0,k} - \frac{P_{{\bf n},k}}{Q_{\bf n}}\|_{\mathcal{K}}^{1/2|{\bf n}|}, \quad {\mathcal{K}} \subset \overline{\mathbb{C}} \setminus  (\mbox{\rm Co}(\supp \sigma_0) \cup \mbox{\rm Co}(\supp \sigma_1)), \quad k=0,\ldots,m.
\]
\end{rem}
\vspace{0,2cm}

{\bf Proof of Theorem \ref{teo:6}.} The existence of the limit claimed in Theorem \ref{teo:6} follows directly from \cite[Theorem 1.4]{FLLS}, but to give an expression of the limit function, we must introduce some notions.

Let $\lambda_1,\lambda_2,$ and $l = (l_1;l_2)$ be as given in Theorem \ref{teo:6}.   Consider the $(m_1+m_2+2)$-sheeted Riemann surface
$$
{\mathcal{R}}=\overline{\bigcup_{k=-m_2-1}^{m_1} {\mathcal{R}}_k} ,
$$
formed by the consecutively ``glued'' sheets
$$
{\mathcal{R}}_{-m_2-1}:=\overline {\mathbb{C}} \setminus
\widetilde{\Delta}_{-m_2},\,\,\, {\mathcal{R}}_k:=\overline
{\mathbb{C}} \setminus (\widetilde{\Delta}_k \cup
\widetilde{\Delta}_{k+1}),\,\, k=-m_2,\ldots,m_1-1,\,\,\, {\mathcal{R}}_{m_1}:=\overline {\mathbb{C}} \setminus
\widetilde{\Delta}_{m_1},
$$
where the upper and lower banks of the slits of two neighboring
sheets are identified. Define
\[(\widetilde{l}_1;\widetilde{l_2}) := (\lambda_1^{-1}(l_1);\lambda_2^{-1}(l_2)).
  \]
Let $\psi^{(\widetilde{l})}$ be a singled valued function
defined on $\mathcal{R}$ onto the extended complex plane
satisfying
\[ \psi^{(\widetilde{l})}(z)=\frac{C_{1}}{z}+{\mathcal{O}}(\frac{1}{z^{2}}),\quad
z\rightarrow\infty ^{(-\widetilde{l}_2-1)},\]
\[ \psi^{(\widetilde{l})}(z)=C_{2}\,z+{\mathcal{O}}(1),\quad z\rightarrow\infty ^{(\widetilde{l}_1)},\]
where $C_{1}$ and $C_{2}$ are nonzero constants. Since the genus
of $\mathcal{R}$ is zero, $\psi^{(\widetilde{l})}$ exists and is uniquely
determined up to a multiplicative constant. Consider the branches
of $\psi^{(\widetilde{l})}$, corresponding to the different sheets
$k=-m_2-1,\ldots,m_1$ of $\mathcal{R}$
\[\psi^{(\widetilde{l})}:=\{\psi^{(\widetilde{l})}_{k}\}_{k=-m_2-1}^{m_1}\,. \]
Given an arbitrary function $F(z)$  which has in a neighborhood of
infinity a Laurent expansion of the form $F(z)= Cz^k +
{\mathcal{O}}(z^{k-1}), C \neq 0,$ and $ k \in {\mathbb{Z}},$ we
denote
\[
\widetilde{F}:= {F}/{C}\,.
\]

Because of Theorem \ref{teo:4}, Lemma \ref{lem:cond}, and the normalization adopted, the sequence $\{{\mathcal{B}_{\bf n}}\}, n \in \Lambda,$ satisfies all the assumptions of \cite[Theorem 6.8]{FLLS}. Consequently,
\[
\lim_{{\bf n}\in\Lambda}\,\frac{{\mathcal{B}}_{{\bf
n}^l}(z)}{{\mathcal{B}}_{{\bf n}}(z)}=
C(\widetilde{l})\widetilde{\psi}^{(\widetilde{l})}_{0}(z),\qquad
{\mathcal{K}} \subset{\mathbb{C}}\setminus(\supp{\rho_0^1} \cup
\supp{\rho_1^1}),
\]
where  $ C(\widetilde{l})$ is a constant, which only depends on
$\widetilde{l}$ and can be determined exactly (see (69), (73), and
(83) in \cite{FLLS}) in terms of  the values of the branches of
$\psi^{(\widetilde{l})}$ at $\infty$. Due to the relation between
${\mathcal{B}}_{{\bf n}}$ and ${\mathcal{A}}_{{\bf n}}$, we obtain
\begin{equation}\label{eq:lim6}
\lim_{{\bf n}\in\Lambda}\,\frac{{\mathcal{A}}_{{\bf n}^l}(z)}{{\mathcal{A}}_{{\bf
n}}(z)}= C(\widetilde{l})\widetilde{\psi}^{(\widetilde{l})}_{0}(z),\qquad  {\mathcal{K}}
\subset{\mathbb{C}}\setminus(\supp{\sigma_0^1} \cup \mbox{Co}(\supp{\sigma_1^1})),
\end{equation}
since $\supp{\sigma_0^1} = \supp{\rho_0^1}$ and $\mbox{Co}(\supp{\rho_1^1}) \subset \mbox{Co}(\supp{\sigma_1^1})$.
\hfill $\Box$ \vspace{0,2cm}

{\bf Proof of Corollary \ref{cor:3}.} As in the proof of Theorem
\ref{teo:5}, for definiteness, in reordering the components of a
given ${\bf n}$, let us take that unique pair of permutations
$(\lambda_1,\lambda_2)$ such that for each $i=1,2,$ whenever
$n_{i,\lambda_i(j)} = n_{i,\lambda_i(k)}$ for some $0 \leq j < k
\leq m_i,$ then $\lambda_i(j) < \lambda_i(k).$ By
$\Lambda(\lambda_1,\lambda_2)$, we denote the set of all
multi-indices in $\Lambda$ whose components ${\bf n}_1,{\bf n}_2,$
are reordered decreasingly with $\lambda_1$ and $\lambda_2$
respectively. We are only interested in those
$\Lambda(\lambda_1,\lambda_2)$ containing an infinite number of
elements of $\Lambda.$ Fix $(\lambda_1,\lambda_2)$ and let
\[ (\widetilde{S}^1,\widetilde{S}^2) = ({\mathcal{N}}(\widehat{s}_{1,\lambda_1(0)}
\rho_0^2,\rho_{1}^1,\ldots,\rho_{m_1}^1),{\mathcal{N}}(\widehat{s}_{1,\lambda_1(0)}
\rho_0^2,\rho_1^2,\ldots,\rho_{m_2}^2))
\]
be the pair of Nikishin systems associated with ${\mathcal{A}}_{\bf n}$ by Theorems \ref{teo:3}-\ref{teo:4} with respect to which ${\mathcal{B}}_{\bf n}$ is a multiple orthogonal linear form.

Let $M$ be the least common multiple of
$m_1+1$ and $m_2+1$, and define
$d_{1}:=M/(m_{1}+1)$, $d_{2}:=M/(m_{2}+1)$. Within the class of
pairs $l=(l_{1};l_{2})$ with $0\leq l_{1}\leq m_{1}$, $0\leq
l_{2}\leq m_{2}$, we distinguish the subclass
\[L:=\{(l_{1};l_{2}): l_{1}\equiv r\,\text{mod}\,(m_{1}+1),\,\,
l_{2}\equiv r\,\text{mod}\,(m_{2}+1)\,\,\mbox{for some}\,\,0\leq
r\leq M-1\}\,.
\]
It is easy to check that for different $r, 0\leq r\leq M-1,$ the
pairs $(l_1,l_2)$ in $L$ are distinct. Let
$\mathbf{p}:=(\mathbf{p}_{1};\mathbf{p}_{2})$, where
${\mathbf{p}}_{1}=(d_{1},\ldots,d_{1})$ and
${\mathbf{p}}_{2}=(d_{2},\ldots,d_{2})$ have $m_{1}+1$ and $m_{2}+1$
components, respectively. By ${\bf n}+{\bf p}$ we denote the
multi-index
$(\mathbf{n}_{1}+\mathbf{p}_{1};\mathbf{n}_{2}+\mathbf{p}_{2})$; that is, ${\bf n}+{\bf p} = \widetilde{\bf n}$.

Given ${\bf n}\in\Lambda(\lambda_1,\lambda_2)$ and $0\leq r\leq M$, let
${\mathbf{n}}(r):={\mathbf{n}}+{\mathbf{q}}(r)$ where
${\mathbf{q}}(r)=({\mathbf{q}}_{1}(r);{\mathbf{q}}_{2}(r))$ is the
multi-index satisfying $({\mathbf{q}}_i(r) = (q_{i,0}(r),\ldots,q_{i,m_i}(r)), i=1,2)$
\[
q_{i,\lambda_i(j)}(r)=
\left\{
\begin{array}{ll}
k_i +1, & j = 0,\ldots,s_i-1, \\
k_i,    & j = s_i,\ldots,m_i,
\end{array}
\right.
\quad
r=k_i(m_{i}+1)+s_i, \quad 0\leq s_i\leq m_{i}\,.
\]
Hence, ${\mathbf{n}}(0)={\mathbf{n}}$, ${\mathbf{n}}(M)={\bf n}+{\bf p} = \widetilde{\bf n}$. It is easy to see that for all $r \in \{0,\ldots,M-1\}$, the same pair $(\lambda_1,\lambda_2)$ reorders the components of ${\mathbf{n}}(r)$
giving rise to the same systems $(\widetilde{S}^1,\widetilde{S}^2)$.

We have
\[
\frac{{\mathcal{A}}_{{\bf n}+{\bf p}}(z)}{{\mathcal{A}}_{{\bf
n}}(z)}=\prod_{r=0}^{M-1}\frac{{\mathcal{A}}_{{\bf n}(r+1)}(z)}{{\mathcal{A}}_{{\bf
n}(r)}(z)}\,.
\]
Due to (\ref{eq:lim6})
\begin{equation}
\label{eq:51}
\lim_{{\bf n}\in\Lambda(\lambda_1,\lambda_2)}
\frac{{\mathcal{A}}_{\widetilde{\bf n}}(z)}{{\mathcal{A}}_{{\bf
n}}(z)}= \prod_{(l_1,l_2) \in L}
C(\widetilde{l})\widetilde{\psi}^{(\widetilde{l})}_{0}(z), \qquad {\mathcal{K}}
\subset{\mathbb{C}}\setminus(\supp{\sigma_0^1} \cup
\mbox{Co}(\supp{\sigma_1^1})),
\end{equation}
where $l=(l_{1};l_{2})$ is precisely the multi-index satisfying
$l_{1}\equiv r\,\text{mod}\,(m_{1}+1)$, $l_{2}\equiv
r\,\text{mod}\,(m_{2}+1)$, and $\widetilde{l} = (\widetilde{l}_1;\widetilde{l}_2) = (\lambda_1^{-1}(l_1);\lambda_2^{-1}(l_2))$.  The limit does not depend on $(\lambda_1,\lambda_2)$ because the set $\widetilde{L} = \{(\widetilde{l}_1;\widetilde{l}_2) : (l_{1};l_{2}) \in L\}$ is the same for all $(\lambda_{1},\lambda_{2})$. The proof is complete. \hfill $\Box$

\begin{rem}
\label{final} The linear forms associated with ${\mathcal{A}}_{\bf n}$ mentioned in the previous remark also satisfy ratio asymptotics in the spirit of the results contained in \cite[Section 6]{FLLS}.
\end{rem}

\end{document}